\documentclass[11pt]{article}

\usepackage{amssymb,amsmath,amscd,amsthm,dsfont,xy,enumitem,mathrsfs,pstricks}
\xyoption{all}

\numberwithin{equation}{section} 

\newtheorem{thm}{Theorem}[section]
\newtheorem{prp}[thm]{Proposition}
\newtheorem{lmm}[thm]{Lemma}

\theoremstyle{definition}

\theoremstyle{remark}
\newtheorem{rmk}[thm]{Remark}

\def\lra{\longrightarrow}

\def\Llra{\Longleftrightarrow}

\def\xlra#1{\xrightarrow{#1}}

\def\BE#1{\begin{equation}\label{#1}}
\def\EE{\end{equation}}
\def\eref#1{(\ref{#1})}

\def\blr#1{\big\langle#1\big\rangle}
\def\wt#1{\widetilde{#1}}
\def\ov#1{\overline{#1}}
\def\tn#1{\textnormal{#1}}
\def\sf#1{\textsf{#1}}

\def\sm#1{\begin{small}#1\end{small}}

\def\cA{\mathcal A}
\def\C{\mathbb C}
\def\cC{\mathcal C}
\def\bE{\mathbb E}

\def\cI{\mathcal I}
\def\cM{\mathcal M}
\def\cN{\mathcal N}
\def\P{\mathbb P}
\def\cP{\mathcal P}
\def\cR{\mathcal R}
\def\Q{\mathbb Q}
\def\R{\mathbb R}

\def\Z{\mathbb Z}

\def\fc{\mathfrak c}

\def\ff{\mathfrak f}

\def\fo{\mathfrak o}

\def\de{\delta}
\def\ep{\epsilon}
\def\ga{\gamma}
\def\io{\iota}
\def\ka{\kappa}

\def\si{\sigma}

\def\vr{\varrho}

\def\Om{\Omega}
\def\Si{\Sigma}

\def\BL{\mathbb{BL}}
\def\BLau{\BL^{\!\!\fc_1}}

\def\BLaua{\BL^{\!\!1}}
\def\BLC{\BL^{\!\!\C}}
\def\BLR{\BL^{\!\!\R}}

\def\cha{\tn{char}}

\def\CR{\tn{CR}}
\def\nd{\tn{d}}

\def\id{\tn{id}}
\def\Id{\tn{Id}}
\def\Im{\tn{Im}}

\def\nod{\tn{nd}}

\def\PD{\tn{PD}}
\def\pt{\tn{pt}}
\def\PS{\tn{PS}}
\def\rk{\tn{rk}}

\def\prl{\parallel}
\def\prt{\partial}
\def\eset{\emptyset}
\def\i{\infty}
\def\bu{\bullet}
\def\1{\mathds{1}}

\begin{document}

\title{The Cohomology Ring of the Deligne-Mumford Moduli Space\\
of Real Rational Curves with Conjugate Marked Points} 
\author{Xujia Chen\thanks{Partially supported by NSF grant DMS 1901979 and Simons Foundation } ,
Penka Georgieva\thanks{Partially supported by 
ANR grant 18-CE40-0009 and ERC Grant ROGW-864919} ,$~$and 
Aleksey Zinger\thanks{Partially supported by NSF grant DMS 1901979 and Simons Foundation}}
\date{\today}
\maketitle

\begin{abstract}
\noindent   
It is a long-established and heavily-used fact that the integral cohomology ring 
of the Deligne-Mumford moduli space of (complex) rational curves 
is the polynomial ring on the boundary divisors modulo the ideal generated by 
the obvious geometric relations between them.
We show that the rational cohomology ring of the Deligne-Mumford moduli space 
of real rational curves with conjugate marked points only 
is the polynomial ring on certain (``complex") boundary divisors
and real boundary hypersurfaces modulo the ideal generated by 
the obvious geometric relations between them and the geometric relation 
in positive dimension and codimension identified in a previous paper.
The proof of the generation result also implies that the integer homology
has only 2-torsion.
\end{abstract}

\tableofcontents

\section{Introduction}
\label{intro_sec}

\noindent
The Deligne-Mumford moduli space~$\ov\cM_{\ell}$ of (complex) rational curves 
with $\ell\!\ge\!3$ marked points is a (complex) blowup of~$(\C\P^1)^{\ell-3}$.
Thus, the $\Z$-(co)homology of~$\ov\cM_{\ell}$ is even-dimensional and torsion-free. 
The smooth complex projective variety~$\ov\cM_{\ell}$ contains  
smooth divisors~$D_{J,K}$ whose generic element is a two-component curve.
These are known as the \sf{boundary divisors} of~$\ov\cM_{\ell}$;
the complement of their union is the open subspace $\cM_{\ell}\!\subset\!\ov\cM_{\ell}$
consisting of smooth curves.
It is a now classical result of Keel~\cite{Keel} that 
\begin{enumerate}[label=(CDM\arabic*),leftmargin=*]

\item\label{Cgen_it}  these divisors generate  $H^*(\ov\cM_{\ell};\Z)$,

\item\label{Crel_it} subject only to some obvious relations; 

\end{enumerate}
see Theorem~\ref{HC_thm}.
The same statements apply with $H^*(\ov\cM_{\ell};\Z)$ replaced by
the Chow ring $A^*(\ov\cM_{\ell})$ of~$\ov\cM_{\ell}$.\\

\noindent
The Deligne-Mumford moduli space~$\R\ov\cM_{0,\ell}$ of real rational curves 
with $\ell\!\ge\!2$ conjugate pairs of marked points is a compact orientable manifold 
of (real) dimension~$2\ell\!-\!1$.
As shown in Section~3 of~\cite{RealEnum}, $H_1(\R\ov\cM_{0,3};\Z)$ contains 2-torsion;
this implies that $H_1(\R\ov\cM_{0,\ell};\Z)$ contains 2-torsion for every $\ell\!\ge\!3$.
The smooth manifold~$\R\ov\cM_{0,\ell}$ contains compact 
real hypersurfaces~$\R E_{J,K}$ and~$\R H_{J,K}$ whose generic elements are 
two-component curves.
It also contains compact submanifolds~$\R D_{I;J,K}$
of (real) codimension~2 whose generic elements are three-component curves.
We call these codimension~1 and~2 submanifolds \sf{boundary hypersurfaces} 
and \sf{boundary divisors}, respectively;
the complement of their union is the open subspace 
$\R\cM_{0,\ell}\!\subset\!\R\ov\cM_{0,\ell}$ consisting of smooth curves.
We show in this paper that 
\begin{enumerate}[label=(RDM\arabic*),leftmargin=*]

\item\label{Rgen_it}  
the submanifolds~$\R E_{J,K}$ and~$\R D_{I;J,K}$ generate $H^*(\R\ov\cM_{0,\ell};\Q)$,

\item\label{Rrel_it}  subject only to some obvious relations and the relation 
of Remark~3.4 in~\cite{RealEnum};

\end{enumerate}
see Theorem~\ref{HR_thm}.
This result is complementary to~\cite{EHKR},
which determines the cohomology of the Deligne-Mumford moduli space of 
stable real rational curves with real marked points only.
Just as Getzler's relation~\cite{Getz} on the Deligne-Mumford moduli space of
genus~1 curves with 4~marked points lies in the middle (co)homology of
this space,
the relation  on~$\R\ov\cM_{0,3}$ established in Remark~3.4 in~\cite{RealEnum}
lies in the first homology (or second cohomology) of
three-dimensional orientable manifold~$\R\ov\cM_{0,3}$.\\

\noindent
Similarly to~\cite{Keel}, our proof of Theorem~\ref{HR_thm} involves
induction on the number of marked points and a sequence of blowups~$X_{\vr}$ 
starting with~$\R\ov\cM_{0,\ell}\!\times\!\C\P^1$ and ending
with~$\R\ov\cM_{0,\ell+1}$.
As in~\cite{Keel}, we use this sequence to show by induction that certain submanifolds 
algebraically generate the cohomologies of the blowups~$X_{\vr}$ and 
thus of the moduli spaces.
However, there are a number of significant differences  with the proof in~\cite{Keel}.
For example, three types of blowups, real, complex, and ``augmented", 
appear in our sequence, as opposed to only one in~\cite{Keel}.
In contrast to~\cite{Keel}, we do not inductively determine the cohomology ring of 
each blowup~$X_{\vr}$, as the description of the relations among the distinguished generators
is already cumbersome in the complex case and would be much more so in the real case.
We instead rely on Propositions~\ref{Rhomol_prp} and~\ref{Ralg_prp}
to inductively determine all relations between
the generators in the cohomology of the top blowup~$\R\ov\cM_{0,\ell+1}$  only.
These propositions are the real analogues of Propositions~\ref{Chomol_prp} and~\ref{Calg_prp}.
Proposition~\ref{Chomol_prp} is equivalent to a homology version of~(3) in \cite[p549]{Keel}
with a typo corrected, but 
is deduced from Theorem~\ref{HC_thm} there (as opposed to the other way around,
as done in the present paper).\\

\noindent
Propositions~\ref{Chomol_prp} and~\ref{Rhomol_prp} readily lead to recursive formulas
for the Betti numbers
$$h_{\ell}^p\equiv\rk_{\Z}H_p\big(\ov\cM_{\ell};\Z\big) \quad\hbox{and}\quad
h_{0,\ell}^p\equiv\dim_{\Q}H_p\big(\R\ov\cM_{0,\ell};\Q\big) $$
of the moduli spaces~$\ov\cM_{\ell}$ and~$\R\ov\cM_{0,\ell}$, respectively.
By $\ov\cM_3\!=\!\{\pt\}$ and Proposition~\ref{Chomol_prp},
\BE{ChomolRec_e}h_3^p=\begin{cases}1,&\hbox{if}~p\!=\!0;\\
0,&\hbox{if}~p\!\neq\!0;\end{cases} \qquad
h_{\ell+1}^p=h_{\ell}^p\!+\!h_{\ell}^{p-2}\!+\!
\frac12\sum_{j=2}^{\ell-2}\binom{\ell}{j}\sum_{q=0}^{p-2}h_{j+1}^qh_{\ell-j+1}^{p-2-q}
~~\forall\,\ell\!\ge\!3\,;\EE
this is equivalent to~(4) in \cite[p550]{Keel} with typos corrected.
By $\R\ov\cM_{0,2}\!=\!\R\P^1$ and Proposition~\ref{Rhomol_prp},
\begin{gather}\notag
h_{0,2}^p=\begin{cases}1,&\hbox{if}~p\!=\!0,1;\\
0,&\hbox{if}~p\!\neq\!0,1;\end{cases}\\
\label{RhomolRec_e} 
h_{0,\ell+1}^p=h_{0,\ell}^p\!+\!h_{0,\ell}^{p-2}\!+\!
2^{\ell-1}\big(h_{\ell+1}^{p-1}\!+\!h_{\ell+1}^{p-2}\big)
+\sum_{i=1}^{\ell-2}2^{\ell-i}\binom{\ell}{i}\sum_{q=0}^{p-2}h_{0,i+1}^qh_{\ell-i+1}^{p-2-q}
~~\forall\,\ell\!\ge\!2\,.
\end{gather}
In particular, $h_{0,\ell}^1\!=\!2^{\ell-1}\!-\!1$.
The Betti numbers of~$\ov\cM_{\ell+1}$ and~$\R\ov\cM_{0,\ell}$
for small values of~$\ell$ are shown in Table~\ref{Betti_tbl}.\\

\begin{table}
\begin{small}
\begin{center}
\begin{tabular}{||c||c|c|c|c|c||}
\hline\hline
$\ell$& 2& 3& 4& 5&6\\
\hline
$\ov\cM_{\ell+1}$& 1& 1,1& 1,5,1&1,16,16,1& 1,42,127,42,1\\
\hline
\!$\R\ov\cM_{0,\ell}$\!&\!1,1\!&\!1,3,3,1\!&\!1,7,20,20,7,1\!&\!1,15,85,171,171,85,15,1\!
&\!1,31,302,1042,2032,2032,1042,302,31,1\!\\
\hline\hline
\end{tabular}
\end{center}
\end{small}
\vspace{-.2in}
\caption{The Betti numbers of $\ov\cM_{\ell+1}$ and $\R\ov\cM_{0,\ell}$.}
\label{Betti_tbl}
\end{table}

\noindent
For illustrative purposes, we demonstrate the approach of this paper on
the complex case treated in~\cite{Keel}.
After recalling Keel's result in Theorem~\ref{HC_thm} in Section~\ref{Cthm_subs},
we formulate its real analogue in Theorem~\ref{HR_thm} in Section~\ref{Rthm_subs}.
Based on geometric considerations stated in these sections,
the homomorphisms~\eref{HCthm_e} and~\eref{HRthm_e} of Theorems~\ref{HC_thm} and~\ref{HR_thm}
are well-defined.
In Section~\ref{CohGener_sec}, we first summarize the relevant topological properties
of the three types of blowup that appear in this paper
and then deduce the surjectivity of the homomorphisms~\eref{HCthm_e} and~\eref{HRthm_e}
from these properties and the sequences of blowups constructed in~\cite{RDMbl}.
In Section~\ref{CohRelat_sec}, we use the generation conclusions of Section~\ref{CohGener_sec}
to inductively describe  the rings~$H^*(\ov\cM_{\ell};\Z)$ and~$H^*(\R\ov\cM_{0,\ell};\Q)$
in Propositions~\ref{Chomol_prp} and~\ref{Rhomol_prp}.
Propositions~\ref{Calg_prp} and~\ref{Ralg_prp} give partial inductive descriptions
of the left-hand sides of the homomorphisms~\eref{HCthm_e} and~\eref{HRthm_e}. 
These four propositions imply the injectivity of the homomorphisms~\eref{HCthm_e} and~\eref{HRthm_e}
and conclude the proofs of Theorems~\ref{HC_thm} and~\ref{HR_thm}.
The approach of this paper should be extendable to the integer cohomology
of~$\R\ov\cM_{0,\ell}$ and to 
the real genus~0 moduli spaces~$\R\ov\cM_{k,\ell}$ with real and conjugate points.


\section{The cohomology of Deligne-Mumford spaces}
\label{thm_sec}

\noindent
For $\ell\!\in\!\Z$, define
\begin{gather*}
[\ell]=\big\{i\!\in\!\Z^+\!\!:1\!\le\!i\!\le\!\ell\big\}, \qquad
[\ell^{\pm}]=\big\{1^+,1^-,\ldots,\ell^+,\ell^-\big\}, \\
\begin{aligned}
\cP(\ell)&=\big\{\!\{J,K\}\!\!:[\ell]\!=\!J\!\sqcup\!K\big\}, 
&\cP_{\bu}(\ell)&=\big\{\!\{J,K\}\!\in\!\cP(\ell)\!\!:2\!\le\!|J|\!\le\!\ell\!-\!2\big\},\\
\wt\cP(\ell)&=\big\{\!(I,\{J,K\})\!\!:[\ell]\!=\!I\!\sqcup\!J\!\sqcup\!K\big\},
&\wt\cP_{\bu}(\ell)&=\big\{\!(I,\{J,K\})\!\in\!\wt\cP(\ell)\!\!:
\,1\!\le\!|I|\!\le\!\ell\!-\!2\big\}.
\end{aligned}
\end{gather*}
For a commutative ring~$R$, a collection~$\cP$, and a $\cP$-tuple $(D_P)_{P\in\cP}$, 
we denote by $R[(D_P)_{P\in\cP}]$ the algebra over~$R$ generated by 
the components~$D_P$ of this tuple.
It becomes a graded algebra once the degree~$|D_P|$ of each~$D_P$ is specified.
If $J,K\!\subset\!\Z^+$ are disjoint subsets that are not both empty, let
$$\ep_{J,K}=\begin{cases}1,&\hbox{if}~\min(J\!\cup\!K)\!\in\!J;\\
-1,&\hbox{if}~\min(J\!\cup\!K)\!\in\!K.
\end{cases}$$
For two pairs $\{J,K\}$ and $\{J',K'\}$ of sets, we write 
\begin{alignat*}{3}
\{J,K\}&\preceq\{J',K'\} &\qquad&\hbox{if}\qquad  \hbox{either}&\quad
&J\subset J'~\hbox{and}~K\subset K' 
\quad\hbox{or} \quad J\!\subset\!K'~\hbox{and}~K\subset J';\\
\{J,K\}&\!\prl\!\{J',K'\} &\qquad&\hbox{if}&\qquad
&\{J,K\}\not\preceq\{J',K'\}~~\hbox{and}~~\{J',K'\}\not\preceq\{J,K\}.
\end{alignat*}
If in addition $J\!\sqcup\!K\!=\!J'\!\sqcup\!K'$, we write 
$$\{J,K\}\!\not\!\cap\,\{J',K'\}  \qquad\hbox{if}\qquad
J\not\subset J',K' ~~\hbox{and}~~ J',K' \not\subset J.$$
For $(I,\{J,K\}),(I',\{J',K'\})\!\in\!\wt\cP(\ell)$, we define
$$\big(I,\{J,K\}\!\big) \!\not\!\cap\big(I',\{J',K'\}\!\big)
 \qquad\hbox{if}\qquad \{J,K\}\!\prl\!\{J',K'\} ~~\hbox{and}~~
J\!\sqcup\!K\not\subset I'.$$

\subsection{The complex case}
\label{Cthm_subs}

\noindent
For $\{J,K\}\!\in\!\cP(\ell)$, we denote~by
$$D_{J,K}\subset\ov\cM_{\ell}$$  
the closure of the subspace consisting of the two-component curves such~that 
one of the components carries the marked points indexed by~$J$
and the other by~$K$.
This complex divisor is empty unless \hbox{$\{J,K\}\!\in\!\cP_{\bu}(\ell)$}.
The intersections of the divisors~$D_{J,K}$ are the \sf{boundary strata} 
in~$\ov\cM_{\ell}$.
We also use $D_{J,K}$ to denote 
the Poincare dual of~$D_{J,K}$ in $H^2(\ov\cM_{\ell};\Z)$.
If $\ell\!=\!4$, the divisors~$D_{J,K}$ are points and represent the three curves
in Figure~\ref{M04rel_fig}.\\

\noindent
Since $D_{J,K}$ and $D_{J',K'}$ are disjoint if $\{J,K\}\!\!\not\!\cap\{J',K'\}$,
\BE{cM04rel_e1a}
D_{J,K}D_{J',K'}=0\in H^4\big(\ov\cM_{\ell};\Z\big)
\qquad\forall~\{J,K\},\!\{J',K'\}\!\in\!\cP_{\bu}(\ell)
~\hbox{s.t.}~\{J,K\}\!\not\!\cap\,\{J',K'\}.\EE
Furthermore, 
\BE{cM04rel_e2} 
\sum_{\begin{subarray}{c} \{J,K\}\in\cP_{\bu}(\ell)\\
a,b\in J,\,c,d\in K \end{subarray}}\hspace{-.22in}D_{J,K}
-\!\!\! \sum_{\begin{subarray}{c} \{J,K\}\in\cP_{\bu}(\ell) \\
a,c\in J,\,b,d\in K \end{subarray}}\hspace{-.22in}D_{J,K}=0
\in H^2\big(\ov\cM_{\ell};\Z\big)
\quad\forall\,a,b,c,d\!\in\![\ell],\,a\!\neq\!b,c\!\neq\!d.\EE 
The $\ell\!=\!4$ cases of this identity, illustrated in Figure~\ref{M04rel_fig},
simply state that three distinct points in $\ov\cM_4\!\approx\!\C\P^1$
determine the same cohomology class.
The non-trivial $\ell\!\ge\!5$ cases of~\eref{cM04rel_e2},
i.e.~with $a\!\neq\!d$ and $b\!\neq\!c$ in addition to the four conditions
stated in~\eref{cM04rel_e2}, are obtained by pulling back
the $\ell\!=\!4$ relations by the forgetting morphisms
$$\ff\!:\ov\cM_{\ell}\lra \ov\cM_4$$
sending $a,b,c,d$ to $1,2,3,4$ and dropping the remaining marked points.\\

\noindent
We denote by
$$\cI_{\ell}\subset \Z\big[(D_{J,K})_{\{J,K\}\in\cP_{\bu}(\ell)}\big]$$
the ideal generated by the left-hand sides of~\eref{cM04rel_e1a} and~\eref{cM04rel_e2}.

\begin{figure}
\begin{pspicture}(-2.8,-1)(10,1.2)
\psset{unit=.3cm}
\pscircle*(5,0){.2}
\psline[linewidth=.07](4,-1)(8,3)\psline[linewidth=.07](4,1)(8,-3)
\pscircle*(6,1){.2}\pscircle*(7,2){.2}\pscircle*(6,-1){.2}\pscircle*(7,-2){.2}
\rput(6,2){\sm{$1$}}\rput(7,3){\sm{$2$}}\rput(6,-2){\sm{$3$}}\rput(7,-3){\sm{$4$}}
\rput(12,0){\begin{Large}$=$\end{Large}}
\pscircle*(17,0){.2}
\psline[linewidth=.07](16,-1)(20,3)\psline[linewidth=.07](16,1)(20,-3)
\pscircle*(18,1){.2}\pscircle*(19,2){.2}\pscircle*(18,-1){.2}\pscircle*(19,-2){.2}
\rput(18,2){\sm{$1$}}\rput(19,3){\sm{$3$}}\rput(18,-2){\sm{$2$}}\rput(19,-3){\sm{$4$}}
\rput(24,0){\begin{Large}$=$\end{Large}}
\pscircle*(29,0){.2}
\psline[linewidth=.07](28,-1)(32,3)\psline[linewidth=.07](28,1)(32,-3)
\pscircle*(30,1){.2}\pscircle*(31,2){.2}\pscircle*(30,-1){.2}\pscircle*(31,-2){.2}
\rput(30,2){\sm{$1$}}\rput(31,3){\sm{$4$}}\rput(30,-2){\sm{$2$}}\rput(31,-3){\sm{$3$}}
\end{pspicture}
\caption{The basic cases of~\eref{cM04rel_e2}, 
equivalences of three points in $H^2(\ov\cM_4;\Z)$.}
\label{M04rel_fig}
\end{figure}

\begin{thm}[{\cite[p550]{Keel}}]\label{HC_thm}
For every $\ell\!\in\!\Z^+$ with $\ell\!\ge\!3$, the homomorphism
\BE{HCthm_e}\Z\big[(D_{J,K})_{\{J,K\}\in\cP_{\bu}(\ell)}\big]\big/\cI_{\ell}
\lra H^*\big(\ov\cM_{\ell};\Z\big), \quad
 D_{J,K}\!+\!\cI_{\ell} \lra D_{J,K}\,,\EE
is a well-defined isomorphism of graded $\Z$-algebras 
if $|D_{J,K}|\!\equiv\!2$ for all \hbox{$\{J,K\}\in\cP_{\bu}(\ell)$}.
\end{thm}

\vspace{.1in}

\noindent 
In light of the relations \eref{cM04rel_e1a} and~\eref{cM04rel_e2}, 
the substance of this theorem is~that
\begin{enumerate}[label=$\C\arabic*$:,ref=$\C\arabic*$,leftmargin=*]

\item\label{HCgen_it} the divisor classes $D_{J,K}\!\in\!H^2(\ov\cM_{\ell};\Z)$
algebraically generate the $\Z$-algebra $H^*(\ov\cM_{\ell};\Z)$ 
or equivalently
the boundary strata of $\ov\cM_{\ell}$ linearly generate
the $\Z$-module $H^*(\ov\cM_{\ell};\Z)$;

\item\label{HCrel_it} the ring homomorphism~\eref{HCthm_e} is injective.

\end{enumerate}
Since $\ov\cM_3\!=\!\{\tn{pt}\}$ and $\ov\cM_4\!\approx\!\C\P^1$,
these statements clearly hold for $\ell\!=\!3,4$.
We re-establish~\ref{HCgen_it} and~\ref{HCrel_it} for arbitrary $\ell\!\ge\!4$ 
in Sections~\ref{GenC_subs} and~\ref{Calg_subs}, respectively, via induction on~$\ell$.
The foundation for relating the moduli spaces~$\ov\cM_{\ell}$
and~$\ov\cM_{\ell+1}$ is provided by the forgetful morphisms
\BE{ffCdfn_e}\ff_{\ell+1}\!: \ov\cM_{\ell+1}\lra\ov\cM_{\ell} 
\qquad\hbox{and}\qquad
\ff_{123,\ell+1}\!: \ov\cM_{\ell+1}\lra\ov\cM_4\EE
dropping the marked point~$z_{\ell+1}$ and 
all marked points other than $z_1,z_2,z_3,z_{\ell+1}$, respectively.
By the Universal Coefficients Theorem, 
Theorem~\ref{HC_thm} also holds with~$\Z$ replaced by any commutative ring
with~unity.

\subsection{The real case}
\label{Rthm_subs}

\noindent
For $\{J,K\}\!\in\!\cP(\ell)$, we denote~by
$$\R E_{J,K}\subset\R\ov\cM_{0,\ell} \qquad 
\big(\hbox{resp.~$\R H_{J,K}\subset\R\ov\cM_{0,\ell}$}\big)$$  
the closure of the subspace consisting of the two-component curves $(\Si,\si)$ such~that 
\begin{enumerate}[label=$\bullet$,leftmargin=*]

 \item the fixed locus~$\Si^{\si}$ of the involution~$\si$ on~$\Si$ is a point
(resp.~a wedge of two circles),

\item one of the components carries the first marked point~$z_j^+$ of each conjugate pair
$(z_j^+,z_j^-)$ with $j\!\in\!J$, and
the other carries the first marked point~$z_k^+$ of each conjugate pair
$(z_k^+,z_k^-)$ with $k\!\in\!K$.

\end{enumerate}
The subspaces $\R E_{J,K},\R H_{J,K}\!\subset\!\R\ov\cM_{0,\ell}$
are smooth real hypersurfaces, which we call \sf{boundary hypersurfaces}.
The hypersurfaces~$\R E_{J,K}$ are orientable, while
the hypersurfaces~$\R H_{J,K}$ are not orientable if $\ell\!\ge\!3$
and $J,K\!\neq\!\eset$ ($\R H_{J,K}\!=\!\eset$ unless $J,K\!\neq\!\eset$).
If $\ell\!=\!2$, the real hypersurfaces~$\R E_{J,K}$ are points and 
represent the two curves on the left side on the first line
in Figure~\ref{TopolRel_fig}.
If $\ell\!=\!3$, there are four real hypersurfaces~$\R E_{J,K}$;
each of them is isomorphic to~$S^2$.\\

\noindent
For $(I,\{J,K\})\!\in\!\wt\cP(\ell)$, we denote~by
$$\R D_{I;J,K}\subset\R\ov\cM_{0,\ell} $$  
the closure of the subspace consisting of the three-component curves $(\Si,\si)$ such~that
\begin{enumerate}[label=$\bullet$,leftmargin=*]

\item one of the components is fixed by the involution~$\si$ and 
carries the pairs~$(z_i^+,z_i^-)$  of marked points with~$i\!\in\!I$,

\item the other two components are interchanged by~$\si$, and 
one of them carries the marked points~$z_j^+$ with $j\!\in\!J$
and~$z_k^-$ with $k\!\in\!K$
(and thus the other component carries the marked points~$z_k^+$ with $k\!\in\!K$
and~$z_j^-$ with $j\!\in\!J$).

\end{enumerate}
The subspaces $\R D_{I;J,K}\!\subset\!\R\ov\cM_{0,\ell}$
are smooth orientable codimension~2 submanifolds, which we call \sf{boundary divisors}.
If $\R D_{I;J,K}\!\neq\!\eset$, then
\hbox{$(I,\{J,K\})\!\in\!\wt\cP_{\bu}(\ell)$}.
The intersections of the manifolds $\R E_{J,K}$, $\R H_{J,K}$, and $\R D_{I;J,K}$
are the \sf{boundary strata} in~$\R\ov\cM_{0,\ell}$.\\

\noindent 
If $\ell\!=\!3$, there are six boundary divisors~$\R D_{I;J,K}$.
Each of them is isomorphic to~$S^1$, is disjoint from
the other five boundary divisors~$\R D_{I';J',K'}$ in~$\R\ov\cM_{0,3}$,
and intersects two of the four real hypersurfaces $\R E_{J',K'}$ at one point each
as indicated by the bottom diagram in Figure~\ref{TopolRel_fig}.
The three curves on the right side on the first line in Figure~\ref{TopolRel_fig}
represent generic elements of three of the submanifolds~$\R D_{I;J,K}$;
the complement of such elements in the corresponding~$\R D_{I;J,K}$ 
is the intersection of the latter with
the hypersurfaces~$\R E_{J',K'}$ and~$\R H_{J',K'}$.\\

\noindent
We orient the entire moduli space $\R\ov\cM_{0,\ell}$,
the boundary hypersurfaces~$\R E_{J,K}$, and the boundary divisors~$\R D_{I;J,K}$
in Section~\ref{Rhomol_subs1}.
We also denote by $\R E_{J,K}$ and $\R D_{I;J,K}$ the corresponding elements
of $H_{2\ell-4}(\R\ov\cM_{0,\ell};\Q)$ and $H_{2\ell-5}(\R\ov\cM_{0,\ell};\Q)$, respectively,
and their Poincare duals in~$H^1(\R\ov\cM_{0,\ell};\Q)$ and~$H^2(\R\ov\cM_{0,\ell};\Q)$.\\

\begin{figure}
\begin{pspicture}(-.5,-.2)(10,5.4)
\psset{unit=.4cm}
\psline[linewidth=.05](2,9)(6,13)\psline[linewidth=.05](2,11)(6,7)
\pscircle*(3.5,10.5){.15}\pscircle*(5,12){.15}
\pscircle*(3.5,9.5){.15}\pscircle*(5,8){.15}
\rput(3.2,8.9){\sm{$1^-$}}\rput(4.5,7.5){\sm{$2^+$}}
\rput(3.5,11.1){\sm{$1^+$}}\rput(5,12.6){\sm{$2^-$}}
\psline[linewidth=.025]{<->}(5.5,8.5)(5.5,11.5)\rput(6,10){$\si$}
\rput(8.5,10){\begin{Large}$=~-$\end{Large}}
\psline[linewidth=.05](10,9)(14,13)\psline[linewidth=.05](10,11)(14,7)
\pscircle*(11.5,10.5){.15}\pscircle*(13,12){.15}
\pscircle*(11.5,9.5){.15}\pscircle*(13,8){.15}
\rput(11.1,8.9){\sm{$1^-$}}\rput(12.5,7.5){\sm{$2^-$}}
\rput(11.5,11.1){\sm{$1^+$}}\rput(13,12.6){\sm{$2^+$}}
\psline[linewidth=.025]{<->}(13.5,8.5)(13.5,11.5)\rput(14,10){$\si$}
\psline[linewidth=.05](19,12.5)(19,7.5)
\psline[linewidth=.02](18.5,12)(22,12)\psline[linewidth=.02](18.5,8)(22,8)
\pscircle*(20,12){.15}\pscircle*(21,12){.15}
\pscircle*(19,11){.15}\pscircle*(19,9){.15}
\rput(20,12.7){\sm{$1^+$}}\rput(21.4,12.7){\sm{$2^+$}}
\rput(18.4,11.3){\sm{$3^+$}}\rput(18.4,9.3){\sm{$3^-$}}
\pscircle*(20,8){.15}\pscircle*(21,8){.15}
\rput(20,7.3){\sm{$1^-$}}\rput(21.4,7.3){\sm{$2^-$}}
\psline[linewidth=.025]{<->}(21.5,8.6)(21.5,11.4)\rput(22,10){$\si$}
\rput(23.7,10){\begin{Large}$-$\end{Large}}
\psline[linewidth=.05](26,12.5)(26,7.5)
\psline[linewidth=.02](25.5,12)(29,12)\psline[linewidth=.02](25.5,8)(29,8)
\pscircle*(27,12){.15}\pscircle*(28,12){.15}
\pscircle*(26,11){.15}\pscircle*(26,9){.15}
\rput(27,12.7){\sm{$1^+$}}\rput(28.4,12.7){\sm{$3^+$}}
\rput(25.4,11.3){\sm{$2^+$}}\rput(25.4,9.3){\sm{$2^-$}}
\pscircle*(27,8){.15}\pscircle*(28,8){.15}
\rput(27,7.3){\sm{$1^-$}}\rput(28.4,7.3){\sm{$3^-$}}
\psline[linewidth=.025]{<->}(28.5,8.6)(28.5,11.4)\rput(29,10){$\si$}
\rput(30.5,10){\begin{Large}$=$\end{Large}}
\psline[linewidth=.05](33,12.5)(33,7.5)
\psline[linewidth=.02](32.5,12)(36,12)\psline[linewidth=.02](32.5,8)(36,8)
\pscircle*(34,12){.15}\pscircle*(35,12){.15}
\pscircle*(33,11){.15}\pscircle*(33,9){.15}
\rput(34,12.7){\sm{$2^+$}}\rput(35.4,12.7){\sm{$3^-$}}
\rput(32.4,11.3){\sm{$1^+$}}\rput(32.4,9.3){\sm{$1^-$}}
\pscircle*(34,8){.15}\pscircle*(35,8){.15}
\rput(34,7.3){\sm{$2^-$}}\rput(35.4,7.3){\sm{$3^+$}}
\psline[linewidth=.025]{<->}(35.5,8.6)(35.5,11.4)\rput(36,10){$\si$}
\psline[linewidth=.05](25,5)(25,0)\psline[linewidth=.05,arrowsize=6pt]{->}(25,5)(25,1)
\psline[linewidth=.05](10,5)(10,0)\psline[linewidth=.05,arrowsize=6pt]{->}(10,0)(10,4)
\rput(25.8,3){$3^-$}\rput(9.4,2){$3^+$}
\psline[linewidth=.05](10,0)(25,0)\psline[linewidth=.05,arrowsize=6pt]{->}(25,5)(12,5)
\rput(15,-.6){$2^-$}\rput(20,5.6){$2^+$}
\psline[linewidth=.05](10,5)(25,5)\psline[linewidth=.05,arrowsize=6pt]{->}(10,0)(23,0)
\psline[linewidth=.05](10,5)(25,0)\psline[linewidth=.05,arrowsize=6pt]{->}(25,0)(12,4.33)
\psline[linewidth=.05](10,0)(17,2.33)
\psline[linewidth=.05](25,5)(18,2.67)\psline[linewidth=.05,arrowsize=6pt]{->}(18,2.67)(23,4.33)
\rput(21.2,2){$1^+$}\rput(14,2){$1^-$}
\pscircle*(10,5){.3}\rput(9.3,5){$\eset$}
\pscircle*(10,0){.3}\rput(9.3,0){$3$}
\pscircle*(25,5){.3}\rput(25.7,5){$2$}
\pscircle*(25,0){.3}\rput(26,0){$23$}
\end{pspicture}
\caption{The first line represents basic cases of~\eref{RcM04rel_e2a} and~\eref{RcM04rel_e2b},
an equivalence of two points in $H^1(\R\ov\cM_{0,2};\Q)$ and 
a relation between three loops in $H^2(\R\ov\cM_{0,3};\Q)$.
The bottom diagram represents the intersection pattern 
of the six boundary divisors $\R D_{I;J,K}\!\approx\!S^1$
with the four boundary hypersurfaces $\R E_{J',K'}\!\approx\!S^2$ in~$\R\ov\cM_{0,3}$.
The former are labeled by the unique element of $I\!\subset\![3]$
and the sign of $(-1)^{|J|}\!=\!(-1)^{|K|}$;
the latter are labeled by the subset $J',K'\!\subset\![3]$ not containing~1.}
\label{TopolRel_fig}
\end{figure}

\noindent
Since the hypersurfaces~$\R E_{J,K}$ and~$\R E_{J',K'}$ are disjoint 
unless $\{J,K\}\!=\!\{J',K'\}$
and the normal bundle to~$\R E_{J,K}$ is trivial,
\BE{RcM04rel_e1a}
\R E_{J,K}\R E_{J',K'}=0\in H^2\big(\R\ov\cM_{0,\ell};\Q\big)
\quad\forall~\{J,K\},\{J',K'\}\!\in\!\cP(\ell)\,.\EE
Since $\R E_{J,K}$ and $\R D_{I';J',K'}$ are disjoint unless 
either $J'\!\subset\!J$ and $K'\!\subset\!K$ or $J'\!\subset\!K$ and $K'\!\subset\!J$,
\BE{RcM04rel_e1b}\begin{split}
\R E_{J,K}\R D_{I';J',K'}=0\in H^3\big(\R\ov\cM_{0,\ell};\Q\big)
\quad\forall~&\{J,K\}\!\in\!\cP(\ell),\,
(I',\{J',K'\})\!\in\!\wt\cP_{\bu}(\ell)\\
&~\hbox{s.t.}~\{J',K'\}\!\not\preceq\!\{J,K\}\,.
\end{split}\EE
Since $\R D_{I;J,K}$ and $\R D_{I';J',K'}$ are disjoint unless 
$\{J,K\}\!\preceq\!\{J',K'\}$, or $\{J',K'\}\!\preceq\!\{J,K\}$,
or \hbox{$J,K\!\subset\!I'$},
\BE{RcM04rel_e1c}\begin{split}
\R D_{I;J,K}\R D_{I';J',K'}=0\in H^4\big(\R\ov\cM_{0,\ell};\Q\big) \quad
\forall~\big(I,\{J,K\}\!\big),\big(I',\{J',K'\}\!\big)\!\in\!\wt\cP_{\bu}(\ell)~~&\\
~\hbox{s.t.}~\big(I,\{J,K\}\!\big)\!\not\!\cap\big(I',\{J',K'\}\!\big)\,.
\end{split}\EE
We also have that 
\BE{RcM04rel_e2a}
\sum_{(J,K)\in\cP(\ell)}\hspace{-.18in}\R E_{J,K}=0\in 
H^1\big(\R\ov\cM_{0,\ell};\Q\big).\EE
Furthermore, 
\BE{RcM04rel_e2b}\begin{split}
\sum_{\begin{subarray}{c} (I,\{J,K\})\in\wt\cP_{\bu}(\ell)\\
a,b\in J,\,c\in I \end{subarray}}\hspace{-.33in}\ep_{J,K}\R D_{I;J,K}
&-\sum_{\begin{subarray}{c} (I,\{J,K\})\in\wt\cP_{\bu}(\ell)\\
a,c\in J,\,b\in I \end{subarray}}\hspace{-.33in}\ep_{J,K}\R D_{I;J,K}
-\!\!\!\sum_{\begin{subarray}{c} (I,\{J,K\})\in\wt\cP_{\bu}(\ell)\\
a\in I,\,b\in J,\,c\in K \end{subarray}}\hspace{-.33in}\ep_{J,K}\R D_{I;J,K}
=0\in H^2\big(\R\ov\cM_{0,\ell};\Q\big)
\end{split}\EE 
for all $a,b,c\!\in\![\ell]$ distinct.\\

\noindent
The $\ell\!=\!2$ case of~\eref{RcM04rel_e2a},
illustrated by the first identity on the top line in Figure~\ref{TopolRel_fig},
simply states that two distinct points in $\R\ov\cM_{0,2}\!\approx\!S^1$
determine the same cohomology class
(the diagram on the left side of this identity represents 
a minus point in $\R\ov\cM_{0,2}$ with our choices of the orientations).
The $\ell\!\ge\!3$ cases of~\eref{RcM04rel_e2a} are obtained by pulling back
the $\ell\!=\!2$ relation by the forgetful morphism
$$\ff_{3^{\pm}}^{\R}\!\circ\!\ldots\!\circ\!\ff_{\ell^{\pm}}^{\R}\!:\R\ov\cM_{0,\ell}\lra \R\ov\cM_{0,2}$$
dropping the conjugate pairs $(z_3^+,z_3^-),\ldots,(z_{\ell}^+,z_{\ell}^-)$ of marked points;
see Section~\ref{Rhomol_subs1}\\

\noindent
If $\ell\!=\!3$, there are three independent cases of~\eref{RcM04rel_e2b}.
They are determined by the value of \hbox{$a\!=\!1,2,3$} and 
correspond to the upper left, lower left, and upper right 
triangles in the bottom diagram of Figure~\ref{TopolRel_fig}
being null-homologous in~$\R\ov\cM_{0,3}$ over~$\Q$.
This implies that the fourth triangle in this diagram is null-homologous as well.
The second identity on the top line in Figure~\ref{TopolRel_fig} 
illustrates the $a\!=\!1$ case of~\eref{RcM04rel_e2b}. 
The three $\ell\!=\!3$ cases of~\eref{RcM04rel_e2b}
are established in Section~3 of~\cite{RealEnum}.
The $\ell\!\ge\!4$ cases of~\eref{RcM04rel_e2b} are obtained by pulling back
the $\ell\!=\!3$ relations by the forgetful morphism
$$\ff_{4^{\pm}}^{\R}\!\circ\!\ldots\!\circ\!\ff_{\ell^{\pm}}^{\R}\!:\R\ov\cM_{0,\ell}\lra \R\ov\cM_{0,3}$$
dropping the conjugate pairs $(z_4^+,z_4^-),\ldots,(z_{\ell}^+,z_{\ell}^-)$ of marked points
and then interchanging conjugate pairs of marked points; see Section~\ref{Rhomol_subs1}.\\

\noindent
We denote by
$$\cI_{0,\ell}\subset \Q\big[(\R E_{J,K})_{\{J,K\}\in\cP(\ell)},
(\R D_{I;J,K})_{(I,\{J,K\})\in\wt\cP_{\bu}(\ell)}\big]$$
the ideal generated by the left-hand sides of~\eref{RcM04rel_e1a}-\eref{RcM04rel_e2b}.

\begin{thm}\label{HR_thm}
For every $\ell\!\in\!\Z^+$ with $\ell\!\ge\!2$, the homomorphism
\BE{HRthm_e}\begin{split}
&\Q\big[(\R E_{J,K})_{\{J,K\}\in\cP(\ell)},
(\R D_{I;J,K})_{(I,\{J,K\})\in\wt\cP_{\bu}(\ell)}\big]\big/\cI_{0,\ell}
\lra H^*\big(\R\ov\cM_{0,\ell};\Q\big), \\
&\hspace{.5in} 
\R E_{J,K}\!+\!\cI_{0,\ell} \lra \R E_{J,K}\,,\quad
\R D_{I;J,K}\!+\cI_{0,\ell} \lra\R D_{I;J,K}\,,
\end{split}\EE
is a well-defined isomorphism of graded $\Z$-algebras 
if $|\R E_{J,K}|\!\equiv\!1$ and \hbox{$|\R D_{I;J,K}|\!\equiv\!2$}
for all \hbox{$\{J,K\}\!\in\!\cP(\ell)$} and 
\hbox{$(I,\{J,K\})\!\in\!\wt\cP_{\bu}(\ell)$}, respectively.
\end{thm}

\vspace{.1in}

\noindent 
In light of the relations~\eref{RcM04rel_e1a}-\eref{RcM04rel_e2b}, 
the substance of this theorem is~that
\begin{enumerate}[label=$\R\arabic*$:,ref=$\R\arabic*$,leftmargin=*]

\item\label{HRgen_it} the classes $\R E_{J,K}\!\in\!H^1(\R\ov\cM_{0,\ell};\Q)$ and
$\R D_{I;J,K}\!\in\!H^2(\R\ov\cM_{0,\ell};\Q)$
algebraically generate the algebra $H^*(\R\ov\cM_{0,\ell};\Q)$ over~$\Q$
or equivalently
the boundary strata obtained as intersections of 
the codimension~1 and~2 submanifolds $\R E_{J,K}$ and $\R D_{I;J,K}$
linearly generate the vector space $H^*(\R\ov\cM_{0,\ell};\Q)$ over~$\Q$;

\item\label{HRrel_it} the ring homomorphism~\eref{HRthm_e} is injective.

\end{enumerate}
Since $\R\ov\cM_{0,2}\!\approx\!S^1$, these statements clearly hold for $\ell\!=\!2$.
For $\ell\!=\!3$, they are essentially established in Section~3 of~\cite{RealEnum}
by explicitly describing the orientable threefold~$\R\ov\cM_{0,3}$.
We establish~\ref{HRgen_it} and~\ref{HRrel_it} for arbitrary $\ell\!\ge\!3$ 
in Sections~\ref{GenR_subs}
and~\ref{Ralg_subs}, respectively, via induction on~$\ell$.
The foundation for relating the moduli spaces~$\R\ov\cM_{0,\ell}$
and~$\R\ov\cM_{0,\ell+1}$ is provided by the forgetful morphisms
\BE{ffRdfn_e}\ff_{(\ell+1)^{\pm}}^{\R}\!:\R\ov\cM_{0,\ell+1}\lra\R\ov\cM_{0,\ell} 
\quad\hbox{and}\quad
\ff_{1^{\pm}2^+(\ell+1)^+}^{\R}\!:\R\ov\cM_{0,\ell+1}\lra\ov\cM_4\EE
dropping the marked points $z_{\ell+1}^+,z_{\ell+1}^-$ and 
all marked points other than $z_1^+,z_1^-,z_2^+,z_{\ell+1}^+$,
respectively.
The proof implies that Theorem~\ref{HC_thm} also holds with~$\Q$ replaced by 
any commutative ring with~unity in which~2 is a unit.\\

\noindent
By the proof of~\ref{HRgen_it} in Sections~\ref{BlHom_subs} and~\ref{GenR_subs},
the cohomology classes $\R E_{J,K}$, $\R H_{J,K}$, and $\R D_{I;J,K}$
algebraically generate $H^*(\R\ov\cM_{0,\ell};\Z)$ or  equivalently
the orientable boundary strata of $\R\ov\cM_{0,\ell}$ linearly generate $H^*(\R\ov\cM_{0,\ell};\Z)$.
However, the orientable intersections with~$\R H_{J,K}$ contribute only torsion.
As shown in \cite[Section~3]{RealEnum}, this contribution is non-trivial if $\ell\!=\!3$
and thus for all $\ell\!\ge\!3$.\\

\noindent
As an example, Theorem~\ref{HR_thm} implies that 
$$H^*\!\big(\R\ov\cM_{0,2};\Q\big)\approx \Q[x]/x^2$$
with $|x|\!=\!1$.
The generator $x$ above can be mapped to either 
$\R E_{12,\eset}\!=\!\R E_{\eset,12}$  or $\R E_{1,2}\!=\!\R E_{2,1}$.
Theorem~\ref{HR_thm} also gives
\begin{gather*}
H^*\!\big(\R\ov\cM_{0,3};\Q\big)\approx 
\Q\big[x_1,x_2,x_3,y_1,y_2,y_3\big]\big/\cI, \qquad\hbox{where}\\
\cI=\Q\big[\{x_ix_j\}_{i,j\in[3]},\{y_iy_j\}_{i,j\in[3]},
\{x_iy_j\}_{i,j\in[3],i\neq j},
\{x_1y_1\!-\!x_2y_2,x_1y_1\!-\!x_3y_3\}\big],
\end{gather*}
with $|x_1|,|x_2|,|x_3|\!=\!1$ and $|y_1|,|y_2|,|y_3|\!=\!2$.
The generators $x_1$, $x_2$, $x_3$, $y_1$, $y_2$, and~$y_3$ could be mapped to 
$\R E_{1,23}$, $\R E_{2,13}$, $\R E_{3,12}$, $\R D_{1;23,\eset}$,  $\R D_{2;13,\eset}$,
and $\R D_{3;12,\eset}$, respectively. 

\begin{figure}
\begin{pspicture}(-.5,-.2)(10,2.2)
\psset{unit=.4cm}
\psline[linewidth=.05](3,5)(3,0)
\psline[linewidth=.02](2.5,4.5)(6,4.5)\psline[linewidth=.02](2.5,.5)(6,.5)
\pscircle*(4,4.5){.15}\pscircle*(5,4.5){.15}
\pscircle*(3,3.5){.15}\pscircle*(3,1.5){.15}
\psline[linewidth=.025]{<->}(5.5,1.1)(5.5,3.9)\rput(6,2.5){$\si$}
\rput(4,5.2){\sm{$1^+$}}\rput(5.4,5.2){\sm{$2^+$}}
\rput(2.4,3.8){\sm{$3^+$}}\rput(2.4,1.8){\sm{$3^-$}}
\pscircle*(4,.5){.15}\pscircle*(5,.5){.15}
\rput(4,-.2){\sm{$1^-$}}\rput(5.4,-.2){\sm{$2^-$}}
\psline[linewidth=.05](11,5)(11,0)
\psline[linewidth=.02](10.5,4.5)(14,4.5)\psline[linewidth=.02](10.5,.5)(14,.5)
\pscircle*(12,4.5){.15}\pscircle*(13,4.5){.15}
\pscircle*(11,3.5){.15}\pscircle*(11,1.5){.15}
\psline[linewidth=.025]{<->}(13.5,1.1)(13.5,3.9)\rput(14,2.5){$\si$}
\rput(12,5.2){\sm{$1^+$}}\rput(13.4,5.2){\sm{$2^-$}}
\rput(10.4,3.8){\sm{$3^+$}}\rput(10.4,1.8){\sm{$3^-$}}
\pscircle*(12,.5){.15}\pscircle*(13,.5){.15}
\rput(12,-.2){\sm{$1^-$}}\rput(13.4,-.2){\sm{$2^+$}}
\rput(8,2.5){\begin{Large}$+$\end{Large}}
\rput(17,2.5){\begin{Large}$=$\end{Large}}
\psline[linewidth=.05](22,5)(22,0)
\psline[linewidth=.02](21.5,4.5)(25,4.5)\psline[linewidth=.02](21.5,.5)(25,.5)
\pscircle*(23,4.5){.15}\pscircle*(24,4.5){.15}
\pscircle*(22,3.5){.15}\pscircle*(22,1.5){.15}
\psline[linewidth=.025]{<->}(24.5,1.1)(24.5,3.9)\rput(25,2.5){$\si$}
\rput(23,5.2){\sm{$1^+$}}\rput(24.4,5.2){\sm{$3^+$}}
\rput(21.4,3.8){\sm{$2^+$}}\rput(21.4,1.8){\sm{$2^-$}}
\pscircle*(23,.5){.15}\pscircle*(24,.5){.15}
\rput(23,-.2){\sm{$1^-$}}\rput(24.4,-.2){\sm{$3^-$}}
\psline[linewidth=.05](30,5)(30,0)
\psline[linewidth=.02](29.5,4.5)(33,4.5)\psline[linewidth=.02](29.5,.5)(33,.5)
\pscircle*(31,4.5){.15}\pscircle*(32,4.5){.15}
\pscircle*(30,3.5){.15}\pscircle*(30,1.5){.15}
\psline[linewidth=.025]{<->}(32.5,1.1)(32.5,3.9)\rput(33,2.5){$\si$}
\rput(31,5.2){\sm{$1^+$}}\rput(32.4,5.2){\sm{$3^-$}}
\rput(29.4,3.8){\sm{$2^+$}}\rput(29.4,1.8){\sm{$2^-$}}
\pscircle*(31,.5){.15}\pscircle*(32,.5){.15}
\rput(31,-.2){\sm{$1^-$}}\rput(32.4,-.2){\sm{$3^+$}}
\rput(27,2.5){\begin{Large}$+$\end{Large}}
\end{pspicture}
\caption{A basic case of~\eref{RcM04rel_e2b2},  
a relation between four loops in $H^2(\R\ov\cM_{0,3};\Q)$.}
\label{TopolRel_fig2}
\end{figure}

\vspace{.1in}

\begin{rmk}\label{M03_rmk}
The relations~\eref{RcM04rel_e2b} imply that 
\BE{RcM04rel_e2b2}
\sum_{\begin{subarray}{c} (I,\{J,K\})\in\wt\cP_{\bu}(\ell)\\
a\in J,\,b\not\in I,\,c\in I \end{subarray}}\hspace{-.33in}\ep_{J,K}\R D_{I;J,K}
-\!\!\sum_{\begin{subarray}{c} (I,\{J,K\})\in\wt\cP_{\bu}(\ell)\\
a\in J,\,b\in I,\,c\not\in I\end{subarray}}\hspace{-.33in}\ep_{J,K}\R D_{I;J,K}
=0\in H^2\big(\R\ov\cM_{0,\ell};\Q\big)\EE 
with $a,b,c\!\in\![\ell]$ such that $a\!\neq\!b,c$.
An $\ell\!=\!3$ case of~\eref{RcM04rel_e2b2} is illustrated in Figure~\ref{TopolRel_fig2}.
The three $\ell\!=\!3$ cases of~\eref{RcM04rel_e2b2}
correspond to the three quadrilaterals in the bottom diagram of Figure~\ref{TopolRel_fig}
being null-homologous in~$\R\ov\cM_{0,3}$ over~$\Q$.
By Remark~3.4 in~\cite{RealEnum}, 
the left-hand side of~\eref{RcM04rel_e2b2} is null-homologous in~$\R\ov\cM_{0,\ell}$
even over~$\Z$.
In contrast, the left-hand side of~\eref{RcM04rel_e2b} in an $\ell\!=\!3$ case
is the unique non-trivial torsion element of $H^*(\R\ov\cM_{0,3};\Z)$.
\end{rmk}

\section{Generation of the (co)homology}
\label{CohGener_sec}

\noindent
Property~\ref{HCgen_it} on page~\pageref{HCgen_it} is a fairly straightforward consequence
of the decomposition of the holomorphic map 
\BE{CMlind_e}\big(\ff_{\ell+1},\ff_{123,\ell+1}\big)\!: 
\ov\cM_{\ell+1}\lra\ov\cM_{\ell}\!\times\!\ov\cM_4\EE
as a sequence of (holomorphic) blowups $\pi_{\vr}\!:X_{\vr}\!\lra\!X_{\vr-1}$
constructed inductively in Section~1 of~\cite{Keel}.
The spaces~$X_{\vr}$ are described as quotients of~$\ov\cM_{\ell+1}$
by explicit equivalence relations in Section~1.1 of~\cite{RDMbl}.
This approach is adapted to explicitly decompose the smooth map 
\BE{RMlind_e}\big(\ff_{(\ell+1)^{\pm}}^{\R},\ff_{1^{\pm}2^+(\ell+1)^+}^{\R}\big)\!: 
\R\ov\cM_{0,\ell+1}\lra \R\ov\cM_{0,\ell}\!\times\!\ov\cM_4\EE
as a sequence of blowups $\pi_{\vr}\!:X_{\vr}\!\lra\!X_{\vr-1}$ of three different types
in Section~1.2 of~\cite{RDMbl}. 
We use this sequence to inductively confirm property~\ref{HRgen_it} 
on page~\pageref{HRgen_it} in Section~\ref{GenR_subs}.

\subsection{Smooth blowups}
\label{blowups_subs}

\noindent
Below we summarize the properties of 
the smooth \sf{real}, \sf{complex}, and \sf{augmented} blowups,
\BE{blowdownmaps_e}\pi\!:\BLR_YX\lra X, \qquad \pi\!:\BLC_YX\lra X, 
\qquad\hbox{and}\qquad \pi\!:\BLau_YX\lra X\EE
of a smooth manifold~$X$ along a closed codimension~$\fc$ submanifold~$Y$ 
that are relevant for our present purposes.
These blowups are described from a local perspective in Section~3 of~\cite{RDMbl} 
and from a global one in~\cite{blowups};
the former is generally easier to use, 
while the latter makes the changes in the topology more evident.
The complex and augmented blowups depend on choices of auxiliary structures 
on a neighborhood of~$Y$ in~$X$;
the augmented blowup also involves a choice of \hbox{$\fc_1\!\in\![\fc\!-\!1]$}.
Such an auxiliary structure determines a complex structure on 
the normal bundle~$\cN_XY$ of~$Y$ in~$X$ in the complex case 
and a distinguished subbundle $\cN_X^{\fc_1}\!Y\!\subset\!\cN_XY$ of corank~$\fc_1$
in the $\fc_1$-augmented case.
The cases relevant for Theorems~\ref{HC_thm} and~\ref{HR_thm}
are the real blowup with $\fc\!=\!3$, the complex blowup with $\fc/2\!=\!2$,
and the augmented blowup with $(\fc,\fc_1)\!=\!(3,1)$.
For $k\!=\!0,1$, let \hbox{$\tau_Y^k\!\lra\!Y$} 
denote the trivial rank~$k$ real line bundle.\\

\noindent
In all three cases, the maps in~\eref{blowdownmaps_e} are surjective.
The restrictions
\begin{gather*}
\pi\!:\BLR_YX\!-\!\pi^{-1}(Y)\lra X\!-\!Y, \qquad 
\pi\!:\BLC_YX\!-\!\pi^{-1}(Y)\lra X\!-\!Y,\\
\hbox{and}\qquad 
\pi\!:\BLau_YX\!-\!\pi^{-1}(Y)\lra X\!-\!Y
\end{gather*}
are diffeomorphisms.
The preimage \hbox{$\pi^{-1}(U)\!\subset\!\wt{X}$} 
of a tubular neighborhood $U\!\subset\!X$ of~$Y$ deformation retracts onto
the \sf{exceptional locus} \hbox{$\pi^{-1}(Y)\!\!\subset\!\wt{X}$}.\\

\noindent
The restrictions of~$\pi$ to the exceptional loci
$$\pi\!:\bE_Y^{\R}X\!\equiv\!\pi^{-1}(Y)\lra Y \qquad\hbox{and}\qquad
\pi\!:\bE_Y^{\C}X\!\equiv\!\pi^{-1}(Y)\lra Y$$
in the first and second cases in~\eref{blowdownmaps_e}
are isomorphic to the real and complex projectivizations, 
$\R\P(\cN_XY)$ and $\C\P(\cN_XY)$, respectively, of~$\cN_XY$.
The normal bundles to the exceptional loci in these two cases are
isomorphic to the real and complex tautological line bundles:
\BE{cNEisom_e}\cN_{\BLR_YX}\!\big(\bE_Y^{\R}X\big)\approx\ga_{\cN_XY}^{\R}\lra\R\P(\cN_XY)
\quad\hbox{and}\quad
\cN_{\BLC_YX}\!\big(\bE_Y^{\C}X\big)\approx\ga_{\cN_XY}^{\C}\lra\C\P(\cN_XY)\,.\EE

\vspace{.15in}

\noindent
In the last case in~\eref{blowdownmaps_e}, the exceptional locus
$$\bE_Y^{\fc_1}X\equiv\pi^{-1}(Y)\subset\wt{X}$$ 
is the union of two closed submanifolds, $\bE^0_YX$ and~$\bE^-_YX$, so that 
$$\big(\bE^0_YX,\bE^0_YX\!\cap\!\bE^-_YX\big)\approx
\big(\R\P\big(\cN_YX/\cN_Y^{\fc_1}X\!\oplus\!\tau_Y^1\big),
\R\P\big(\cN_YX/\cN_Y^{\fc_1}X\!\oplus\!\tau_Y^0\big)\big)$$
as fiber bundle pairs over~$Y$, while 
$$\pi\!:\bE^-_YX\lra Y$$
factors through an $S^{\fc-\fc_1}$-fiber bundle $\bE^-_YX\!\lra\!\bE^0_YX\!\cap\!\bE^-_YX$.
In the relevant $(\fc,\fc_1)\!=\!(3,1)$ case, 
$\bE^-_YX$ is the $S^2$-sphere bundle over~$Y$ obtained by collapsing 
each fiber of the $\R\P^1$-fiber subbundle 
\hbox{$\R\P(\cN_Y^{\fc_1}X)\!\subset\!\R\P(\cN_YX)$} to a point.
The section of \hbox{$\bE_Y^{\fc_1}X\!\lra\!Y$}
determined by this collapse is~$\bE^0_YX\!\cap\!\bE^-_YX$.\\

\noindent
Suppose $X$ is oriented.
A complex blowup~$\BLC_YX$ of~$X$ along~$Y$ and 
the corresponding exceptional divisor~$\bE_Y^{\C}X$ then inherit orientations from~$X$. 
We denote the homology class of the latter in~$\bE_Y^{\C}\!X$ and in~$\BLC_YX$ 
by~$[\bE_Y^{\C}\!X]$.
A real blowup~$\BLR_YX$ and 
a $\fc_1$-augmented blowup~$\BLau_YX$ are orientable if and only if 
the codimension~$\fc$ of~$Y$ in~$X$ is odd;
if so, they inherit orientations from~$X$.
If $\fc\!\not\in\!2\Z$ and $\fc_1\!=\!1$, 
orientations on~$Y$, $\cN_XY$, and~$\cN_Y^{\fc_1}\!X$ determine orientations
on~$\bE_Y^0\!X$ and~$\bE_Y^-X$ and thus homology classes
$[\bE_Y^0X]$ in~$\bE_Y^0X$ and 
$[\bE_Y^-X]$ in~$\bE_Y^-X$, respectively.
We denote the images of these homology classes
in~$\bE_Y^1X$ and~$\BLaua_YX$ in the same ways.\\

\noindent
If $f\!:Z\!\lra\!X$ is a smooth map transverse to~$Y$,
then $f^{-1}(Y)\!\subset\!Z$ is a smooth submanifold with the normal bundle 
$\cN_Z(f^{-1}(Y)\!)$ isomorphic to $f^*\cN_XY$ via~$\nd f$.
An auxiliary structure on the normal neighborhood of~$Y$ in~$X$ needed
for a complex or augmented blowup pulls back to such a structure on
the normal neighborhood of \hbox{$f^{-1}(Y)\!\subset\!Z$}.
In each of the three settings of~\eref{blowdownmaps_e}, $f$ thus lifts to a smooth~map 
so that the corresponding diagram 
$$\xymatrix{\BLR_{f^{-1}(Y)}Z \ar[d]^{\pi}\ar[r]^>>>>>>{\BLR_Y\!f}& \BLR_YX\ar[d]^{\pi} &
\BLC_{f^{-1}(Y)}Z \ar[d]^{\pi}\ar[r]^>>>>>{\BLC_Y\!f}& \BLC_YX\ar[d]^{\pi}&
\BLau_{f^{-1}(Y)}Z \ar[d]^{\pi}\ar[r]^>>>>>{\BLau_Y\!f}& \BLau_YX\ar[d]^{\pi}\\
Z \ar[r]^f& X& Z \ar[r]^f& X & Z \ar[r]^f& X}$$
commutes.
The resulting lift~$\wt{f}$ of~$f$ is transverse to the exceptional locus in the real and complex cases
and to its components~$\bE^0_YX$ and~$\bE^-_YX$ in the augmented case.\\

\noindent
We say that smooth submanifolds $Y,Y'\!\subset\!X$ \sf{intersect cleanly} if 
$Y\!\cap\!Y'\!\subset\!X$  is also a submanifold~and 
\BE{PrInter_e}T(Y\!\cap\!Y')=TY|_{Y\cap Y'}\!\cap\!TY'|_{Y\cap Y'}\subset TX|_{Y\cap Y'}.\EE
This is equivalent to the condition that the homomorphism
$$\cN_{Y'}(Y\!\cap\!Y')\!\equiv\!\frac{TY'|_{Y\cap Y'}}{T(Y\!\cap\!Y')}
\lra \frac{TX|_{Y\cap Y'}}{TY|_{Y\cap Y'}}
\!\equiv\!\cN_XY|_{Y\cap Y'}$$
induced by the inclusions is injective.
If $Y$ contains~$Y'$ or intersects~it transversely, then
$Y$ and~$Y'$ intersect cleanly, but these are just extreme examples of clean intersection.\\

\noindent
If $\pi\!:\wt{X}\!\lra\!X$ is a blowup along a closed submanifold $Y\!\subset\!X$
with exceptional locus~$\bE$
and $Y'\!\subset\!X$ is a submanifold intersecting~$Y$ cleanly,
we call a submanifold $\wt{Y}'\!\subset\!\wt{X}$ \sf{$\pi$-equivalent} to~$Y'$
if $\pi|_{\wt{Y}'}$ is a bijection onto~$Y'$, 
$\wt{Y}'$ is disjoint from $\bE^0_YX\!\cap\!\bE^-_YX$ if 
$\pi$ is an augmented blowup and intersects~$\bE$ cleanly in all three cases, and 
$$T_{\wt{x}}\wt{Y}'\cap\ker\nd_{\wt{x}}\big\{\pi|_{\bE}\big\}=\{0\}
\quad\forall\,\wt{x}\!\in\!\bE\!\cap\!\wt{Y}'.$$
These conditions imply that $\pi|_{\wt{Y}'}$ is a diffeomorphism onto~$Y'$.
The next statement follows from the definitions of the three blowups 
through direct local considerations.

\begin{lmm}\label{BlTrans_lmm}
Suppose \hbox{$\pi\!:\wt{X}\!\lra\!X$} is a blowup 
along a closed codimension~$\fc$ submanifold~$Y$ as in~\eref{blowdownmaps_e}
with exceptional locus~$\bE$,
$Y'\!\subset\!X$ is a submanifold intersecting~$Y$ cleanly,
and $\wt{Y}'\!\subset\!\wt{X}$ is a submanifold  $\pi$-equivalent to~$Y'$.
If $f\!:Z\!\lra\!X$ is a smooth map transverse to~$Y$, $Y'$, and~$Y\!\cap\!Y'$,
then its lift~$\wt{f}$ to~$\wt{X}$ is transverse to~$\bE$, $\wt{Y}'$, 
and~$\bE\!\cap\!\wt{Y}'$.
\end{lmm}

\subsection{Blowups and (co)homology}
\label{BlHom_subs}

\noindent
Let $X$ be a smooth manifold.
For a continuous map $f\!:Z\!\lra\!X$,  let 
$$\Om(f)=\bigcap_{K\subset Z\text{~cmpt}}\!\!\!\!\!\!\!\!\ov{f(Z\!-\!K)}$$
be \sf{the limit set of~$f$}.
A \sf{pseudocycle} into~$X$ is a smooth map \hbox{$f\!:Z\!\lra\!X$} 
from an oriented manifold, without boundary, 
so that the closure of $f(Z)$ in~$X$ is compact and
there exists a smooth map $h\!:Z'\!\lra\!X$ such~that 
\BE{fhcond_e}\dim\,Z'\le \dim\,Z\!-\!2  \qquad\hbox{and}\qquad \Om(f)\subset h(Z')\,.\EE
The \sf{dimension} of such a pseudocycle~$f$ is $\dim\,Z$.
We denote by~$-f$ the pseudocycle with the opposite orientation on the domain~$Z$ of~$f$.\\

\noindent
A smooth map \hbox{$\wt{f}\!:\wt{Z}\!\lra\!X$} from an oriented manifold, 
possibly with boundary, 
is a \sf{bordered pseudocycle with boundary} $f\!:Z\!\lra\!X$ if
the closure of $\wt{f}(\wt{Z})$ in~$X$ is compact, 
$$\prt\wt{Z}=Z, \qquad \wt{f}|_Z=f,$$
and there exists a smooth map $\wt{h}\!:\wt{Z}'\!\lra\!X$ such~that 
$$\dim\,\wt{Z}'\le \dim\,\wt{Z}\!-\!2\qquad\hbox{and}\qquad
\Om(\wt{f})\subset\wt{h}(\wt{Z}') \,.$$
An \sf{equivalence} between two pseudocycles 
\hbox{$f_1\!:Z_1\!\lra\!X$} and \hbox{$f_2\!:Z_2\!\lra\!X$}  
is a bordered pseudocycle  \hbox{$\wt{f}\!:\wt{Z}\!\lra\!X$}
with boundary $(-f_1)\!\sqcup\!f_2$.
The set~$\PS_k(X)$ of the equivalence classes of dimension~$k$ pseudocycles
forms an abelian group under the disjoint union operation.
By \cite[Theorem~1.1]{pseudo}, this abelian group is naturally isomorphic to~$H_k(X;\Z)$.\\

\noindent
A pseudocycle $f\!:Z\!\lra\!X$ is \sf{transverse} to a smooth submanifold \hbox{$Y\!\subset\!X$}
if $f$ is transverse to~$Y$ 
and there exists a smooth map $h\!:Z'\!\lra\!X$ as above which is also transverse to~$Y$.
A bordered pseudocycle $\wt{f}\!:\wt{Z}\!\lra\!X$ is \sf{transverse} to~$Y$
if $\wt{f}$ and $\wt{f}|_{\prt\wt{Z}}$ are transverse to~$Y$ and 
there exists a smooth map $\wt{h}\!:\wt{Z}'\!\lra\!X$ as above 
which is also transverse to~$Y$.
If $Y'\!\subset\!X$ is another submanifold so that $Y\!\cap\!Y'\!\subset\!X$ is a smooth submanifold
as~well, 
we denote~by $\PS_k^{Y,Y'}(X)$ the set of the equivalence classes of 
the dimension~$k$ pseudocycles into~$X$ transverse to $Y$, $Y'$, and~$Y\!\cap\!Y'$ 
modulo the pseudocycle equivalences transverse to~$Y$, $Y'$, and~$Y\!\cap\!Y'$. 
This set forms an abelian group under the disjoint union operation.\\

\noindent
Let $Y,Y'\!\subset\!X$ be smooth submanifolds intersecting cleanly.
By the proof of \cite[Theorem~1.1]{pseudo}, a generic pseudocycle representing a given
homology class is then transverse to~$Y$, $Y'$, and~$Y\!\cap\!Y'$.
Furthermore, an equivalence between two equivalent pseudocycles that are transverse 
to~$Y$, $Y'$, and~$Y\!\cap\!Y'$ can be chosen to be also transverse to~$Y$, $Y'$, and~$Y\!\cap\!Y'$.
This yields Lemma~\ref{pseudo_lmm} below, which refines \cite[Theorem~1.1]{pseudo}.
The isomorphism of this lemma is the composition of the natural homomorphisms
$$\PS_k^{Y,Y'}(X) \lra \PS_k(X) \lra H_k(X;\Z).$$

\begin{lmm}\label{pseudo_lmm}
If $Y,Y'\!\subset\!X$ are closed submanifolds intersecting cleanly,
then the group $\PS_k^{Y,Y'}(X)$ is naturally isomorphic to $H_k(X;\Z)$. 
\end{lmm}

\noindent
Suppose $R$ is a commutative ring with unity and
$Y\!\subset\!X$ is a closed codimension~$\fc$ submanifold.
If $f\!:Z\!\lra\!X$ is a dimension~$k$ pseudocycle and
$h$ is a smooth map as in~\eref{fhcond_e} that are transverse to~$Y$,
then \eref{fhcond_e} holds with~$(f,h)$ replaced by $(f|_{f^{-1}(Y)},h|_{h^{-1}(Y)})$.
If the normal bundle~$\cN_XY$ of~$X$ in~$Y$ is oriented, then $f^{-1}(Y)$ inherits 
an orientation from~$Z$.
The map $f|_{f^{-1}(Y)}$ is then a pseudocycle into~$Y$ of dimension $k\!-\!\fc$ and 
determines a homology class in~$Y$.
By Lemma~\ref{pseudo_lmm} with $Y'\!=\!\eset$, this homology class is independent of the choice 
of admissible~$f$ in its homology class~$[f]$ in~$X$.
Thus, we obtain a \sf{restriction homomorphism}
\BE{Ycapdfn_e}Y\!\cap\!: H_*(X;R)\lra H_{*-\fc}(Y;R)\EE
if $\cN_XY$ is oriented or $\cha(R)\!=\!2$.
If $X$ is oriented or $\cha(R)\!=\!2$, 
this homomorphism equals \hbox{$\PD_Y\!\circ\!i^*\!\circ\!\PD_X$},
where $i\!:Y\!\lra\!X$ is the inclusion, and thus respects the ring structures.
If $X$ is also compact, the surjectivity of~\eref{Ycapdfn_e} is equivalent 
to the injectivity of the homomorphism $i_*\!:H_*(Y;R)\!\lra\!H_*(X;R)$.
We call the submanifold $Y\!\subset\!X$ \sf{$R$-incompressible} if
the homomorphism~\eref{Ycapdfn_e} is surjective.\\

\noindent
We call a collection $\cA$ of compact submanifolds of $X$ \sf{$R$-complete} if
\begin{enumerate}[label=$\bullet$,leftmargin=*]

\item either $X$ and all elements of~$\cA$ are oriented or $\cha(R)\!=\!2$ and 

\item the elements of~$\cA$ generate $H_*(X;R)$ as a ring.

\end{enumerate}
For a smooth map $\pi\!:\wt{X}\!\lra\!X$,
we call collections~$\cA$ and~$\wt\cA$ of compact submanifolds of~$X$ and~$\wt{X}$,
respectively, \sf{$(\pi,R)$-related} if
\begin{enumerate}[label=$\bullet$,leftmargin=*]

\item either all elements of~$\cA$ and of~$\wt\cA$
are oriented or $\cha(R)\!=\!2$ and 

\item $\{\pi_*([\wt{Y}'])\!\in\!H_*(X;R)\!:\wt{Y}'\!\in\!\wt\cA\}\!=\!
\{[Y']\!\in\!H_*(X;R)\!:Y'\!\in\!\cA\}$.

\end{enumerate}
The next statement collects the topological properties of 
the blowups~\eref{blowdownmaps_e} needed to establish 
the generation properties~\ref{HCgen_it} and~\ref{HRgen_it}.

\begin{prp}\label{BlHom_prp}
Let $R$ be a commutative ring with unity.
Suppose \hbox{$\pi\!:\wt{X}\!\lra\!X$} is a blowup 
along a closed codimension~$\fc$ submanifold~$Y$ as in~\eref{blowdownmaps_e},
$\cA$ is an $R$-complete collection of submanifolds of~$X$,
$\wt\cA$ is a collection of submanifolds of~$\wt{X}$ $(\pi,R)$-related to~$\cA$,
$Y'\!\subset\!X$ is an $R$-incompressible submanifold intersecting~$Y$ cleanly,
and $\wt{Y}'\!\subset\!\wt{X}$ is 
a submanifold  $\pi$-equivalent to~$Y'$.

\begin{enumerate}[label=(\arabic*),leftmargin=*]

\setcounter{enumi}{-1}

\item\label{BlHomR_it} 
If $\pi$ is a real blowup,
2 is a unit in~$R$, and $\fc\!\not\in\!2\Z$,
then the collection~$\wt\cA$ is $R$-complete 
and the submanifold $\wt{Y}'\!\subset\!\wt{X}$ is $R$-incompressible.

\end{enumerate}
Suppose in addition $Y\!\subset\!X$ is compact and $R$-incompressible.
\begin{enumerate}[label=(\arabic*),leftmargin=*]

\item\label{BlHomR_it2} 
If $\pi$ is a real blowup and $\cha(R)\!=\!2$, then the collection
$\wt\cA\!\cup\!\{\bE^{\R}_Y\!X\}$ is $R$-complete 
and the submanifold $\wt{Y}'\!\subset\!\wt{X}$ is $R$-incompressible.

\item\label{BlHomC_it} If $\pi$ is a complex blowup, then
the collection $\wt\cA\!\cup\!\{\bE_Y^{\C}\!X\}$ is $R$-complete 
and the submanifold $\wt{Y}'\!\subset\!\wt{X}$ is $R$-incompressible.

\item\label{BlHomRR_it} 
If $\pi$ is a $\fc_1\!=\!1$ augmented blowup
and either $\cha(R)\!=\!2$ or $\fc\!\not\in\!2\Z$ and the vector bundles 
$\cN_XY$ and $\cN^{\fc_1}_YX$ are oriented,
then 
the collection $\wt\cA\!\cup\!\{\bE_Y^0\!X,\bE_Y^-\!X\}$ is $R$-complete 
and the submanifold $\wt{Y}'\!\subset\!\wt{X}$ is $R$-incompressible.

\end{enumerate}
\end{prp}

\noindent
Let $\pi\!:\bE\!\lra\!Y$ be a smooth fiber bundle with compact oriented $m$-dimensional fibers.
If $f\!:Z\!\lra\!Y$ is a dimension~$k$ pseudocycle, then 
$$\bE|_f\!:f^*\bE\lra \bE, \qquad \bE|_f(z,e)=e,$$
is a $(k\!+\!m)$-pseudocycle and determines a homology class~$[\bE|_f]$ on~$\bE$.
This homology class is independent of the choice 
of~$f$ in the homology class~$[f]$ in~$Y$ determined by~$f$.
Thus, we obtain a homomorphism
$$\bE|_{\cdot}\!:H_*(Y;R)\lra H_{*+m}(\bE;R)$$
if the fibers of $\pi$ are (continuously) oriented or $\cha(R)\!=\!2$.
If $Y$ is oriented  or $\cha(R)\!=\!2$, 
this homomorphism equals $\PD_{\bE}\!\circ\!\pi^*\!\circ\!\PD_Y$.
If $\pi'\!:\bE'\!\lra\!Y$ is another continuous map such that
$\bE\!\subset\!\bE'$ and $\pi'|_{\bE}\!=\!\pi$, 
we denote the composition of the homomorphism~$\bE|_{\cdot}$
with the natural homomorphism \hbox{$H_*(\bE;R)\!\lra\!H_*(\bE';R)$}
also by~$\bE|_{\cdot}$.

\begin{proof}[Proof of Proposition~\ref{BlHom_prp}]
If $f\!:Z\!\lra\!X$ and $h'\!:Z\!\lra\!X$ are smooth map as in~\eref{fhcond_e} 
which are transverse to~$Y$,
then \eref{fhcond_e} holds with~$(f,h)$ replaced by 
$(\BLR_Y\!f,\BLR_Y\!h)$, $(\BLC_Y\!f,\BLC_Y\!h)$, and $(\BLau_Y\!f,\BLau_Y\!h)$.
The domain of~$\BLC_Y\!f$ inherits an orientation from that of~$Z$.
The same is the case for $\BLR_Y\!f$ and $\BLau_Y\!f$
if the rank of the normal bundle~$\cN_XY$ of~$Y$ in~$X$ is odd.
Thus, $\BLR_Y\!f$, $\BLC_Y\!f$, and $\BLau_Y\!f$ are dimension~$k$ pseudocycles 
if $f$ is and the appropriate orientability conditions for the relevant blowup hold.
In such a case, they determine homology classes in the corresponding blowups of~$X$.
Along with Lemmas~\ref{pseudo_lmm} and~\ref{BlTrans_lmm}, 
this implies the $R$-incompressible claims of the proposition.\\

\noindent
It remains to establish the $R$-complete claims of the proposition.
We denote~by
$$\wt{i}\!:\bE\lra\wt{X}$$
the inclusion of the exceptional locus.
In light of Lemma~\ref{pseudo_lmm} with $Y'\!=\!\eset$, 
the lifted homology classes above are independent of the choice 
of admissible~$f$ in its homology class~$[f]$ in~$X$.
Thus, we obtain a \sf{blowup homomorphism}
\BE{BLdfn_e}\BL\!:H_*(X;R)\lra H_*(\wt{X};R) \qquad\hbox{s.t.}\quad
\pi_*\!\circ\!\BL\!=\!\Id\!:H_*(X;R)\lra H_*(X;R),\EE
whenever the appropriate orientability conditions hold or $\cha(R)\!=\!2$. 
If $X$ is oriented or $\cha(R)\!=\!2$,
this homomorphism is the composition $\PD_{\wt{X}}\!\circ\!\pi^*\!\circ\!\PD_X$
and thus respects the ring structures.
In the real and complex cases,
\BE{BLrestr_e1} \wt{i}_*\big(\bE|_{Y\cap \mu'}\big)=
[\bE]\!\cdot\!\BL(\mu') \qquad\forall~\mu'\!\in\!H_*(X;R),\EE
whenever the appropriate orientability conditions hold or $\cha(R)\!=\!2$. 
In the augmented case,
\BE{BLrestr_e2} \wt{i}_*\big(\bE_Y^0X|_{Y\cap \mu'}\big)=[\bE_Y^0X]\!\cdot\!\BL(\mu'),
~~
\wt{i}_*\big(\bE_Y^-X|_{Y\cap \mu'}\big)=[\bE_Y^-X]\!\cdot\!\BL(\mu') 
\quad\forall~\mu'\!\in\!H_*(X;R)\EE
under the assumptions in Claim~\ref{BlHomRR_it}.\\

\noindent
We omit the coefficient ring~$R$ below.
Let $U\!\subset\!X$ of~$Y$ be a tubular neighborhood 
and $\wt{U}\!\subset\!\wt{X}$ be its preimage under~$\pi$. 
Since $\bE$ is a deformation retract of~$\wt{U}$, the Mayer-Vietoris sequences
for $\wt{X}\!=\!(\wt{X}\!-\!\bE)\!\cup\!\wt{U}$ and $X\!=\!(X\!-\!Y)\!\cup\!U$
induce a commutative diagram 
$$\xymatrix{\ldots\ar[r]& H_*(\wt{U}\!-\!\bE)\ar[r]\ar[d]^{\id} &
H_*(\wt{X}\!-\!\bE)\!\oplus\!H_*(\bE) \ar[r]\ar[d]|{\id\oplus\pi_*}&
H_*(\wt{X}) \ar[r]\ar[d]^{\pi_*}&  H_{*-1}(\wt{U}\!-\!\bE)\ar[d]^{\id}\ar[r]&\ldots\\
\ldots\ar[r]& H_*(U\!-\!Y)\ar[r]& H_*(X\!-\!Y)\!\oplus\!H_*(Y) \ar[r]&
H_*(X) \ar[r]&  H_{*-1}(U\!-\!Y)\ar[r]&\ldots}$$
of exact sequences.
It gives
\BE{BlHom_e5}\ker\!\big(\pi_*\!: H_*(\wt{X})\!\lra\!H_*(X)\!\big)=
\wt{i}_*\!\big(\!\ker\!\big(\pi_*\!:H_*(\bE)\!\lra\!H_*(Y)\!\big)\!\big).\EE 
Along with~\eref{BLdfn_e}, this implies that it is sufficient to show that
the right-hand side of~\eref{BlHom_e5} is contained in the ring generated
by the image of~$\BL$ in the setting of~\ref{BlHomR_it},
by the image of~$\BL$ and~$[\bE]$ in the settings of~\ref{BlHomR_it2}
and~\ref{BlHomC_it},
and by the image of~$\BL$, $[\bE_Y^0\!X]$, and~$[\bE_Y^-\!X]$
in the setting of~\ref{BlHomRR_it}.\\

\noindent
Claim~\ref{BlHomR_it} follows from~\eref{BlHom_e5} and 
the first statement of Lemma~\ref{BlHom_lmm}\ref{BlHom2R_it} below.
By~\eref{cNEisom_e}, 
the normal bundle of $\bE\!=\!\R\P(\cN_XY)$ (resp.~$\bE\!=\!\C\P(\cN_XY)$)
in~$\wt{X}$ in the setting of Claim~\ref{BlHomR_it2} (resp.~\ref{BlHomC_it})
is the tautological line bundle~$\ga_{\cN_XY}^{\R}$ (resp.~$\ga_{\cN_XY}^{\C}$).
The second statement of Lemma~\ref{BlHom_lmm}\ref{BlHom2R_it},
the $R$-incompressibility of~$Y\!\subset\!X$, and~\eref{BLrestr_e1} then give
\begin{equation*}\begin{split}
\wt{i}_*\!\big(\!\ker\!\big(\pi_*\!:H_*(\bE)\!\lra\!H_*(Y)\!\big)\!\big)
&=\Big\{\sum_{r=1}^{\fc-1}\wt{i}_*\big(\!
\big(w_1(\cN_{\wt{X}}\bE)\!\big)^{\fc-1-r}\!\cap\!\big(\bE|_{\mu_r}\big)\!\big)\!:
\mu_r\!\in\!H_{*-r}(Y)\!\Big\}\\
&=\Big\{\sum_{r=1}^{\fc-1}[\bE]^{\fc-r}\!\cdot\!\BL(\mu_r')\!:
\mu_r'\!\in\!H_{*-r+\fc}(X)\!\Big\}
\subset H_*(\wt{X}).
\end{split}\end{equation*}
Along with \eref{BlHom_e5}, this implies Claim~\ref{BlHomR_it2}. 
Claims~\ref{BlHomC_it} and~\ref{BlHomRR_it} are obtained similarly 
from~\ref{BlHom2C_it} and~\ref{BlHom2RR_it}, respectively, in Lemma~\ref{BlHom_lmm}
(without the need to consider $\cN_{\wt{X}}\bE$ in the last case),
along with~\eref{BLrestr_e1} and~\eref{BLrestr_e2}.
\end{proof}

\noindent
For a real vector bundle $\cN\!\lra\!Y$ of rank~$\fc\!\ge\!2$ and
a corank~1 subbundle $\cN^1\!\subset\!\cN$, let
\BE{bEcN1mdfn_e}\bE_{\cN^1}^-\cN\lra Y\EE
be the $S^{\fc-1}$-fiber bundle obtained by collapsing 
each fiber of the $\R\P^{\fc-2}$-fiber subbundle \hbox{$\R\P\cN^1\!\subset\!\R\P\cN$}
to a point.
Let
$$\bE_{\cN^1}\cN=
\R\P\big(\cN/\cN^1\!\oplus\!\tau_Y^1\big)\!\cup\!\bE_{\cN^1}^-\cN$$
be the topological space obtained by identifying $\R\P(\cN/\cN^1\!\oplus\!\tau_Y^0)$
with the section~$\bE_{\cN^1}^{-0}\cN$
of~\eref{bEcN1mdfn_e} obtained from the above collapse.

\begin{lmm}\label{BlHom_lmm}
Let $\cN\!\lra\!Y$ be a real vector bundle of rank~$\fc$.
\begin{enumerate}[label=(\arabic*),leftmargin=*]

\item\label{BlHom2R_it}
If 2 is a unit in~$R$ and $\fc\!\not\in\!2\Z$, then the homomorphism 
$$\pi_*\!:H_*(\R\P\cN;R)\lra H_*(Y;R)$$
is an isomorphism.
If $\cha(R)\!=\!2$, then
$$\ker\pi_*=\Big\{\sum_{r=1}^{\fc-1}\!
\big(w_1(\ga_{\cN}^{\R})\!\big)^{\fc-1-r}\!\cap\!\big(\R\P\cN|_{\mu_r}\big)
\!:\mu_r\!\in\!H_{*-r}(Y;R)\!\Big\}.$$

\item\label{BlHom2C_it} If $\cN$ is a complex vector bundle, then
$${}\!\!\!\!\!\!\!\!\!\ker\!\big(\pi_*\!:H_*(\C\P\cN;R)\!\lra\!H_*(Y;R)\!\big)
=\Big\{\sum_{r=1}^{\fc/2-1}\!\!\!
\big(c_1(\ga_{\cN}^{\C})\!\big)^{\fc/2-r}\!\cap\!\big(\C\P\cN|_{\mu_r}\big)
\!:\mu_r\!\in\!H_{*-2r}(Y;R)\!\Big\}.$$

\item\label{BlHom2RR_it} If $\fc\!\ge\!2$,
$\cN^1\!\subset\!\cN$ is a vector bundle of corank~1,
and either $\cha(R)\!=\!2$ or 
the vector bundles~$\cN$ and~$\cN^1$ are oriented, then
\begin{equation*}\begin{split}
&{}\!\!\!\ker\!\big(\pi_*\!:H_*(\bE_{\cN^1}\cN;R)\!\lra\!H_*(Y;R)\!\big)\\
&\hspace{.5in}
=\big\{\R\P(\cN/\cN^1\!\oplus\!\tau_Y^1)|_{\mu_1}\!+\!
\bE_{\cN^1}^-\cN|_{\mu_2}\!:
\mu_1\!\in\!H_{*-1}(Y;R),\,\mu_2\!\in\!H_{*-\fc+1}(Y;R)\!\big\}.
\end{split}\end{equation*}
  
\end{enumerate}
\end{lmm}

\begin{proof}
We denote by $\bE$ the bundle $\R\P\cN$, $\P\cN$, or $\bE_{\cN^1}\cN$ over~$Y$,
depending on the case.
In the setting of~\ref{BlHom2RR_it}, the orientations of~$\cN$ and~$\cN^1$
determine orientation classes 
\begin{equation*}\begin{split}
u_1&\in H^1\big(\bE,\bE_{\cN^1}^-\cN;R\big)
\!=\!H^1\big(\R\P(\cN/\cN^1\!\oplus\!\tau_Y^1),
\R\P(\cN/\cN^1\!\oplus\!\tau_Y^0);R\big) 
\qquad\hbox{and}\\
u_2&\in H^{\fc-1}\big(\bE,\R\P(\cN/\cN^1\!\oplus\!\tau_Y^1);R\big)
\!=\!H^{\fc-1}\big(\bE_{\cN^1}^-\cN,\bE_{\cN^1}^{-0}\cN;R\big).
\end{split}\end{equation*}
We denote their restrictions to $\bE$ by $u_1'$ and $u_2'$, respectively.\\

\noindent
For every $y\!\in\!Y$, the collections  
$$\big\{1|_{\bE_y}\big\} \quad\hbox{and}\quad
\big\{1|_{\bE_y},w_1(\ga_{\cN}^{\R})|_{\bE_y},\ldots,w_1(\ga_{\cN}^{\R})^{\fc-1}|_{\bE_y}\big\}$$
are bases for~$H^*(\bE_y;R)$ in the two cases of~\ref{BlHom2R_it}.
By \cite[Theorem~5.7.9]{Spanier}, the corresponding homomorphisms
\begin{alignat*}{2}
H_*(\bE;R)&\lra H_*(Y;R), &\qquad \eta&\lra\pi_*(\eta),\\
H_*(\bE;R)&\lra \bigoplus_{r=0}^{\fc-1}\!H_{*-r}(Y;R),
&\qquad
\eta&\lra\big(\pi_*\big(\!(w_1(\ga_{\cN}^{\R})\!)^r\!\cap\eta\big)_{r=0,\ldots,\fc-1}\big),
\end{alignat*}
are thus isomorphisms.
The collections
$$\big\{1|_{\bE_y},c_1(\ga_{\cN}^{\C})|_{\bE_y},\ldots,
c_1(\ga_{\cN}^{\C})^{\fc/2-1}|_{\bE_y}\big\}
\quad\hbox{and}\quad
\big\{1|_{\bE_y},u_1'|_{\bE_y},u_2'|_{\bE_y}\big\}$$
are bases for~$H^*(\bE_y;R)$ in the settings of~\ref{BlHom2C_it} and~\ref{BlHom2RR_it}, respectively.
By \cite[Theorem~5.7.9]{Spanier}, the corresponding homomorphisms
\begin{alignat*}{2}
H_*(\bE;R)&\lra \bigoplus_{r=0}^{\fc/2-1}\!\!H_{*-2r}(Y;R),
&\quad
\eta&\lra\big(\pi_*\big(\!(c_1(\ga_{\cN}^{\C})\!)^r\!\cap\eta\big)_{r=0,\ldots,\fc/2-1}\big),\\
H_*(\bE;R)&\lra H_*(Y;R)\!\oplus\!H_{*-1}(Y;R)\!\oplus\!H_{*-\fc+1}(Y;R),
&\quad
\eta&\lra\big(\pi_*(\eta),\pi_*\big(u_1'\!\cap\!\eta\big),\pi_*\big(u_2'\!\cap\!\eta\big)\!\big),
\end{alignat*}
are thus isomorphisms.
This implies the claims.
\end{proof}

\subsection{The complex case}
\label{GenC_subs}

\noindent
Let $\ell\!\in\!\Z^+$ with $\ell\!\ge\!3$.
For the purposes of establishing the generation property~\ref{HCgen_it} 
on page~\pageref{HCgen_it}, we index the divisors~$D_{I,J}$
in~$\ov\cM_{\ell}$ and~$\ov\cM_{\ell+1}$ by the~sets
$$\cA_{\ell}\equiv\big\{\vr\!\subset\![\ell]\!: \big|\vr\!\cap\![3]\!\big|\!\ge\!2,~
\big|[\ell]\!-\!\vr\big|\!\ge\!2\big\}
\quad\hbox{and}\quad
\wt\cA_{\ell}\equiv\cA_{\ell}\!\sqcup\!
\big\{\!\big[\ell\big]\!-\!\{i\}\!:i\!\in\![\ell]\!\big\}.$$
For each $\vr\!\in\!\cA_{\ell}$ and $i\!\in\![\ell]$, let 
\begin{gather*}
D_{\vr}\equiv D_{\vr,[\ell]-\vr}\subset\ov\cM_{\ell}, \quad
\wt{D}_{[\ell]-\{i\}}^0\equiv D_{[\ell]-\{i\},\{i,\ell+1\}}\subset\ov\cM_{\ell+1},
\quad \wt{D}_{[\ell]-\{i\}}^+=\eset,\\
\wt{D}_{\vr}^+= D_{\vr\cup\{\ell+1\},[\ell]-\vr}\subset\ov\cM_{\ell+1},
\qquad
\wt{D}_{\vr}^0\equiv D_{\vr,([\ell]-\vr)\cup\{\ell+1\}}\subset\ov\cM_{\ell+1}\,.
\end{gather*}
The set~$\cA_{\ell}$ is partially ordered by the inclusion~$\subsetneq$ of subsets of~$[\ell]$.
We extend this partial order to a strict order~$<$ on $\{0\}\!\sqcup\!\cA_{\ell}$
so that~0 is the smallest element and define \hbox{$\vr\!-\!1\!\in\!\{0\}\!\cup\!\cA_{\ell}$}
to be the predecessor of $\vr\!\in\!\cA_{\ell}$.
For \hbox{$\vr^*\!\in\!\{0\}\!\sqcup\!\cA_{\ell}$}, define 
$$\cA_{\ell}(\vr^*)=\big\{\vr\!\in\!\cA_{\ell}\!:\vr\!>\!\vr^*\big\}.$$

\vspace{.15in}

\noindent 
For $\vr\!\in\!\cA_{\ell}$ and $\wt\cC_1,\wt\cC_2\!\in\!\ov\cM_{\ell+1}$, 
we define $\wt\cC_1\!\sim_{\vr}\!\wt\cC_2$ if either $\wt\cC_1\!=\!\wt\cC_2$ or 
$$\ff_{\ell+1}(\wt\cC_1)=\ff_{\ell+1}(\wt\cC_2)\in \ov\cM_{\ell}
\quad\hbox{and}\quad \wt\cC_1,\wt\cC_2\in \wt{D}_{\vr}^0.$$
For every $\vr^*\!\in\!\{0\}\!\sqcup\!\cA_{\ell}$
the union of all equivalence relations~$\sim_{\vr}$ on~$\ov\cM_{\ell+1}$ with 
$\vr\!\in\!\cA_{\ell}(\vr^*)$ is again an equivalence relation;
see Section~1.1 in~\cite{RDMbl}.
The quotient~$X_{\vr^*}$ of~$\ov\cM_{\ell+1}$ by the last equivalence relation
is a compact Hausdorff space for every \hbox{$\vr^*\!\in\!\{0\}\!\sqcup\!\cA_{\ell}$}.
Let
\BE{pivrdfn_e}
q_{\vr^*}\!:\ov\cM_{\ell+1}\lra X_{\vr^*} \qquad\hbox{and}\qquad
\pi_{\vr^*}\!:X_{\vr^*}\lra X_{\vr^*-1}\EE
be the quotient projection and the map induced by~$q_{\vr^*-1}$
whenever $\vr^*\!\in\!\cA_{\ell}$, respectively.
For $\vr\!\in\!\wt\cA_{\ell}$, define
$$Y_{\vr^*;\vr}^+=q_{\vr^*}\big(\wt{D}_{\vr}^+\big)
\quad\hbox{and}\quad
Y_{\vr^*;\vr}^0= q_{\vr^*}\big(\wt{D}_{\vr}^0\big).$$
These subspaces of $X_{\vr^*}$ are compact.
By Theorem~2.1($\C1$) in~\cite{RDMbl},
$X_{\vr^*}$ is a complex manifold containing~$Y_{\vr^*;\vr}^+,Y_{\vr^*;\vr}^0$
as complex submanifolds.\\

\noindent
We assume that property~\ref{HCgen_it} holds as stated and
with~$\ell$ replaced by any smaller integer.
We will establish the following statements by induction 
on $\vr^*\!\in\!\{0\}\!\cup\!\cA_{\ell}$. 
\begin{enumerate}[label=($\C1\alph*$),leftmargin=*]

\item\label{CYgen_it} The submanifolds $Y_{\vr^*;[\ell]-\{1\}}^0\!\subset\!X_{\vr^*}$,
$Y_{\vr^*;\vr}^+$ with $\vr\!\in\!\cA_{\ell}$, 
and $Y_{\vr^*;\vr}^0$ with \hbox{$\vr\!\in\!\cA_{\ell}\!-\!\cA_{\ell}(\vr^*)$}
generate $H_*(X_{\vr^*};\Z)$ as a~ring. 

\item\label{CYTincomp_it}
For every $\vr\!\in\!\cA_{\ell}(\vr^*)$,
the submanifold $Y_{\vr^*;\vr}^0\!\subset\!X_{\vr^*}$ is $\Z$-incompressible.

\end{enumerate}
If $\vr^*\!\in\!\cA_{\ell}$ is the maximal element, 
\ref{CYgen_it} is a stronger version of property~\ref{HCgen_it}
with~$\ell$ replaced by~$\ell\!+\!1$.\\

\noindent
By Theorem~2.1($\C1$) in~\cite{RDMbl}, the map~\eref{CMlind_e} descends to a biholomorphism 
$$\Psi_0\!:X_0\lra \ov\cM_{\ell}\!\times\!\ov\cM_4\,.$$
We identify~$X_0$ with $\ov\cM_{\ell}\!\times\!\ov\cM_4$ via~$\Psi_0$.
Under this identification,
\begin{gather*}
Y_{0;\vr}^+\equiv\big\{\!(\ff_{\ell+1},\ff_{123,\ell+1})\!\big\}(\wt{D}_{\vr}^+)=
D_{\vr}\!\times\!\ov\cM_4 \quad\forall~\vr\!\in\!\cA_{\ell},\\
Y_{0;[\ell]-\{1\}}^0\equiv\big\{\!(\ff_{\ell+1},\ff_{123,\ell+1})\!\big\}
\big(\wt{D}_{[\ell]-\{1\}}^0\big)=\ov\cM_{\ell}\!\times\!D_{14,23}.
\end{gather*}
Since the divisors~$D_{\vr}$ with $\vr\!\in\!\cA_{\ell}$ generate $H_*(\ov\cM_{\ell};\Z)$ as a ring,
the submanifolds $Y_{0;[\ell]-\{1\}}^0,Y_{0;\vr}^+\!\subset\!X_0$ generate $H_*(X_0;\Z)$ as a ring.
This confirms the $\vr^*\!=\!0$ case of~\ref{CYgen_it}.\\

\noindent
For a finite set~$S$, we denote by $\ov\cM_S$ the Deligne-Mumford moduli space
of stable rational $S$-marked curves. 
If $J\!\sqcup\!K\!=\!S$, we define  the divisor 
\hbox{$D_{J,K}\!\subset\!\ov\cM_S$} as in Section~\ref{Cthm_subs}.
Let $\vr\!\in\!\cA_{\ell}$ and $\vr_{\ell}^c\!=\![\ell]\!-\!\vr$.
Since $|\vr_{\ell}^c\!\cap\![3]|\!\le\!1$,  the forgetful morphism
$$\ff_{\vr}\!\equiv\!\ff_{123j}\!:D_{\vr}\lra \ov\cM_4$$
is independent of the choice of $j\!\in\!\vr_{\ell}^c$ distinct from $1,2,3$.
If $[3]\!\subset\!\vr$, $\ff_{\vr}$ is the composition of the natural identification
\BE{DvrCiden_e}D_{\vr}\xlra{~\approx~} 
\ov\cM_{\{\nod\}\sqcup \vr}\!\times\!\ov\cM_{\{\nod\}\sqcup\vr_{\ell}^c}\EE
with the projection to the first component and the forgetful morphism sending 
the marked points $z_1$, $z_2$, $z_3$, and $z_0$ to $1,2,3,4$, respectively. 
Otherwise, $\ff_{\vr}$ is the constant map so that $\ff_{\vr}(D_{\vr})$
is the two-component curve in~$\ov\cM_4$ with one of the components carrying
the marked points corresponding to $\vr\!\cap\![3]$.\\

\noindent
For every $\vr\!\in\!\cA_{\ell}$, the~map
$$s_{\vr}\!=\!\big(\io_{\vr},\ff_{\vr}\big)\!:
D_{\vr}\lra \ov\cM_{\ell}\!\times\!\ov\cM_4\!=\!X_0,$$
where $\io_{\vr}\!:D_{\vr}\!\lra\!\ov\cM_{\ell}$ is the inclusion,
is a holomorphic embedding and
$$Y_{0;\vr}^0\equiv\big\{\!(\ff_{\ell+1},\ff_{123,\ell+1})\!\big\}(\wt{D}_{\vr}^0)=
s_{\vr}(D_{\vr}).$$
If in addition $\vr'\!\in\!\cA_{\ell}$,
\BE{CYvrvrcap_e}Y_{0;\vr}^0\!\cap\!Y_{0;\vr'}^+\stackrel{s_{\vr}}{\approx}
D_{\vr}\!\cap\!D_{\vr'}=\begin{cases}
D_{\vr',\{\nod\}\sqcup(\vr-\vr')}
\!\times\!\ov\cM_{\{\nod\}\sqcup \vr^c_{\ell}}, &\hbox{if}~\vr'\!\subsetneq\!\vr;\\
\ov\cM_{\{\nod\}\sqcup\vr}\!\times\!
D_{\{\nod\}\sqcup(\vr^c_{\ell}-\vr'^c_{\ell}),\vr'^c_{\ell}}, 
&\hbox{if}~\vr'\!\supsetneq\!\vr;
\end{cases}\EE
under the identification~\eref{DvrCiden_e}.
Since $Y_{0;\vr}^0$ and $Y_{0;\vr'}^+$ intersect transversely in~$X_0$ for $\vr'\!\neq\!\vr$,
the K\"unneth formula and the inductive assumption for property~\ref{HCgen_it}
thus imply that the restriction homomorphism
$$Y_{0;\vr}^0\cap\!: H_*(X_0;\Z)\lra H_{*-4}(Y_{0;\vr}^0;\Z)$$
is surjective, i.e.~$Y_{0;\vr}^0\!\subset\!X_0$ is $\Z$-incompressible.
This confirms the $\vr^*\!=\!0$ case of~\ref{CYTincomp_it}.\\

\noindent
Suppose $\vr^*\!\in\!\cA_{\ell}$ is such that~\ref{CYgen_it} and~\ref{CYTincomp_it} 
with~$\vr^*$ replaced by~$\vr^*\!-\!1$ hold.
By Theorem~2.1($\C2$) in~\cite{RDMbl}, 
the second map in~\eref{pivrdfn_e} 
is the holomorphic blowup of~$X_{\vr^*-1}$ along~$Y_{\vr^*-1;\vr^*}^0$
with the exceptional locus~$Y_{\vr^*;\vr^*}^0$.
By the first statement of Theorem~2.1($\C3$) in~\cite{RDMbl}, 
each submanifold $Y_{\vr^*;\vr}^{\bu}\!\subset\!X_{\vr^*}$ 
with $(\vr,\bu)\!\neq\!(\vr^*,0)$ is 
$(\pi_{\vr^*},\Z)$-related to~$Y_{\vr^*-1;\vr}^{\bu}$.
Along with~\ref{CYgen_it} with~$\vr^*$ replaced by~$\vr^*\!-\!1$,
\ref{CYTincomp_it} with~$(\vr^*,\vr)$ replaced by~$(\vr^*\!-\!1,\vr^*)$,
and the first claim of Proposition~\ref{BlHom_prp}\ref{BlHomC_it},
this implies that \ref{CYgen_it} as stated holds as~well.
By the second statement of Theorem~2.1($\C3$) in~\cite{RDMbl}, 
each submanifold \hbox{$Y_{\vr^*;\vr}^0\!\subset\!X_{\vr^*}$} with $\vr\!\in\!\cA_{\ell}(\vr^*)$
is $\pi_{\vr^*}$-equivalent to the submanifold 
\hbox{$Y_{\vr^*-1;\vr}^0\!\subset\!X_{\vr^*-1}$}.
Along with~\ref{CYTincomp_it} with~$(\vr^*,\vr)$ replaced by~$(\vr^*\!-\!1,\vr)$
and the second claim of Proposition~\ref{BlHom_prp}\ref{BlHomC_it},
this implies that \ref{CYTincomp_it} as stated holds as~well.

\subsection{The real case}
\label{GenR_subs}

\noindent
Let $\ell\!\in\!\Z^+$ with $\ell\!\ge\!2$ and
$$\cA_{\ell}^{\pm}=\big\{\vr\!\subset\![\ell^{\pm}]\!: \big|\vr\!\cap\!\{1^+,1^-,2^+\}\!\big|\!\ge\!2,~
\big|[\ell^{\pm}]\!-\!\vr\big|\!\ge\!2\big\}.$$
For a finite set $S$, we denote by $\R\ov\cM_{0,S}$ the Deligne-Mumford moduli space
of stable real rational curves with conjugate pairs of points marked by the~set
$$S^{\pm}\equiv\big\{s^+\!:s\!\in\!S\big\}\!\sqcup\!\big\{s^-\!:s\!\in\!S\big\}.$$
If $J\!\sqcup\!K\!=\!S$ and \hbox{$I'\!\sqcup\!J'\!\sqcup\!K'\!=\!S$}, we define
$$\R E_{J,K},\R H_{J,K}\subset\R\ov\cM_{0,S} \qquad\hbox{and}\qquad
\R D_{I';J',K'}\subset\R\ov\cM_{0,S}$$
as in Section~\ref{Rthm_subs}.
For any $\{J,K\}\!\in\!\cP(\ell)$ and \hbox{$(I',\{J',K'\})\!\in\!\wt\cP_{\bu}(\ell)$},
there are natural identifications
\BE{Rboundident_e}
\R E_{J,K}\approx \ov\cM_{\{\nod\}\sqcup[\ell]} \quad\hbox{and}\quad
\R D_{I';J',K'}\approx \R\ov\cM_{0,\{\nod\}\sqcup I'}\!\times\!
\ov\cM_{\{\nod\}\sqcup(J'\cup K')}.\EE
Thus, all boundary hypersurfaces $\R E_{J,K}$ and boundary divisors $\R D_{I';J',K'}$
of $\R\ov\cM_{0,\ell}$ (and of~$\R\ov\cM_{0,\ell+1}$) are connected and orientable.\\

\noindent
For $\{J,K\}\!\in\!\cP(\ell)$ and $(I',\{J',K'\})\!\in\!\wt\cP_{\bu}(\ell)$
with $J'\!\subset\!J$ and $K'\!\subset\!K$
(and thus \hbox{$I'\!=\!(J\!-\!J')\!\cup\!(K\!-\!K')$}),
\BE{REinter_e}\begin{split}
\R E_{J,K}\!\cap\!\R D_{I';J',K'}&\approx D_{\{\nod\}\sqcup I',J'\cup K'} \\
&\approx \R E_{\{\nod\}\sqcup (J-J'),K-K'}\!\times\!
\ov\cM_{\{\nod\}\sqcup(J'\cup K')}
\end{split}\EE
under the two identifications in~\eref{Rboundident_e}.
For $(I,\{J,K\}),(I',\{J',K'\})\!\in\!\wt\cP_{\bu}(\ell)$ distinct,
\BE{RDinter_e}\R D_{I';J',K'}\!\cap\!\R D_{I;J,K}=
\begin{cases}\R\ov\cM_{0,\{\nod\}\sqcup I'}\!\times\!
D_{\{\nod\}\sqcup(J'\cup K'-J\cup K),J\cup K},
&\hbox{if}~\{J,K\}\!\prec\!\{J',K'\};\\
\R D_{I;\{\nod\}\sqcup(J-J'),K-K'}\!\times\!
\ov\cM_{\{\nod\}\sqcup(J'\cup K')},
&\hbox{if}~J'\!\subset\!J,\,K'\!\subset\!K;\\
\R  D_{\{\nod\}\sqcup(I'-J\cup K);J,K}\!\times\!
\ov\cM_{\{\nod\}\sqcup(J'\cup K')},
&\hbox{if}~J,K\!\subset\!I';
\end{cases}\EE
under the second identification in~\eref{Rboundident_e}.
For $\{J,K\}\!\in\!\cP(\ell)$ 
and $(I',\{J',K'\})\!\in\!\wt\cP_{\bu}(\ell)$ with $J,K\!\neq\!\eset$
and $K\!\subset\!I'$,
\BE{RDinter_e2}\R D_{I';J',K'}\!\cap\!\R H_{J,K}\approx
\R H_{\{\nod\}\sqcup(I'-K),K}\!\times\! \ov\cM_{\{\nod\}\sqcup(J'\cup K')}\EE
under the second identification in~\eref{Rboundident_e}.
The intersections of~$\R E_{J,K}$ and $\R D_{I';J',K'}$ with other
boundary hypersurfaces and divisors are empty.\\

\noindent
For $\vr\!\subset\![\ell^{\pm}]$, define
$$\vr_{\ell^{\pm}}^c=[\ell^{\pm}]\!-\!\vr, \qquad 
\ov\vr=\big\{i^+\!:i^-\!\in\!\vr\big\}\!\cup\!\big\{i^-\!:i^+\!\in\!\vr\big\},
\quad\hbox{and}\quad 
\vr\!\cap\![\ell]=\big\{i\!\in\![\ell]\!:i^+\!\in\!\vr\big\}.$$
For the purposes of establishing the generation property~\ref{HRgen_it} 
on page~\pageref{HRgen_it}, 
we index the boundary hypersurfaces and divisors of~$\R\ov\cM_{0,\ell}$ 
and~$\R\ov\cM_{0,\ell+1}$ by the~sets
\begin{gather*}
\begin{aligned}
\cA_{\ell;1}^D&=
\big\{\vr\!\in\!\cA_{\ell}^{\pm}\!:\ov\vr\!\subsetneq\!\vr_{\ell^{\pm}}^c\big\}, 
&\quad
\cA_{\ell;2}^D&=
\big\{\vr\!\in\!\cA_{\ell}^{\pm}\!: \ov\vr\!\supsetneq\!\vr_{\ell^{\pm}}^c\big\},
&\quad
\cA_{\ell}^D&=\cA_{\ell;1}^D\!\cup\!\cA_{\ell;2}^D,\\
\cA_{\ell}^H&=\big\{\vr\!\in\!\cA_{\ell}^{\pm}\!: \ov\vr\!=\!\vr\big\}, 
&\quad
\cA_{\ell}^E&=\big\{\vr\!\in\!\cA_{\ell}^{\pm}\!:\ov\vr\!=\!\vr_{\ell^{\pm}}^c\big\},
&\quad
\cA_{\ell}^{\R}&=
\cA_{\ell}^D\!\cup\!\cA_{\ell}^H\!\cup\!\cA_{\ell}^E,
\end{aligned}\\
\cA_{\ell}^{\star}=\big\{\vr\!\in\!\cA_{\ell}^{\R}\!:
\ov\vr^c_{\ell^{\pm}}\!\not\in\!\cA_{\ell;1}^D\big\},
\quad
\wt\cA_{\ell}^{\R}=\cA_{\ell}^{\R}\!\sqcup\!
\big\{\!\big[\ell^{\pm}\big]\!-\!\{i\}\!:i\!\in\!\big[\ell^{\pm}\big]\!\big\}.
\end{gather*}
If $\vr\!\in\!\cA_{\ell;1}^D$, then $\ov\vr^c_{\ell^{\pm}}\!\in\!\cA_{\ell;2}^D$.
If $\vr\!\in\!\cA_{\ell;2}^D$, then either $\ov\vr^c_{\ell^{\pm}}\!\in\!\cA_{\ell;1}^D$
or $\ov\vr\!\in\!\cA_{\ell;2}^D$, but not both.\\

\noindent
For $\vr\!\subset\![\ell^{\pm}]$, we denote by 
$$D_{\ell;\vr}\subset\R\ov\cM_{0,\ell}$$
the subspace parametrizing the $[\ell^{\pm}]$-marked curves~$\cC$ with a node $\nod_{\vr}(\cC)$ 
that separates~$\cC$ into two topological components, 
$\cC_{\vr}'$ and~$\cC_{\vr}''$, so that $\cC_{\vr}'$ (resp.~$\cC_{\vr}''$) carries
the marked points indexed by~$\vr$ (resp.~$[\ell^{\pm}]\!-\!\vr$).
This subspace is empty unless $\vr\!\in\!\cA_{\ell}^{\R}$.
It is an $H$-hypersurface (resp.~$E$-hypersurface, boundary divisor) if and only if 
$\vr\!\in\!\cA_{\ell}^H$ (resp.~$\vr\!\in\!\cA_{\ell}^E$, $\vr\!\in\!\cA_{\ell}^D$).
Every boundary hypersurface and divisor of~$\R\ov\cM_{0,\ell}$ equals~$D_{\ell;\vr}$
some $\vr\!\in\!\cA_{\ell}^{\R}$.
We note~that 
$$D_{\ell;\vr}=\begin{cases}D_{\ell;\ov\vr_{\ell^{\pm}}^c},
&\hbox{if}~\vr\!\in\!\cA_{\ell;1}^D;\\
D_{\ell;\ov\vr},&\hbox{if}~\vr,\ov\vr\!\in\!\cA_{\ell;2}^D.
\end{cases}$$
These are the only cases of 
distinct $\vr,\vr'\!\in\cA_{\ell}^{\R}$  with $D_{\ell;\vr}\!=\!D_{\ell;\vr'}$.\\

\noindent
We now define boundary hypersurfaces and divisors 
$$\wt{D}_{\vr}^{\bu}\subset \big\{\ff_{(\ell+1)^{\pm}}^{\R}\big\}^{-1}(D_{\ell;\vr})
\subset\R\ov\cM_{0,\ell+1} \qquad\hbox{with}\quad \bu=+,0,-.$$
These strata are illustrated in Figure~\ref{DT_fig}.
For $\vr\!\in\!\cA_{\ell}^E\!\cup\!\cA_{\ell;1}^D$, let
$$\wt{D}_{\vr}^+=D_{\ell+1;\vr\cup\{(\ell+1)^+\}}, \quad
\wt{D}_{\vr}^0=D_{\ell+1;\vr}, \quad
\wt{D}_{\vr}^-=D_{\ell+1;\vr\cup\{(\ell+1)^-\}}.$$
If $\vr\!\in\!\cA_{\ell}^H$ (resp.~$\vr\!\in\!\cA_{\ell;2}^D$), 
\hbox{$D_{\ell+1;\vr\cup\{(\ell+1)^-\}}\!=\!\eset$}
(resp.~$D_{\ell+1;\vr}\!=\!\eset$).
In these cases, we~set
$$\wt{D}_{\vr}^+=D_{\ell+1;\vr\cup\{(\ell+1)^+,(\ell+1)^-\}}\,,
\quad
\wt{D}_{\vr}^0,\wt{D}_{\vr}^-=
\begin{cases}
D_{\ell+1;\vr},&\hbox{if}~\vr\!\in\!\cA_{\ell}^H;\\
D_{\ell+1;\vr\cup\{(\ell+1)^-\}},&\hbox{if}~
\vr\!\in\!\cA_{\ell;2}^D.
\end{cases} $$
If $\vr\!\in\!\cA_{\ell;1}^D$, then
$$\wt{D}_{\vr}^0=\wt{D}_{\ov\vr_{\ell^{\pm}}^c}^+ \quad\hbox{and}\quad
\wt{D}_{\vr}^-=\wt{D}_{\ov\vr_{\ell^{\pm}}^c}^0,\wt{D}_{\ov\vr_{\ell^{\pm}}^c}^-\,.$$
If $\vr,\ov\vr\!\in\!\cA_{\ell;2}^D$, then $\wt{D}_{\vr}^+\!=\!\wt{D}_{\ov\vr}^+$.
These are the only cases of pairs~$(\vr,\bu)$ and~$(\vr',\circ)$ 
with distinct $\vr,\vr'\!\in\!\cA_{\ell}^{\R}$ such that 
$\wt{D}_{\vr}^{\bu}\!=\!\wt{D}_{\vr'}^{\circ}$.
For $i\!\in\![\ell^{\pm}]$, let 
\BE{wtDidfn_e}\wt{D}_{[\ell^{\pm}]-\{i\}}^+=\eset\subset\ov\cM_{0,\ell+1}\,,\quad
\wt{D}_{[\ell^{\pm}]-\{i\}}^0,\wt{D}_{[\ell^{\pm}]-\{i\}}^-=
D_{\ell+1;([\ell^{\pm}]-\{i\})\cup\{(\ell+1)^-\}}.\EE

\vspace{.15in}

\begin{figure}
\begin{pspicture}(-1,-9.2)(10,5.6)
\psset{unit=.4cm}
\rput(0,7.5){\framebox{$\vr\!\in\!\cA_{\ell}^H$}}
\pscircle[linewidth=.05](2.5,11){1.5}\pscircle[linewidth=.05](5.5,11){1.5}
\pscircle*(1.44,12.06){.15}\pscircle*(1.44,9.94){.15}
\rput(1.1,12.7){\sm{$1^+$}}\rput(1.1,9.3){\sm{$1^-$}}
\rput(1.6,11.1){\sm{$\vr$}}\rput(6.3,11.1){\sm{$\vr_{\ell^{\pm}}^c$}}
\psline[linewidth=.025]{<->}(7.5,9.5)(7.5,12.5)\rput(8,12){$\si$}
\rput(5.5,7.5){\sm{$D_{\ell;\vr}$}}
\pscircle[linewidth=.05](12.5,11){1.5}\pscircle[linewidth=.05](15.5,11){1.5}
\pscircle*(11.44,12.06){.15}\pscircle*(11.44,9.94){.15}
\pscircle*(12.5,12.5){.15}\pscircle*(12.5,9.5){.15}
\rput(11.1,12.7){\sm{$1^+$}}\rput(11.1,9.3){\sm{$1^-$}}
\rput(11.6,11.1){\sm{$\vr$}}\rput(16.3,11.1){\sm{$\vr_{\ell^{\pm}}^c$}}
\psline[linewidth=.025]{<->}(17.5,9.5)(17.5,12.5)\rput(18,11){$\si$}
\rput(13.5,13.2){\sm{$(\ell\!+\!1)^+$}}\rput(13.5,8.8){\sm{$(\ell\!+\!1)^-$}}
\rput(14.5,7.5){\sm{$\wt{D}_\vr^+\!\subset\!\R\ov\cM_{0,\ell+1}$}}
\pscircle[linewidth=.05](22.5,11){1.5}\pscircle[linewidth=.05](25.5,11){1.5}
\pscircle*(21.44,12.06){.15}\pscircle*(21.44,9.94){.15}
\pscircle*(25.5,12.5){.15}\pscircle*(25.5,9.5){.15}
\rput(21.1,12.7){\sm{$1^+$}}\rput(21.1,9.3){\sm{$1^-$}}
\rput(21.6,11.1){\sm{$\vr$}}\rput(26.3,11.1){\sm{$\vr_{\ell^{\pm}}^c$}}
\psline[linewidth=.025]{<->}(27.5,9.5)(27.5,12.5)\rput(28,11){$\si$}
\rput(25.5,13.2){\sm{$(\ell\!+\!1)^+$}}\rput(25.5,8.8){\sm{$(\ell\!+\!1)^-$}}
\rput(24.5,7.5){\sm{$\wt{D}_\vr^0,\wt{D}_\vr^-\!\subset\!\R\ov\cM_{0,\ell+1}$}}
\rput(0,-2.5){\framebox{$\vr\!\in\!\cA_{\ell}^E$}}
\psline[linewidth=.05](2,1)(6,5)\psline[linewidth=.05](2,3)(6,-1)
\pscircle*(3.5,2.5){.15}\pscircle*(5,4){.15}
\pscircle*(3.5,1.5){.15}\pscircle*(5,0){.15}
\rput(3.5,3.1){\sm{$1^{\pm}$}}\rput(3.2,.9){\sm{$1^{\mp}$}}
\rput(4.9,4.6){\sm{$2^+$}}\rput(4.7,-.7){\sm{$2^-$}}
\rput(6.5,5.1){\sm{$\vr$}}\rput(6.7,-1){\sm{$\vr_{\ell^{\pm}}^c$}}
\psline[linewidth=.025]{<->}(5.5,.5)(5.5,3.5)\rput(6,2){$\si$}
\rput(5.5,-2.5){\sm{$D_{\ell;\vr}$}}
\psline[linewidth=.05](12,1)(16,5)\psline[linewidth=.05](12,3)(16,-1)
\pscircle*(13.5,2.5){.15}\pscircle*(15,4){.15}\pscircle*(14.25,3.25){.15}
\pscircle*(13.5,1.5){.15}\pscircle*(15,0){.15}\pscircle*(14.25,.75){.15}
\rput(13.3,3){\sm{$1^{\pm}$}}\rput(13.2,.9){\sm{$1^{\mp}$}}
\rput(15,3){\sm{$2^+$}}\rput(15,1.2){\sm{$2^-$}}
\rput(13.5,-.5){\sm{$(\ell\!+\!1)^-$}}
\rput(14,4.6){\sm{$(\ell\!+\!1)^+$}}
\psline[linewidth=.025]{<->}(16.5,.5)(16.5,3.5)\rput(17,2){$\si$}
\rput(16.5,5.1){\sm{$\vr$}}\rput(16.7,-1){\sm{$\vr_{\ell^{\pm}}^c$}}
\rput(14.5,-2.5){\sm{$\wt{D}_\vr^+\!\subset\!\R\ov\cM_{0,\ell+1}$}}
\psline[linewidth=.05](23,4.5)(23,-.5)
\psline[linewidth=.02](22.5,4)(26,4)\psline[linewidth=.02](22.5,0)(26,0)
\pscircle*(24,4){.15}\pscircle*(25.5,4){.15}\pscircle*(25.5,0){.15}
\pscircle*(23,3){.15}\pscircle*(23,1){.15}\pscircle*(24,0){.15}
\rput(24,4.7){\sm{$1^{\pm}$}}\rput(24,-.7){\sm{$1^{\mp}$}}
\rput(25.5,4.7){\sm{$2^+$}}\rput(25.5,-.7){\sm{$2^-$}}
\psline[linewidth=.025]{<->}(25.5,.6)(25.5,3.4)\rput(26,2){$\si$}
\rput(26.5,4.1){\sm{$\vr$}}\rput(26.7,0){\sm{$\vr_{\ell^{\pm}}^c$}}
\rput(21.5,3.2){\sm{$(\ell\!+\!1)^+$}}\rput(21.5,1.2){\sm{$(\ell\!+\!1)^-$}}
\rput(24.5,-2.5){\sm{$\wt{D}_\vr^0\!\subset\!\R\ov\cM_{0,\ell+1}$}}
\psline[linewidth=.05](32,1)(36,5)\psline[linewidth=.05](32,3)(36,-1)
\pscircle*(33.5,2.5){.15}\pscircle*(35,4){.15}\pscircle*(34.25,3.25){.15}
\pscircle*(33.5,1.5){.15}\pscircle*(35,0){.15}\pscircle*(34.25,.75){.15}
\rput(33.3,3){\sm{$1^{\pm}$}}\rput(33.2,.9){\sm{$1^{\mp}$}}
\rput(35,3){\sm{$2^+$}}\rput(35,1.2){\sm{$2^-$}}
\rput(33.5,-.5){\sm{$(\ell\!+\!1)^+$}}
\rput(34,4.6){\sm{$(\ell\!+\!1)^-$}}
\psline[linewidth=.025]{<->}(36.5,.5)(36.5,3.5)\rput(37,2){$\si$}
\rput(36.5,5.1){\sm{$\vr$}}\rput(36.7,-1){\sm{$\vr_{\ell^{\pm}}^c$}}
\rput(34.5,-2.5){\sm{$\wt{D}_\vr^-\!\subset\!\R\ov\cM_{0,\ell+1}$}}
\rput(0,-12.5){\framebox{$\vr\!\in\!\cA_{\ell;1}^D$}}
\psline[linewidth=.05](3,-5.5)(3,-10.5)
\psline[linewidth=.02](2.5,-6)(6,-6)\psline[linewidth=.02](2.5,-10)(6,-10)
\pscircle*(3.5,-6){.15}\pscircle*(4.5,-6){.15}
\pscircle*(3.5,-10){.15}\pscircle*(4.5,-10){.15}
\rput(3.7,-5.3){\sm{$1^{\pm}$}}\rput(3.7,-10.7){\sm{$1^{\mp}$}}
\rput(4.7,-6.6){\sm{$2^+$}}\rput(4.7,-9.3){\sm{$2^-$}}
\psline[linewidth=.025]{<->}(5.8,-9.4)(5.8,-6.6)\rput(6.3,-8){$\si$}
\rput(6.5,-5.9){\sm{$\vr$}}\rput(6.7,-10){\sm{$\ov\vr$}}
\rput(1.6,-8){\sm{$\vr_{\ell^{\pm}}^c\!-\!\ov\vr$}}
\rput(5.5,-12.5){\sm{$D_{\ell;\vr}$}}
\psline[linewidth=.05](13,-5.5)(13,-10.5)
\psline[linewidth=.02](12.5,-6)(16,-6)\psline[linewidth=.02](12.5,-10)(16,-10)
\pscircle*(13.5,-6){.15}\pscircle*(14.5,-6){.15}\pscircle*(15.5,-6){.15}
\pscircle*(13.5,-10){.15}\pscircle*(14.5,-10){.15}\pscircle*(15.5,-10){.15}
\rput(13.7,-5.3){\sm{$1^{\pm}$}}\rput(13.7,-10.7){\sm{$1^{\mp}$}}
\rput(14.7,-6.6){\sm{$2^+$}}\rput(14.7,-9.3){\sm{$2^-$}}
\rput(16.5,-5.3){\sm{$(\ell\!+\!1)^+$}}\rput(16.5,-10.7){\sm{$(\ell\!+\!1)^-$}}
\psline[linewidth=.025]{<->}(15.8,-9.4)(15.8,-6.6)\rput(16.3,-8){$\si$}
\rput(16.5,-6.2){\sm{$\vr$}}\rput(16.7,-9.7){\sm{$\ov{\vr}$}}
\rput(11.6,-8){\sm{$\vr_{\ell^{\pm}}^c\!-\!\ov{\vr}$}}
\rput(14.5,-12.5){\sm{$\wt{D}_\vr^+\!\subset\!\R\ov\cM_{0,\ell+1}$}}
\psline[linewidth=.05](23,-5.5)(23,-10.5)
\psline[linewidth=.02](22.5,-6)(26,-6)\psline[linewidth=.02](22.5,-10)(26,-10)
\pscircle*(23.5,-6){.15}\pscircle*(24.5,-6){.15}
\pscircle*(23.5,-10){.15}\pscircle*(24.5,-10){.15}
\pscircle*(23,-6.5){.15}\pscircle*(23,-9.5){.15}
\rput(23.7,-5.3){\sm{$1^{\pm}$}}\rput(23.7,-10.7){\sm{$1^{\mp}$}}
\rput(24.7,-6.6){\sm{$2^+$}}\rput(24.7,-9.3){\sm{$2^-$}}
\psline[linewidth=.025]{<->}(25.8,-9.4)(25.8,-6.6)\rput(26.3,-8){$\si$}
\rput(26.5,-5.9){\sm{$\vr$}}\rput(26.7,-10){\sm{$\ov{\vr}$}}
\rput(21.6,-8){\sm{$\vr_{\ell^{\pm}}^c\!-\!\ov{\vr}$}}
\rput(21.3,-6.8){\sm{$(\ell\!+\!1)^+$}}\rput(21.4,-9.2){\sm{$(\ell\!+\!1)^-$}}
\rput(24.5,-12.5){\sm{$\wt{D}_\vr^0\!\subset\!\R\ov\cM_{0,\ell+1}$}}
\psline[linewidth=.05](33,-5.5)(33,-10.5)
\psline[linewidth=.02](32.5,-6)(36,-6)\psline[linewidth=.02](32.5,-10)(36,-10)
\pscircle*(33.5,-6){.15}\pscircle*(34.5,-6){.15}\pscircle*(35.5,-6){.15}
\pscircle*(33.5,-10){.15}\pscircle*(34.5,-10){.15}\pscircle*(35.5,-10){.15}
\rput(33.7,-5.3){\sm{$1^{\pm}$}}\rput(33.7,-10.7){\sm{$1^{\mp}$}}
\rput(34.7,-6.6){\sm{$2^+$}}\rput(34.7,-9.3){\sm{$2^-$}}
\rput(36.5,-5.3){\sm{$(\ell\!+\!1)^-$}}\rput(36.5,-10.7){\sm{$(\ell\!+\!1)^+$}}
\psline[linewidth=.025]{<->}(35.8,-9.4)(35.8,-6.6)\rput(36.3,-8){$\si$}
\rput(36.5,-6.2){\sm{$\vr$}}\rput(36.7,-9.7){\sm{$\ov{\vr}$}}
\rput(31.6,-8){\sm{$\vr_{\ell^{\pm}}^c\!-\!\ov{\vr}$}}
\rput(34.5,-12.5){\sm{$\wt{D}_\vr^-\!\subset\!\R\ov\cM_{0,\ell+1}$}}
\rput(0,-22.5){\framebox{$\vr\!\in\!\cA_{\ell;2}^D$}}
\psline[linewidth=.05](3,-15.5)(3,-20.5)
\psline[linewidth=.02](2.5,-16)(6,-16)\psline[linewidth=.02](2.5,-20)(6,-20)
\psline[linewidth=.025]{<->}(5.8,-19.4)(5.8,-16.6)\rput(6.3,-18){$\si$}
\rput(6.7,-15.9){\sm{$\ov{\vr}_{\ell^{\pm}}^c$}}\rput(6.7,-20){\sm{$\vr_{\ell^{\pm}}^c$}}
\rput(1.6,-18){\sm{$\ov{\vr}\!-\!\vr_{\ell^{\pm}}^c$}}
\rput(5.8,-22.5){\sm{$D_{\ell;\vr}$}}
\psline[linewidth=.05](13,-15.5)(13,-20.5)
\psline[linewidth=.02](12.5,-16)(16,-16)\psline[linewidth=.02](12.5,-20)(16,-20)
\pscircle*(13,-16.5){.15}\pscircle*(13,-19.5){.15}
\psline[linewidth=.025]{<->}(15.8,-19.4)(15.8,-16.6)\rput(16.3,-18){$\si$}
\rput(16.7,-15.9){\sm{$\ov{\vr}_{\ell^{\pm}}^c$}}\rput(16.7,-20){\sm{$\vr_{\ell^{\pm}}^c$}}
\rput(11.6,-18){\sm{$\ov{\vr}\!-\!\vr_{\ell^{\pm}}^c$}}
\rput(11.3,-16.8){\sm{$(\ell\!+\!1)^+$}}\rput(11.4,-19.2){\sm{$(\ell\!+\!1)^-$}}
\rput(14.5,-22.5){\sm{$\wt{D}_\vr^+\!\subset\!\R\ov\cM_{0,\ell+1}$}}
\psline[linewidth=.05](33,-15.5)(33,-20.5)
\psline[linewidth=.02](32.5,-16)(36,-16)\psline[linewidth=.02](32.5,-20)(36,-20)
\pscircle*(34.5,-16){.15}\pscircle*(34.5,-20){.15}
\rput(35.5,-15.3){\sm{$(\ell\!+\!1)^-$}}\rput(35.5,-20.7){\sm{$(\ell\!+\!1)^+$}}
\psline[linewidth=.025]{<->}(35.8,-19.4)(35.8,-16.6)\rput(36.3,-18){$\si$}
\rput(36.7,-16.5){\sm{$\ov\vr_{\ell^{\pm}}^c$}}\rput(36.7,-19.5){\sm{$\vr_{\ell^{\pm}}^c$}}
\rput(31.6,-18){\sm{$\ov\vr\!-\!\vr_{\ell^{\pm}}^c$}}
\rput(34.5,-22.5){\sm{$\wt{D}_\vr^0,\wt{D}_\vr^-\!\subset\!\R\ov\cM_{0,\ell+1}$}}
\end{pspicture}
\caption{Generic representatives of the top boundary strata $D_{\vr}$ in $\R\ov\cM_{0,\ell}$ 
and $\wt{D}_{\vr}^{\bu}$ in $\R\ov\cM_{0,\ell+1}$.
Each line or circle represents~$\C\P^1$.
Each double-headed arrow labeled~$\si$ indicates
the involution on the corresponding real curve.}
\label{DT_fig}
\end{figure}

\noindent
The set~$\cA_{\ell}^{\R}$ is partially ordered by the inclusion~$\subsetneq$ of 
subsets of~$[\ell^{\pm}]$.
We extend this partial order to a strict order~$<$ on $\{0\}\!\sqcup\!\cA_{\ell}^{\R}$
so that~0 is the smallest element and define \hbox{$\vr\!-\!1\!\in\!\{0\}\!\cup\!\cA_{\ell}^{\R}$}
to be the predecessor of $\vr\!\in\!\cA_{\ell}^{\R}$.
For \hbox{$\vr^*\!\in\!\{0\}\!\sqcup\!\cA_{\ell}^{\R}$}, let 
$$\cA_{\ell}^{\R}(\vr^*)=\big\{\vr\!\in\!\cA_{\ell}^{\R}\!:\vr\!>\!\vr^*\big\}.$$

\vspace{.15in}

\noindent 
For $\vr\!\in\!\cA_{\ell}^{\R}$ and $\wt\cC_1,\wt\cC_2\!\in\!\R\ov\cM_{0,\ell+1}$, 
we define $\wt\cC_1\!\sim_{\vr}\!\wt\cC_2$ if either $\wt\cC_1\!=\!\wt\cC_2$ or 
$$\ff_{\ell+1}^{\R}(\wt\cC_1)=\ff_{\ell+1}^{\R}(\wt\cC_2)\in\R\ov\cM_{0,\ell}
\quad\hbox{and}\quad \wt\cC_1,\wt\cC_2\in \wt{D}_{\vr}^0\!\cup\!\wt{D}_{\vr}^-.$$
For every $\vr^*\!\in\!\{0\}\!\sqcup\!\cA_{\ell}^{\R}$
the union of all equivalence relations~$\sim_{\vr}$ on~$\R\ov\cM_{0,\ell+1}$ with 
\hbox{$\vr\!\in\!\cA_{\ell}^{\R}(\vr^*)$} is again an equivalence relation;
see Section~1.2 in~\cite{RDMbl}.
The quotient~$X_{\vr^*}$ of~$\R\ov\cM_{0,\ell+1}$ by the last equivalence relation
is a compact Hausdorff space for every \hbox{$\vr^*\!\in\!\{0\}\!\sqcup\!\cA_{\ell}^{\R}$}.
Let
\BE{Rpivrdfn_e}
q_{\vr^*}\!:\R\ov\cM_{0,\ell+1}\lra X_{\vr^*} \qquad\hbox{and}\qquad
\pi_{\vr^*}\!:X_{\vr^*}\lra X_{\vr^*-1}\EE
be the quotient projection and the map induced by~$q_{\vr^*-1}$
whenever $\vr^*\!\in\!\cA_{\ell}^{\R}$, respectively.
For $\vr\!\in\!\wt\cA_{\ell}^{\R}$ and $\bu\!=\!+,0,-$, define
$$Y_{\vr^*;\vr}^{\bu}=q_{\vr^*}\big(\wt{D}_{\vr}^{\bu}\big).$$
This subspace of $X_{\vr^*}$ is compact.
If $\vr\!\in\!\cA_{\ell}^{\R}(\vr^*)$, $Y_{\vr^*;\vr}^0\!=\!Y_{\vr^*;\vr}^-$.\\

\noindent
By Theorem~2.2($\R1$) in~\cite{RDMbl},
$X_{\vr^*}$ is a smooth manifold,
$Y_{\vr^*;\vr}^{\bu}\!\subset\!X_{\vr^*}$ is a smooth submanifold, and $q_{\vr^*}$ is a smooth map.
By the construction of the smooth structure on~$X_{\vr^*}$ at the beginning
of Section~5.4 in~\cite{RDMbl}, the restriction 
\BE{RblowupPf_e0}
q_{\vr^*}\!:\R\ov\cM_{0,\ell+1}-\bigcup_{\vr\in\cA_{\ell}^{\R}(\vr^*)}\!\!\!\!\!\!\!\!
\big(\wt{D}_{\vr}^0\!\cup\!\wt{D}_{\vr}^-\big)
\lra X_{\vr^*}-\bigcup_{\vr\in\cA_{\ell}^{\R}(\vr^*)}\!\!\!\!\!\!Y_{\vr^*;\vr}^0\EE
is a diffeomorphism.
Since the fibers of~$q_{\vr^*}$ over a submanifold~$Y_{\vr^*;\vr}^0\!\subset\!X_{\vr^*}$
with $\vr\!\in\!\cA_{\ell}^{\R}(\vr^*)$ are positive-dimensional,
the codimension of~$Y_{\vr^*;\vr}^0$ is at least~2 (in fact, it is~3 or~4).
Since the left-hand side in~\eref{RblowupPf_e0} is orientable,
it follows that the manifold~$X_{\vr^*}$ is also orientable.\\

\noindent
For each $\wt\cC\!\in\!\R\ov\cM_{0,\ell+1}$, the set
$$\cA_{\wt\cC}^{\R}(\vr^*)\equiv\big\{\vr'\!\in\!\cA_{\ell}^{\R}(\vr^*)\!:
\wt\cC\!\in\!\wt{D}_{\vr'}^0\!\cup\!\wt{D}_{\vr'}^-\big\}\subset
\cA^{\R}_{\ell}(\vr^*)$$
is ordered by the inclusion of subsets of~$[\ell^{\pm}]$.
If $\cA_{\wt\cC}^{\R}(\vr^*)\!\neq\!\eset$, 
let $\vr_{\min}(\wt\cC)\!\in\!\cA_{\wt\cC}^{\R}(\vr^*)$ be the minimal element.
By the definition of~$q_{\vr^*}$,
\BE{qvrpreim_e}q_{\vr^*}^{\,-1}\!\big(q_{\vr^*}(\wt\cC)\!\big)=
\begin{cases}\{\wt\cC\},&\hbox{if}~\cA_{\wt\cC}^{\R}(\vr^*)\!=\!\eset;\\
\big\{\ff_{(\ell+1)^{\pm}}^{\R}\big\}^{\!-1}\!
\big(\ff_{(\ell+1)^{\pm}}^{\R}(\wt\cC)\!\big)
\!\cap\!\big(\wt{D}_{\!\vr_{\min}(\wt\cC)}^0\!\cup\!\wt{D}_{\!\vr_{\min}(\wt\cC)}^-\big),
&\hbox{if}~\cA_{\wt\cC}^{\R}(\vr^*)\!\neq\!\eset.
\end{cases}\EE
For $\vr\!\in\!\wt\cA_{\ell}^{\R}$ and $\bu\!=\!+,0,-$, 
let $\cA_{\vr}^{\bu}(\vr^*)\!\subset\!\cA_{\ell}^{\R}(\vr^*)$ be 
the subset of elements~$\vr'$ such~that the restriction
$$\ff_{(\ell+1)^{\pm}}^{\R}\!:
\wt{D}_{\vr}^{\bu}\!\cap\!\big(\wt{D}_{\vr'}^0\!\cup\!\wt{D}_{\vr'}^-\big)
\lra \R\ov\cM_{0,\ell}$$
is not injective.
The quotient projection~$q_{\vr^*}$ restricts to a homeomorphism
$$q_{\vr^*}\!:\wt{D}_{\vr}^{\bu}-
\bigcup_{\vr'\in\cA_{\vr}^{\bu}(\vr^*)}\hspace{-.22in}
\big(\wt{D}_{\vr'}^0\!\cup\!\wt{D}_{\vr'}^-\big)
\lra Y_{\vr^*;\vr}^{\bu}\!-\!\bigcup_{\vr'\in\cA_{\vr}^{\bu}(\vr^*)}\hspace{-.2in}
Y_{\vr^*;\vr'}^0\,.$$
By~\eref{qvrpreim_e} and the paragraph containing~\eref{REinter_e}-\eref{RDinter_e2}
with~$\ell$ replaced by~$\ell\!+\!1$,
the fibers of the restriction
$$q_{\vr^*}\!:\wt{D}_{\vr}^{\bu}\cap
\bigcup_{\vr'\in\cA_{\vr}^{\bu}(\vr^*)}\hspace{-.22in}
\big(\wt{D}_{\vr'}^0\!\cup\!\wt{D}_{\vr'}^-\big)
\lra Y_{\vr^*;\vr}^{\bu}\!\cap\!\bigcup_{\vr'\in\cA_{\vr}^{\bu}(\vr^*)}\hspace{-.2in}
Y_{\vr^*;\vr'}^0$$
are connected unions of submanifolds of~$\wt{D}_{\vr}^+$ for every $\vr'\!\in\!\cA_{\ell}^{\R}$.\\

\noindent
If $\vr\!\in\!\cA_{\ell}^{\star}$, then 
$\wt{D}_{\vr}^+\!\cap\!\wt{D}_{\vr'}^0$ and $\wt{D}_{\vr}^+\!\cap\!\wt{D}_{\vr'}^-$ 
are submanifolds of~$\wt{D}_{\vr}^+$
of positive codimension for every \hbox{$\vr'\!\in\!\cA_{\ell}^{\R}$}.
Thus, 
$$\dim\wt{D}_{\vr}^+\!\cap\!\wt{D}_{\vr'}^0,\dim\wt{D}_{\vr}^+\!\cap\!\wt{D}_{\vr'}^-
<\dim\wt{D}_{\vr}^+=\dim Y_{\vr^*;\vr}^+
\quad\forall\,\vr'\!\in\!\cA_{\ell}^{\R}\,.$$
Along with the previous paragraph, this implies that an orientation on~$\wt{D}_{\vr}^+$,
if any, induces an orientation on~$Y_{\vr^*;\vr}^+$.
By~\eref{Rboundident_e} with~$\ell$ replaced by~$\ell\!+\!1$,
the submanifold $Y_{\vr^*;\vr}^+\!\subset\!X_{\vr^*}$ is thus orientable for every
\hbox{$\vr\!\in\!\cA_{\ell}^{\star}\!-\!\cA_{\ell}^H$}.
By the same reasoning, the submanifolds
$$Y_{\vr^*;\vr}^0 ~~\hbox{with}~~
\vr\!\in\!\wt\cA_{\ell}^{\R}\!-\!\cA_{\ell}^H\!-\!\cA_{\ell}^{\R}(\vr^*)
\quad\hbox{and}\quad
Y_{\vr^*;\vr}^- ~~\hbox{with}~~
\vr\!\in\!\cA_{\ell}^E\!-\!\cA_{\ell}^{\R}(\vr^*)$$ 
are orientable as~well.\\

\noindent
We assume that property~\ref{HCgen_it} on page~\pageref{HCgen_it} holds with~$\ell$ replaced 
by each $\ell'\!\le\!\ell\!+\!1$
and that property~\ref{HRgen_it} on page~\pageref{HRgen_it}
holds with~$\ell$ replaced by each $\ell'\!\le\!\ell$.
We will establish the following statements by induction 
on \hbox{$\vr^*\!\in\!\{0\}\!\cup\!\cA_{\ell}^{\R}$}. 
\begin{enumerate}[label=($\R1\alph*$),leftmargin=*]

\item\label{RYgen_it}  
The submanifolds $Y_{\vr^*;[\ell^{\pm}]-\{1^+\}}^0\!\subset\!X_{\vr^*}$,
$Y_{\vr^*;\vr}^+$  with $\vr\!\in\!\cA_{\ell}^{\star}\!-\!\cA_{\ell}^H$,
$Y_{\vr^*;\vr}^0$ with
\hbox{$\vr\!\in\!\cA_{\ell}^{\R}\!-\!\cA_{\ell}^H\!-\!\cA_{\ell}^{\R}(\vr^*)$},
and~$Y_{\vr^*;\vr}^-$ with
$\vr\!\in\!\cA_{\ell}^E\!-\!\cA_{\ell}^{\R}(\vr^*)$
generate $H_*(X_{\vr^*};\Q)$ as a~ring. 

\item\label{RYTincomp_it}
If $\vr\!\in\!\cA_{\ell}^{\R}(\vr^*)\!-\!\cA_{\ell}^H$, 
the submanifold $Y_{\vr^*;\vr}^0\!\subset\!X_{\vr^*}$ is 
orientable and $\Q$-incompressible.

\end{enumerate}
If $\vr^*\!\in\!\cA_{\ell}^{\R}$ is the maximal element, 
\ref{RYgen_it} is a stronger version of
property~\ref{HRgen_it} with~$\ell$ replaced by~$\ell\!+\!1$.\\

\noindent
By Theorem~2.2($\R1$) in~\cite{RDMbl} and the construction of the smooth structure on~$X_{\vr^*}$ 
at the beginning of Section~5.4 in~\cite{RDMbl}, the map~\eref{RMlind_e} descends to a diffeomorphism 
$$\Psi_0\!:X_0\lra \R\ov\cM_{0,\ell}\!\times\!\ov\cM_4\,.$$
We identify~$X_0$ with \hbox{$\R\ov\cM_{0,\ell}\!\times\!\ov\cM_4$} via~$\Psi_0$.
Under this identification,
\begin{gather*}
Y_{0;\vr}^+\equiv\big\{\!(\ff_{\ell+1}^{\R},\ff_{1^{\pm}2^+(\ell+1)^+}^{\R})\!\big\}
(\wt{D}_{\vr}^+)=D_{\vr}\!\times\!\ov\cM_4 \quad\forall~
\vr\!\in\!\cA_{\ell}^{\star},\\
Y_{0;[\ell^{\pm}]-\{1^+\}}^0=\big\{\!(\ff_{\ell+1}^{\R},\ff_{1^{\pm}2^+(\ell+1)^+}^{\R})\!\big\}
\big(\wt{D}_{[\ell^{\pm}]-\{1^+\}}^0\big)
=\R\ov\cM_{0,\ell}\!\times\!D_{14,23}.
\end{gather*}
Since the submanifolds~$D_{\vr}$ with $\vr\!\in\!\cA_{\ell}^{\star}\!-\!\cA_{\ell}^H$
generate $H_*(\R\ov\cM_{0,\ell};\Q)$ as a ring,
the submanifolds
$Y_{0;[\ell^{\pm}]-\{1^+\}}^0,Y_{0;\vr}^+\!\subset\!X_0$ generate $H_*(X_0;\Q)$ as a ring.
This confirms the $\vr^*\!=\!0$ case of~\ref{RYgen_it}.\\

\noindent
Let $\vr\!\in\!\cA_{\ell}^{\R}$.
As in the complex case of Section~\ref{GenC_subs}, the forgetful morphism
$$\ff_{\vr}\!\equiv\!\ff_{1^{\pm}2^+j}^{\R}\!:D_{\vr}\lra \ov\cM_4$$
is independent of the choice of $j\!\in\!\vr_{\ell^{\pm}}^c$ distinct from $1^{\pm},2^+$.
The~map
$$s_{\vr}\!=\!\big(\io_{\vr},\ff_{\vr}\big)\!:
D_{\vr}\lra \R\ov\cM_{0,\ell}\!\times\!\ov\cM_4\!=\!X_0,$$
where $\io_{\vr}\!:D_{\vr}\!\lra\!\R\ov\cM_{0,\ell}$ is the inclusion,
is again a smooth embedding and 
$$Y_{0;\vr}^0\equiv\big\{\!
(\ff_{\ell+1}^{\R},\ff_{1^{\pm}2^{\pm}(\ell+1)^+}^{\R})\!\big\}(\wt{D}_{\vr}^0)
=s_{\vr}(D_{\vr}).$$
Thus, the submanifolds $Y_{0;\vr}^0\!\subset\!X_0$ with 
$\vr\!\in\!\cA_{\ell}^{\R}\!-\!\cA_{\ell}^H$ are  orientable.\\

\noindent
If in addition $\vr'\!\in\!\cA_{\ell}^{\star}$, the restriction
$$s_{\vr}\!:D_{\vr}\!\cap\!D_{\vr'}\lra Y_{0;\vr}^0\!\cap\!Y_{0;\vr'}^+ $$
is a diffeomorphism.
Since $Y_{0;\vr}^0$ and $Y_{0;\vr'}^+$ intersect transversely in~$X_0$ for $\vr\!\neq\!\vr'$,
the first identification in~\eref{REinter_e} and property~\ref{HCgen_it} imply that
the restriction homomorphism
$$Y_{0;\vr}^0\cap\!: H_*(X_0;\Q)\lra H_{*-3}(Y_{0;\vr}^0;\Q)$$
is surjective if $\vr\!\in\!\cA_{\ell}^E$, 
i.e.~$Y_{0;\vr}^0\!\subset\!X_0$ is $\Q$-incompressible in this case.
If $\vr\!\in\!\cA_{\ell}^D$, the same conclusion follows from
the second identification in~\eref{REinter_e}, \eref{RDinter_e},
the K\"unneth formula,  property~\ref{HCgen_it},
and the inductive assumption for property~\ref{HRgen_it}.
This confirms the $\vr^*\!=\!0$ case of~\ref{RYTincomp_it}.\\

\noindent
Suppose $\vr^*\!\in\!\cA_{\ell}^{\R}$ is such that~\ref{RYgen_it} and~\ref{RYTincomp_it} 
with~$\vr^*$ replaced by~$\vr^*\!-\!1$ hold.
By Theorem~2.2($\R2$) in~\cite{RDMbl}, 
the second map in~\eref{Rpivrdfn_e} is a real (resp.~complex, 1-augmented) blowup
 of~$X_{\vr^*-1}$ along~$Y_{\vr^*-1;\vr^*}^0$
with the exceptional locus~$Y_{\vr^*;\vr^*}^0\!\cup\!Y_{\vr^*;\vr^*}^-$
if $\vr^*\!\in\!\cA_{\ell}^H$ (resp.~$\vr^*\!\in\!\cA_{\ell}^D$, $\vr^*\!\in\!\cA_{\ell}^E$).
If \hbox{$\vr^*\!\in\!\cA_{\ell}^H\cup\!\cA_{\ell}^D$}, 
\hbox{$Y_{\vr^*;\vr^*}^-\!\subset\!Y_{\vr^*;\vr^*}^0$}.
By the first statement of Theorem~2.2($\R3$) in~\cite{RDMbl}, 
each submanifold $Y_{\vr^*;\vr}^{\bu}\!\subset\!X_{\vr^*}$ 
with \hbox{$\vr\!\in\!\wt\cA_{\ell}^{\R}\!-\!\cA_{\ell}^H$}, 
$$(\vr,\bu)\neq(\vr^*,0),(\vr^*,-), ~~
(\vr,\bu)\neq(\ov{\vr^*}_{\ell^{\pm}}^c,+)~\hbox{if}~\vr^*\!\in\!\cA_{\ell;1}^D, 
~~\hbox{and}~~
(\vr,\bu)\neq(\ov{\vr^*}_{\ell^{\pm}}^c,-)~\hbox{if}~\vr^*\!\in\!\cA_{\ell;2}^D$$
is $(\pi_{\vr^*},\Q)$-related to~$Y_{\vr^*-1;\vr}^{\bu}$.
Along with~\ref{RYgen_it} with~$\vr^*$ replaced by~$\vr^*\!-\!1$,
\ref{RYTincomp_it} with~$(\vr^*,\vr)$ replaced by~$(\vr^*\!-\!1,\vr^*)$
if $\vr^*\!\in\!\cA_{\ell}^D\!\cup\!\cA_{\ell}^E$,
and the first claim of Proposition~\ref{BlHom_prp}\ref{BlHomR_it}
(resp.~\ref{BlHom_prp}\ref{BlHomC_it}, \ref{BlHom_prp}\ref{BlHomRR_it}),
this implies that \ref{RYgen_it} as stated holds as~well.
By the second statement of Theorem~2.2($\R3$) in~\cite{RDMbl}, 
each submanifold \hbox{$Y_{\vr^*;\vr}^0\!\subset\!X_{\vr^*}$} with $\vr\!\in\!\cA_{\ell}^{\R}(\vr^*)$
is $\pi_{\vr^*}$-equivalent to the submanifold 
\hbox{$Y_{\vr^*-1;\vr}^0\!\subset\!X_{\vr^*-1}$}.
Along with~\ref{RYTincomp_it} with~$(\vr^*,\vr)$ replaced by~$(\vr^*\!-\!1,\vr)$
and the second claim of Proposition~\ref{BlHom_prp}\ref{BlHomR_it}
(resp.~\ref{BlHom_prp}\ref{BlHomC_it}, \ref{BlHom_prp}\ref{BlHomRR_it}),
this implies that \ref{RYTincomp_it} as stated holds as~well.

\section{Recursive relations}
\label{CohRelat_sec}

\noindent
It remains to show that the homomorphisms~\eref{HCthm_e} and~\eref{HRthm_e} are injective.
We do so indirectly by establishing recursive relations for the domains and 
targets of these maps separately.
In the proofs of Propositions~\ref{Chomol_prp} and~\ref{Rhomol_prp}
in Sections~\ref{Chomol_subs} and~\ref{Rhomol_subs2}, respectively, 
we show that certain homomorphisms between the homology groups of 
the moduli spaces of complex and real curves are injective.
In the proofs of Propositions~\ref{Calg_prp} and~\ref{Ralg_prp}
in Sections~\ref{Calg_subs} and~\ref{Ralg_subs}, respectively, 
we show that the ``corresponding" homomorphisms between the quotient algebras
in the domains of~\eref{HCthm_e} and~\eref{HRthm_e} are surjective.
The homomorphisms between the homology groups and the quotient algebras
``correspond" in the sense that the diagrams~\eref{HCdiag_e} and~\eref{HRdiag_e} commute. 
Along with the inductive assumptions on the homomorphisms~\eref{HCthm_e} and~\eref{HRthm_e}
and the surjectivity of these homomorphisms provided by Sections~\ref{GenC_subs} and~\ref{GenR_subs},
the four propositions and the commutativity of the diagrams~\eref{HCdiag_e} and~\eref{HRdiag_e}
imply that all homomorphisms in these diagrams are in fact isomorphisms;
this concludes the proofs of Theorems~\ref{HC_thm} and~\ref{HR_thm}.
For the convenience of use, such as in~\eref{ChomolRec_e} and~\eref{RhomolRec_e},
we state Propositions~\ref{Chomol_prp} and~\ref{Rhomol_prp} as providing
isomorphisms between the relevant homology groups,
even though the surjectivity of these maps is provided by the later 
Propositions~\ref{Calg_prp} and~\ref{Ralg_prp}, respectively.
In Section~\ref{Rhomol_subs1}, we orient the real moduli spaces $\R\ov\cM_{0,\ell}$,
the boundary hypersurfaces~$\R E_{J,K}$, and the boundary divisors~$\R D_{I;J,K}$  
defined in Section~\ref{Rthm_subs} and confirm that the general cases of~\eref{RcM04rel_e2a}
and~\eref{RcM04rel_e2b} indeed follow from the basic cases of these homology relations, 
as claimed in Section~\ref{Rthm_subs}.

\subsection{The complex case: topology}
\label{Chomol_subs}

\noindent 
Let $\ell\!\in\!\Z^+$ with $\ell\!\ge\!3$.
For $\ka\!\in\!H_p(\ov\cM_{\ell};\Z)$, let
\begin{equation*}\begin{split}
\ov\cM_{\ell+1}\big|_{\ka}
&=\PD_{\ov\cM_{\ell+1}}\!\big(\ff_{\ell+1}^{\,*}\big(\PD_{\ov\cM_{\ell}}\!(\ka)\!\big)\!\big)
\in H_{p+2}\big(\ov\cM_{\ell+1};\Z\big),\\
D_{\{1,\ell+1\},[\ell]-\{1\}}\!\big|_{\ka}
&=D_{\{1,\ell+1\},[\ell]-\{1\}}\!\cap\!\big(\ov\cM_{\ell+1}\big|_{\ka}\big)
\in H_p\big(\ov\cM_{\ell+1};\Z\big).
\end{split}\end{equation*}
Suppose $\{J,K\}\!\in\!\cP_{\bu}(\ell)$. 
If $J\!\ni\!1$, define
\BE{wtDdfn_e}\wt{D}_{J,K}^+=D_{J\cup\{\ell+1\},K}, \quad
\wt{D}_{J,K}^0=D_{J,K\cup\{\ell+1\}}\,.\EE
If $\ka_{J,K}\!\in\!H_{p-2}(D_{J,K};\Z)$ and $\bu\!=\!+,0$, let
$$\wt\ka_{J,K}^{\bu}=\big\{\wt\io_{J,K}^{\bu}\big\}_{\!*}
\Big(\PD_{\wt{D}_{J,K}^{\bu}}\!\big(
\big\{\ff_{\ell+1}|_{\wt{D}_{J,K}^{\bu}}\big\}^{\!*}(\PD_{D_{J,K}}\!(\ka_{J,K})\!\big)\!\big)
\!\Big) \in H_p\big(\ov\cM_{\ell+1};\Z\big),$$
where $\wt\io_{J,K}^{\bu}\!:\wt{D}_{J,K}^{\bu}\!\lra\!\ov\cM_{\ell+1}$ is the inclusion map.

\begin{prp}\label{Chomol_prp}
For every $\ell\!\in\!\Z^+$ with $\ell\!\ge\!3$ and $p\!\in\!\Z$, the homomorphism
\BE{Chomol_e}\begin{split}
&\Psi\!:H_p\big(\ov\cM_{\ell};\Z\big)\!\oplus\!H_{p-2}\big(\ov\cM_{\ell};\Z\big)
\!\oplus\!\bigoplus_{\{J,K\}\in\cP_{\bu}(\ell)}\hspace{-.24in}
H_{p-2}\big(D_{J,K};\Z\big)\lra H_p\big(\ov\cM_{\ell+1};\Z\big),\\
&\hspace{.1in}
\Psi\big(\ka_0,\ka,(\ka_{J,K})_{\{J,K\}\in\cP_{\bu}(\ell)}\big)
=D_{\{1,\ell+1\},[\ell]-\{1\}}\big|_{\ka_0}\!+\!\ov\cM_{\ell+1}\big|_{\ka}
\!+\!\sum_{\{J,K\}\in\cP_{\bu}(\ell)}\hspace{-.24in}\wt\ka_{J,K}^0\,,
\end{split}\EE
of $\Z$-modules is an isomorphism.
\end{prp}

\begin{proof} By Proposition~\ref{Calg_prp}, 
the top homomorphism in the diagram~\eref{HCdiag_e} is surjective.
By the generation property~\ref{HCgen_it} on page~\pageref{HCgen_it} with~$\ell$ replaced 
by~$\ell\!+\!1$, so is the right homomorphism.
By Theorem~\ref{HC_thm} with~$\ell$ replaced by each $\ell'\!\le\!\ell$,
the left homomorphism is an isomorphism.
Along with the commutativity of the diagram~\eref{HCdiag_e},
this implies the surjectivity of the bottom homomorphism~\eref{Chomol_e}. 
Below we establish the injectivity of~\eref{Chomol_e}.\\

\noindent
The restriction of~$\Psi$ to~$H_p(\ov\cM_{\ell};\Z)$ is injective because
$$\big\{\ff_{\ell+1}\big\}_{\!*}\!\circ\!\Psi\!=\!\id\!:
H_p\big(\ov\cM_{\ell};\Z\big)\lra H_p\big(\ov\cM_{\ell};\Z\big)\,.$$
Since the homomorphism $\{\ff_{\ell+1}\}_*$ vanishes on the images of~$H_{p-2}(\ov\cM_{\ell};\Z)$
and $H_{p-2}(D_{J,K};\Z)$ under~$\Psi$, it also follows that
$$\Psi\big(H_p(\ov\cM_{\ell};\Z)\!\big)\cap
\Psi\!\bigg(\!\!H_{p-2}\big(\ov\cM_{\ell};\Z\big)
\!\oplus\!\bigoplus_{\{J,K\}\in\cP_{\bu}(\ell)}\hspace{-.24in}
H_{p-2}\big(D_{J,K};\Z\big)\!\!\!\bigg)=\{0\}.$$
The restriction of~$\Psi$ to~$H_{p-2}(\ov\cM_{\ell};\Z)$ is injective because
$$\big\{\ff_{\ell+1}\big\}_{\!*}
\big(D_{\{1,\ell+1\},[\ell]-\{1\}}\!\cap\!\Psi(\cdot)\!\big)\!=\!\id\!:
H_{p-2}\big(\ov\cM_{\ell};\Z\big)\lra H_{p-2}\big(\ov\cM_{\ell};\Z\big)\,.$$
Since the divisors $D_{\{1,\ell+1\},[\ell]-\{1\}},\wt{D}_{J,K}^0\!\subset\!\ov\cM_{\ell+1}$
are disjoint for all $\{J,K\}\!\in\!\cP_{\bu}(\ell)$,
$$\big\{\ff_{\ell+1}\big\}_{\!*}
\big(D_{\{1,\ell+1\},[\ell]-\{1\}}\!\cap\!\Psi(\cdot)\!\big)\!=\!0\!:
H_{p-2}\big(D_{J,K};\Z\big)\lra H_{p-2}\big(\ov\cM_{\ell};\Z\big)\,.$$
It follows that 
$$\Psi\big(H_{p-2}(\ov\cM_{\ell};\Z)\!\big)\cap
\Psi\!\bigg(\bigoplus_{\{J,K\}\in\cP_{\bu}(\ell)}\hspace{-.24in}
H_{p-2}\big(D_{J,K};\Z\big)\!\!\!\bigg)=\{0\}.$$
The injectivity of the restriction
$$\Psi\!:\bigoplus_{\{J,K\}\in\cP_{\bu}(\ell)}\hspace{-.24in}
H_{p-2}\big(D_{J,K};\Z\big)\lra \ov\cM_{\ell+1}$$
follows from Lemma~\ref{Chomol_lmm} below and the Poincare Duality for~$D_{J,K}$.
\end{proof}

\noindent
We denote by $\cdot_{D_{J,K}}$ and $\cdot_{\ov\cM_{\ell+1}}$ the homology intersection pairings
on~$D_{J,K}$ and $\ov\cM_{\ell+1}$, respectively.

\begin{lmm}\label{Chomol_lmm}
Suppose $\ell\!\in\!\Z^+$ with $\ell\!\ge\!3$, $p\!\in\!\Z$, 
and \hbox{$\{J,K\},\{J',K'\}\!\in\!\cP_{\bu}(\ell)$}.
If \hbox{$\ka\!\in\!H_{p-2}(D_{J,K};\Z)$} and 
$\ka'\!\in\!H_{2\ell-6-p}(D_{J',K'};\Z)$, then
$$\wt\ka^0\!\cdot_{\ov\cM_{\ell+1}}\!\wt\ka'^+
=\begin{cases}\ka\!\cdot_{D_{J,K}}\!\ka',&\hbox{if}~\{J,K\}\!=\!\{J',K'\};\\
0,&\hbox{if}~\{J,K\}\!\neq\!\{J',K'\}.
\end{cases}$$
\end{lmm}

\begin{proof} If $\{J,K\}\!\neq\!\{J',K'\}$, the divisors 
$D_{J,K},D_{J',K'}\!\subset\!\ov\cM_{\ell}$ are transverse.
We can then choose pseudocycle representatives
\BE{Chomol_e21}f\!:Z\lra D_{J,K} \qquad\hbox{and}\qquad f'\!:Z'\lra D_{J',K'}\EE
for~$\ka$ and~$\ka'$, respectively, 
so that they are transverse in~$\ov\cM_{\ell}$ and thus disjoint.
The pseudocycle representatives
\BE{Chomol_e25}\wt{D}_{J,K}^0\big|_f\!:f^*\wt{D}_{J,K}^0\lra \wt{D}_{J,K}^0\subset \ov\cM_{\ell+1}
\quad\hbox{and}\quad
\wt{D}_{J',K'}^+\big|_{f'}\!:f'^*\wt{D}_{J',K'}^+\lra \wt{D}_{J',K'}^+\subset \ov\cM_{\ell+1}\EE
for~$\wt\ka^0$ and~$\wt\ka'^+$, respectively, 
are then disjoint as well.
This implies the claim in the second case.\\

\noindent
If $\{J,K\}\!=\!\{J',K'\}$, we can choose pseudocycle representatives as in~\eref{Chomol_e21}
so that they are transverse in~$D_{J,K}$ and thus intersect at a finite numbers
of signed pairs of points \hbox{$(z_r,z_r')\!\in\!Z\!\times\!Z'$}
with $f(z_r)\!=\!f'(z_r')$.
The induced pseudocycles~\eref{Chomol_e25} then intersect at the points 
$$\big(\!(z_r,\wt\cC_r),(z_r',\wt\cC_r)\!\big)\in 
f^*\wt{D}_{J,K}^0\!\times\!f'^*\wt{D}_{J',K'}^+
\quad\hbox{with}\quad \{\wt\cC_r\}
=\big(\wt{D}_{J,K}^0\!\cap\!\wt{D}_{J',K'}^+\big)\!\big|_{f(z_r)}\,.$$
Since the isomorphism on the normal bundles
$$\cN_{\wt{D}_{J,K}^0}\!\big(\wt{D}_{J,K}^0\!\cap\!\wt{D}_{J',K'}^+\big)\!\big|_{\wt\cC_r}
\!\oplus\!
\cN_{\wt{D}_{J',K'}^0}\!\big(\wt{D}_{J,K}^0\!\cap\!\wt{D}_{J',K'}^+\big)\!\big|_{\wt\cC_r}
\lra
\cN_{\ov\cM_{\ell+1}}\!\big(\wt{D}_{J,K}^0\!\cap\!\wt{D}_{J',K'}^+\big)\!\big|_{\wt\cC_r}$$
induced by the inclusions of the tangent bundles is orientation-preserving and 
the two summands on the left-hand side are of even (real) ranks,
the sign of the intersection point $\big(\!(z_r,\wt\cC_r),(z_r',\wt\cC_r)\!\big)$
is the same as the sign of~$(z_r,z_r')$.
This implies the claim in the first case.
\end{proof}

\subsection{The complex case: algebra}
\label{Calg_subs}

\noindent
For a finite set~$S$, let 
$$\cP(S)=\big\{\!\{J,K\}\!\!:S\!=\!J\!\sqcup\!K\big\}, \quad
\cP_{\bu}(S)=\big\{\!\{J,K\}\!\in\!\cP(S)\!\!:2\!\le\!|J|\!\le\!|S|\!-\!2\big\}.$$
Denote by $\cI_S\!\subset\!\Z[(D_{J,K})_{\{J,K\}\in\cP_{\bu}(S)}]$
the ideal generated by the left-hand sides of~\eref{cM04rel_e1a} and~\eref{cM04rel_e2}
with~$\cP_{\bu}(\ell)$ replaced by~$\cP_{\bu}(S)$ and $a,b,c,d\!\in\!S$ distinct.
Let
$$\cR_S=\Z\big[(D_{J,K})_{\{J,K\}\in\cP_{\bu}(S)}\big]\big/\cI_S\,.$$
Thus, $\cI_{\ell}\!=\!\cI_{[\ell]}$ and the left-hand side of~\eref{HCthm_e}
is the ring $\cR_{\ell}\!\equiv\!\cR_{[\ell]}$.\\

\noindent
Let $\ell\!\in\!\Z^+$ with $\ell\!\ge\!3$.
Define a $\Z$-algebra homomorphism
\BE{FdfnC_e}F\!:\Z\big[(D_{J,K})_{\{J,K\}\in\cP_{\bu}(\ell)}\big]\lra
\Z\big[(D_{J,K})_{\{J,K\}\in\cP_{\bu}(\ell+1)}\big], ~~
F\big(D_{J,K}\big)=\wt{D}_{J,K}^+\!+\!\wt{D}_{J,K}^0\,,\EE
with  $\wt{D}_{J,K}^+,\wt{D}_{J,K}^0\!\in\!\Z[(D_{J,K})_{\{J,K\}\in\cP_{\bu}(\ell+1)}]$
defined as in~\eref{wtDdfn_e}.
For $\{J,K\}\!\in\!\cP_{\bu}(\ell)$, define a $\Z$-algebra homomorphism
\begin{gather*}
F_J\!:\Z\big[(D_{J',K'})_{\{J',K'\}\in\cP_{\bu}(\{\nod\}\sqcup K)}\big]
\lra \Z\big[(D_{J',K'})_{\{J',K'\}\in\cP_{\bu}(\ell)}\big] \qquad\hbox{by}\\
F_J\big(D_{J',K'}\big)=D_{J\cup(J'-\{\nod\}),K'} \quad\hbox{if}~~J'\!\ni\!\nod.
\end{gather*}
We denote by 
$$F_{J,K}\!:\Z\big[(D_{J',K'})_{\{J',K'\}\in\cP_{\bu}(\{\nod\}\sqcup J)}\big]
\!\otimes_{\Z}\!
\Z\big[(D_{J',K'})_{\{J',K'\}\in\cP_{\bu}(\{\nod\}\sqcup K)}\big]
\lra \Z\big[(D_{J',K'})_{\{J',K'\}\in\cP_{\bu}(\ell)}\big]$$
the $\Z$-algebra homomorphism induced by~$F_K$ and~$F_J$.

\begin{prp}\label{Calg_prp}
For every $\ell\!\in\!\Z^+$ with $\ell\!\ge\!3$, the homomorphism
\BE{Calg_e}\begin{split}
&\hspace{.7in}
\Phi\!:\cR_{\ell}\!\oplus\!\cR_{\ell}\!\oplus\!
\bigoplus_{\{J,K\}\in\cP_{\bu}(\ell)}\hspace{-.24in}
\cR_{\{\nod\}\sqcup J}\!\otimes_{\Z}\!\cR_{\{\nod\}\sqcup K}\lra\cR_{\ell+1}\,,\\
&\Phi\big([\ka_0],[\ka],([\ka_{J,K}])_{\{J,K\}\in\cP_{\bu}(\ell)}\big)
=D_{\{1,\ell+1\},[\ell]-\{1\}}F(\ka_0)\!+\!F(\ka)\\
&\hspace{2.4in}
+\sum_{\{J,K\}\in\cP_{\bu}(\ell)}\hspace{-.24in}
\wt{D}_{J,K}^0F\big(F_{J,K}(\ka_{J,K})\!\big)+\cI_{\ell+1}\,,
\end{split}\EE
of $\Z$-modules is well-defined and surjective.
\end{prp}

\noindent
With $\Psi$ denoting the composition of the homomorphism~\eref{Chomol_e}
with the Poincare duals on both sides, the diagram
\BE{HCdiag_e}\begin{split}
\xymatrix{\cR_{\ell}\!\oplus\!\cR_{\ell}\!\oplus\!
\bigoplus\limits_{\{J,K\}\in\cP_{\bu}(\ell)}\hspace{-.24in}
\cR_{\{\nod\}\sqcup J}\!\otimes_{\Z}\!\cR_{\{\nod\}\sqcup K}
\ar[r]^>>>>>>>>>>>>>>>{\Phi}\ar[d]_{\eref{HCthm_e}}
&\cR_{\ell+1}\ar[d]^{\eref{HCthm_e}} \\
H^*\big(\ov\cM_{\ell};\Z\big)\!\oplus\!H^*\big(\ov\cM_{\ell};\Z\big)
\!\oplus\!\bigoplus\limits_{\{J,K\}\in\cP_{\bu}(\ell)}\hspace{-.24in}
H^*\big(D_{J,K};\Z\big)\ar[r]^>>>>>>{\Psi}
& H^*\big(\ov\cM_{\ell+1};\Z\big)}
\end{split}\EE
commutes. 
The commutativity on the $\{J,K\}$ summand in the domain of~$\Phi$ follows from 
the second identity in~\eref{CYvrvrcap_e} and 
$$\big\{\wt\io_{J,K}^0\big\}_{\!*}
\Big(\PD_{\wt{D}_{J,K}^0}\!\big(\!\{\wt\io_{J,K}^0\}^{\!*}\wt\ka\big)\!\!\Big)
=\wt\ka\!\cap\!\wt{D}_{J,K}^0\in H_*\big(\ov\cM_{\ell+1};\Z\big)
\quad\forall~\wt\ka\!\in\!H^*\big(\ov\cM_{\ell+1};\Z\big).$$
By Proposition~\ref{Calg_prp} and the generation property~\ref{HCgen_it} on 
page~\pageref{HCgen_it} with~$\ell$ replaced by~$\ell\!+\!1$, 
the top and right homomorphisms in this diagram are surjective.
By Theorem~\ref{HC_thm} with~$\ell$ replaced by each $\ell'\!\le\!\ell$
and Proposition~\ref{Chomol_prp}, 
the left and bottom homomorphisms are isomorphisms.
Along with the commutativity of the diagram~\eref{HCdiag_e},
this implies that the top and right homomorphisms are isomorphisms as well
and establishes the injectivity property~\ref{HCrel_it} on 
page~\pageref{HCrel_it} with~$\ell$ replaced by~$\ell\!+\!1$.

\begin{proof}[{\bf{\emph{Proof of Proposition~\ref{Calg_prp}: well-defined}}}]
Since $F$ sends the left-hand side of~\eref{cM04rel_e1a} to
the sum of four analogous terms in its target
and the  left-hand side of~\eref{cM04rel_e2} to its analogue in the target,
\hbox{$F(\cI_{\ell})\!\subset\!\cI_{\ell+1}$}.
Thus, $\Phi$ is well-defined on the two $\cR_{\ell}$-summands in~\eref{Calg_e}.
Below we denote by~$\eref{cM04rel_e2}_{abcd}^S$ 
the left-hand side of~\eref{cM04rel_e2} with~$[\ell]$ replaced by a finite set~$S$.
For $s_0\!\in\!S$ and \hbox{$a,b,c,d\!\in\!S\!-\!\{s_0\}$} distinct,
\BE{cM04rel_e2gen}  \eref{cM04rel_e2}_{abcd}^S=
\eref{cM04rel_e2}_{s_0adc}^S- \eref{cM04rel_e2}_{s_0adb}^S\,.\EE
In particular, the collection of 
the left-hand sides of~\eref{cM04rel_e2} with~$[\ell]$ replaced by a finite set~$S$
is linearly generated by its subcollection with $a\!=\!s_0$.\\

\noindent
Let $\{J,K\}\!\in\!\cP_{\bu}(\ell)$. 
Since~$F_J$ maps the left-hand side of~\eref{cM04rel_e1a}
in its domain to an analogous product in the target,
$F\!\circ\!F_{J,K}$ sends this product tensored with~1 to~$\cI_{\ell+1}$.
If $b,c,d\!\in\!K$ distinct and $a\!\in\!J$,
$$F_J\Big(\eref{cM04rel_e2}_{(\nod)bcd}^{\{\nod\}\sqcup K}\Big)
-\eref{cM04rel_e2}_{abcd}^{[\ell]}
=-\!\!\!\sum_{\begin{subarray}{c} \{J',K'\}\in\cP_{\bu}(\ell)\\
a,b\in J',\,c,d\in K'\\ \{J,K\}\not\cap\,\{J',K'\} \end{subarray}}\hspace{-.32in}D_{J',K'}
+\!\!\! \sum_{\begin{subarray}{c} \{J',K'\}\in\cP_{\bu}(\ell) \\
a,c\in J',\,b,d\in K'\\ \{J,K\}\not\cap\,\{J',K'\}\end{subarray}}\hspace{-.32in}D_{J',K'}\,.$$
The product of~$\wt{D}_{J,K}^0$ with the image of the right-hand side
of the last equation under~$F$ lies in the linear span
of the left-hand sides of~\eref{cM04rel_e1a} with~$\ell$ replaced by~$\ell\!+\!1$.
Along with~\eref{cM04rel_e2gen}, this implies~that
$$\Phi\big(
\Z\big[(D_{J',K'})_{\{J',K'\}\in\cP_{\bu}(\{\nod\}\sqcup J)}\big]
\!\otimes_{\Z}\!\cI_{\{\nod\}\sqcup K}\big)\subset \cI_{\ell+1}\,.$$
Similarly,
$$\Phi\big(\cI_{\{\nod\}\sqcup J}\!\otimes_{\Z}\!
\Z\big[(D_{J',K'})_{\{J',K'\}\in\cP_{\bu}(\{\nod\}\sqcup K)}\big]\big)
\subset \cI_{\ell+1}\,.$$
Thus, $\Phi$ is well-defined on each 
$\cR_{\{\nod\}\sqcup J}\!\otimes_{\Z}\!\cR_{\{\nod\}\sqcup K}$-summand.
\end{proof}

\begin{proof}[{\bf{\emph{Proof of Proposition~\ref{Calg_prp}: surjective}}}]
The $\Z$-algebra $\cR_{\ell+1}$ is algebraically generated by the equivalence classes of 
\BE{Calg_e15}F(D_{J,K}),\wt{D}_{J,K}^0,D_{\{i,\ell+1\},[\ell]-\{i\}}\in
\Z\big[(D_{J',K'})_{\{J',K'\}\in\cP_{\bu}(\ell+1)}\big] 
~~\hbox{with}~\{J,K\}\!\in\!\cP_{\bu}(\ell),\,i\!\in\![\ell].\EE
Since the left-hand side of~\eref{cM04rel_e2} with $(\ell,a,b,d)$
replaced by $(\ell\!+\!1,1,\ell\!+\!1,i)$, $i\!\neq\!1$, and any $c\!\in\![\ell]$ distinct from~$1,i$
lies in~$\cI_{\ell+1}$, it is enough to take $i\!=\!1$ in the above list of algebraic generators.
Since $D_{\{1,\ell+1\},[\ell]-\{1\}}\wt{D}_{J,K}^0$ lies in~$\cI_{\ell+1}$ as~well,
it also follows that 
$\cR_{\ell+1}$ is {\it linearly} generated by arbitrary products of the elements of the first two types
in~\eref{Calg_e15} and arbitrary products of the elements of the first type
with a {\it single} factor of $D_{\{1,\ell+1\},[\ell]-\{1\}}$.
Since $F$ is a ring homomorphism, the algebraic span of the elements of the first type
in~\eref{Calg_e15} is~$\Im\,F$.\\

\noindent
Suppose $\{J,K\},\{J',K'\}\!\in\!\cP_{\bu}(\ell)$ and $1\!\in\!J,J'$.
If $\{J,K\}\!\not\!\cap\{J',K'\}$, then 
$$\wt{D}_{J,K}^{\bu}\wt{D}_{J',K'}^{\circ}\in\cI_{\ell+1}\quad\forall~\bu,\circ=+,0.$$
If \hbox{$J\!\supsetneq\!J'$} (and thus $K\!\subsetneq\!K'$), then
$$\wt{D}_{J,K}^0\wt{D}_{J',K'}^0-\wt{D}_{J,K}^0F\big(D_{J',K'}\big)
=-\wt{D}_{J,K}^0\wt{D}_{J',K'}^+\in\cI_{\ell+1}\,.$$
If $b\!\in\!J\!-\!\{1\}$ and $d\!\in\!K$, then
\begin{equation*}\begin{split}
\wt{D}_{J,K}^+\wt{D}_{J,K}^0-\!F\big(D_{J,K}\big)D_{\{1,\ell+1\},[\ell]-\{1\}}
-\wt{D}_{J,K}^+\eref{cM04rel_e2}_{1b(\ell+1)d}^{[\ell+1]}&\\
+\wt{D}_{J,K}^+\bigg(\!\!
\sum_{\begin{subarray}{c} \{J',K'\}\in\cP_{\bu}(\ell)\\
1,b\in J'\neq J,\,d\in K'  \end{subarray}}\hspace{-.32in}\wt{D}_{J',K'}^0
-\!\!\! \sum_{\begin{subarray}{c} \{J',K'\}\in\cP_{\bu}(\ell) \\
1\in J'\neq J,\,b,d\in K' \end{subarray}}\hspace{-.32in}\wt{D}_{J',K'}^+\!\!\bigg)
&=-\wt{D}_{J,K}^0D_{\{1,\ell+1\},[\ell]-\{1\}}\in\cI_{\ell+1}.
\end{split}\end{equation*}
By~\eref{FdfnC_e},
the last product on the left-hand side lies in the $(\Im\,F)$-submodule
linearly generated by  
of~1, $\wt{D}_{J',K'}^0$ with \hbox{$\{J',K'\}\!\in\!\cP_{\bu}(\ell)$},
and $\wt{D}_{J,K}^0\wt{D}_{J',K'}^0$ with \hbox{$\{J',K'\}\!\in\!\cP_{\bu}(\ell)$} {\it distinct} 
from~$\{J,K\}$.
By the preceding discussion,
$\cR_{\ell+1}$ is thus a module over $\Im\,F$ linearly generated by 
the equivalence classes of~1, $D_{\{1,\ell+1\},[\ell]-\{1\}}$, and~$\wt{D}_{J,K}^0$.
Along with Lemma~\ref{Calg_lmm} below, 
this establishes the surjectivity of the homomorphism~\eref{Calg_e}.
\end{proof}

\begin{lmm}\label{Calg_lmm}
Let $\ell\!\in\!\Z^+$ with $\ell\!\ge\!3$.
For every $\{J,K\}\!\in\!\cP_{\bu}(\ell)$,
$$\wt{D}_{J,K}^0(\Im\,F)\subset\wt{D}_{J,K}^0\big(\Im\{F\!\circ\!F_{J,K}\}\!\big)+\cI_{\ell+1}\,.$$
\end{lmm}

\begin{proof}
By~\eref{cM04rel_e2},
the domain of $F$ is algebraically generated by $D_{J',K'}$ with \hbox{$\{J',K'\}\!\in\!\cP_{\bu}(\ell)$}
distinct from~$\{J,K\}$ and~$\cI_{\ell}$.
If \hbox{$\{J,K\}\!\not\!\cap\{J',K'\}$}, then 
$\wt{D}_{J,K}^0F(D_{J',K'})\!\in\!\cI_{\ell+1}$ by~\eref{cM04rel_e1a}.
If $J'\!\subsetneq\!J$,
$$D_{J',K'}=F_{J,K}\big(D_{J',\{\nod\}\sqcup(J-J')}\!\otimes\!1\big)$$
by the definition of~$F_{J,K}$.
Since $F(\cI_{\ell})\!\subset\!\cI_{\ell+1}$, the claim follows.
\end{proof}

\subsection{The real case: orientations}
\label{Rhomol_subs1}

\noindent
We identify the one-dimensional moduli space 
$\R\ov\cM_{0,2}\!\subset\!\ov\cM_4$ with $\R\P^1\!\subset\!\C\P^1$ via the cross ratio 
\BE{cM02ident_e}
\CR\!:\R\ov\cM_{0,2}\lra \R\P^1, ~~
\CR\big([\cC]\big)= 
\frac{z_2^+(\cC)\!-\!z_1^-(\cC)}{z_2^-(\cC)\!-\!z_1^-(\cC)}:
\frac{z_2^+(\cC)\!-\!z_1^+(\cC)}{z_2^-(\cC)\!-\!z_1^+(\cC)}
=\pm\frac{|1\!-\!z_1^+(\cC)/z_2^-(\cC)|^2}{|z_1^+(\cC)\!-\!z_2^+(\cC)|^2}\,,\EE
where $z_i^{\pm}(\cC)\!\in\!\C\P^1$ is the marked point of the curve~$\cC$ 
indexed by~$(i,\pm)$ if $\cC$ is smooth ($\CR$ extends smoothly over all of~$\R\ov\cM_{0,2}$).
The cross ratio~$\CR$ takes positive values on the subspace 
\hbox{$\cM_{0,2}^{\tau}\!\subset\!\R\ov\cM_{0,2}$}
of smooth curves with nonempty real locus, as indicated in Figure~\ref{M02_fig},
and negative values on the subspace $\cM_{0,2}^{\eta}\!\subset\!\R\ov\cM_{0,2}$
of smooth curves with~empty real locus.
The identification~\eref{cM02ident_e} determines an orientation~$\fo_{0,2}$ on~$\R\ov\cM_{0,2}$, 
which is the {\it opposite} of the orientation induced by~\cite[(3.1)]{RealEnum} 
and~\cite[(1.12)]{RealGWsII}.\\

\noindent
For $\ell\!\ge\!2$, the restriction of the forgetful morphism~$\ff_{(\ell+1)^{\pm}}^{\R}$ 
in~\eref{ffRdfn_e} to the preimage of the subspace $\R\cM_{0,\ell}\!\subset\!\R\ov\cM_{0,\ell}$
of smooth curves is a $\C\P^1$-fiber bundle.
Thus, it induces a short exact sequence 
\BE{cMorientR_e}
0\lra \ker\nd\ff_{(\ell+1)^{\pm}}^{\R}\lra T\big(\R\ov\cM_{0,\ell+1}\big)
\xlra{\nd\ff_{(\ell+1)^{\pm}}^{\R}} 
\big\{\ff_{(\ell+1)^{\pm}}^{\R}\big\}^{\!*}T\big(\R\ov\cM_{0,\ell}\big)\lra0\EE
of vector bundles over $\{\ff_{(\ell+1)^{\pm}}^{\R}\}^{-1}(\R\cM_{0,\ell})$.
The vertical tangent bundle above is canonically oriented by the complex orientation of 
the fiber~$\C\P^1$ at the marked point~$z_{\ell+1}^+$.
For $\ell\!\ge\!2$, we inductively define an orientation~$\fo_{0,\ell+1}$ on~$\R\ov\cM_{0,\ell+1}$
so that the short exact sequence~\eref{cMorientR_e} respects the canonical orientation
on the vertical tangent bundle,
the orientation~$\fo_{0,\ell+1}$ on~$\R\ov\cM_{0,\ell+1}$,
and the orientation~$\fo_{0,\ell}$ on~$\R\ov\cM_{0,\ell}$.
The orientation~$\fo_{0,\ell}$ is preserved by the automorphism of~$\R\ov\cM_{0,\ell}$ 
interchanging the conjugate pairs~$(z_i^+,z_i^-)$ and~$(z_j^+,z_j^-)$ of marked points
of each marked curve.
Thus, it determines an orientation~$\fo_{0,S}$ on~$\R\ov\cM_{0,S}$ for any set~$S$ with 
$|S|\!=\!\ell$.\\

\begin{figure}  
\begin{pspicture}(-3.2,-1.5)(10,2)
\psset{unit=.4cm}
\pscircle[linewidth=.1](-3,0){1}\pscircle[linewidth=.1](-3,2){1}
\pscircle*(-3,1){.2}\pscircle*(-3.7,-.7){.2}\pscircle*(-2.3,-.7){.2}
\pscircle*(-3.7,2.7){.2}\pscircle*(-2.3,2.7){.2}
\rput(-3.6,-1.5){\sm{$1^-$}}\rput(-1.9,-1.5){\sm{$2^+$}}
\rput(-3.6,3.5){\sm{$1^+$}}\rput(-1.9,3.5){\sm{$2^-$}}
\psline[linewidth=.1,linestyle=dotted](-1,1)(1,1)
\pscircle[linewidth=.1](3.5,1){1.5}
\psarc[linewidth=.05](3.5,3.6){3}{240}{300}
\psarc[linewidth=.05,linestyle=dashed](3.5,-1.6){3}{60}{120}
\pscircle*(2.44,-.06){.2}\pscircle*(4.56,-.06){.2}\pscircle*(2.44,2.06){.2}\pscircle*(4.56,2.06){.2}
\rput(2.3,-.9){\sm{$1^-$}}\rput(4.7,-.9){\sm{$2^+$}}
\rput(2.3,2.9){\sm{$1^+$}}\rput(4.7,2.9){\sm{$2^-$}}
\psline[linewidth=.1,linestyle=dotted](6.5,1)(8.5,1)
\pscircle[linewidth=.1](11,1){1}\pscircle[linewidth=.1](13,1){1}
\pscircle*(12,1){.2}\pscircle*(10.3,.3){.2}\pscircle*(13.7,.3){.2}
\pscircle*(10.3,1.7){.2}\pscircle*(13.7,1.7){.2}
\rput(10.16,-.54){\sm{$1^-$}}\rput(13.84,-.54){\sm{$2^-$}}
\rput(10.16,2.54){\sm{$1^+$}}\rput(13.84,2.54){\sm{$2^+$}}
\psline[linewidth=.1,linestyle=dotted](15.5,1)(17.5,1)
\pscircle[linewidth=.1](20.5,1){1.5}
\psarc[linewidth=.05](20.5,3.6){3}{240}{300}
\psarc[linewidth=.05,linestyle=dashed](20.5,-1.6){3}{60}{120}
\pscircle*(19.44,-.06){.2}\pscircle*(21.56,-.06){.2}\pscircle*(19.44,2.06){.2}\pscircle*(21.56,2.06){.2}
\rput(19.3,-.9){\sm{$1^-$}}\rput(21.7,-.9){\sm{$2^-$}}
\rput(19.3,2.9){\sm{$1^+$}}\rput(21.7,2.9){\sm{$2^+$}}
\psline[linewidth=.1,linestyle=dotted](23,1)(25,1)
\pscircle[linewidth=.1](27,0){1}\pscircle[linewidth=.1](27,2){1}
\pscircle*(27,1){.2}\pscircle*(26.3,-.7){.2}\pscircle*(27.7,-.7){.2}
\pscircle*(26.3,2.7){.2}\pscircle*(27.7,2.7){.2}
\rput(26.4,-1.5){\sm{$1^-$}}\rput(28.1,-1.5){\sm{$2^-$}}
\rput(26.4,3.5){\sm{$1^+$}}\rput(28.1,3.5){\sm{$2^+$}}
\psline[linewidth=.1](-3,-3)(27,-3)
\pscircle*(-3,-3){.2}\pscircle*(12,-3){.2}\pscircle*(27,-3){.2}
\rput(-3,-3.8){0}\rput(12,-3.8){1}\rput(27,-3.8){$\i$}
\rput(-3,-3.8){0}\rput(12,-3.8){1}\rput(27,-3.8){$\i$}
\end{pspicture}
\caption{The structure and orientation of the Deligne-Mumford compactification 
$\ov\cM_{0,2}^{\tau}\!\subset\!\R\ov\cM_{0,2}$ of~$\cM_{0,2}^{\tau}$.}
\label{M02_fig}
\end{figure}

\noindent
Let $\ell\!\in\!\Z^+$ with $\ell\!\ge\!2$. 
For $\{J,K\}\!\in\!\cP(\ell)$ with $J\!\in\!1$,
we take the first identification in~\eref{Rboundident_e} to be given by the projection
\BE{RboundidentE_e}\R E_{J,K} \stackrel{\approx}{\lra} \ov\cM_{\{\nod\}\sqcup[\ell]} \EE
to the minimal connected union~$\cC^+$ of the irreducible components of each curve~$\cC$
carrying the marked points $z_j^+(\cC)$ with $j\!\in\!J$,
the marked points $z_k^-(\cC)$ with $k\!\in\!K$, and the real node of~$\cC$
(which is sent to the marked point~$z_{\nod}(\cC^+)$ of~$\cC^+$).
We define the orientation on~$\R E_{J,K}$ to be induced via 
the identification~\eref{RboundidentE_e} from the complex orientation
of~$\ov\cM_{\{\nod\}\sqcup[\ell]}$ if $|K|\!\in\!2\Z$
and from the opposite orientation of~$\ov\cM_{\{\nod\}\sqcup[\ell]}$ if $|K|\!\not\in\!2\Z$.\\

\noindent 
Thus, the first (resp.~second) diagram in Figure~\ref{TopolRel_fig} on page~\pageref{TopolRel_fig}
represents a minus (resp.~plus) point in~$\R\ov\cM_{0,2}$;
this yields the $\ell\!=\!2$ case of~\eref{RcM04rel_e2a}.
Furthermore, the short exact sequences
\BE{ffEses_e}\begin{split}
0&\lra \ker\nd\ff_{(\ell+1)^{\pm}}^{\R}\lra T\big(\R E_{J\cup\{\ell+1\},K}\big)
\xlra{\nd\ff_{(\ell+1)^{\pm}}^{\R}} 
\big\{\ff_{(\ell+1)^{\pm}}^{\R}\big\}^{\!*}T\big(\R E_{J,K}\big)\lra0,\\
0&\lra \ker\nd\ff_{(\ell+1)^{\pm}}^{\R}\lra T\big(\R E_{J,K\cup\{\ell+1\}}\big)
\xlra{\nd\ff_{(\ell+1)^{\pm}}^{\R}} 
\big\{\ff_{(\ell+1)^{\pm}}^{\R}\big\}^{\!*}T\big(\R E_{J,K}\big)\lra0,
\end{split}\EE
of vector bundles over the main strata of $\R E_{J\cup\{\ell+1\},K}$ and $\R E_{J,K\cup\{\ell+1\}}$
respect the canonical orientation on the vertical tangent bundle
and the specified orientations on~$\R E_{J,K}$, $\R E_{J\cup\{\ell+1\},K}$, 
and~$\R E_{J,K\cup\{\ell+1\}}$.
Since
\begin{gather*}
\big\{\ff_{(\ell+1)^{\pm}}^{\R}\big\}^{\!-1}\big(\R E_{J,K}\big)
=\R E_{J\cup\{\ell+1\},K}\!\cup\!\R D_{\{\ell+1\};J,K}\!\cup\!\R E_{J,K\cup\{\ell+1\}}\,,\\
\dim\R D_{\{\ell+1\};J,K}<\dim\R E_{J\cup\{\ell+1\},K},\dim\R E_{J,K\cup\{\ell+1\}}\,,
\end{gather*}
and~\eref{cMorientR_e} respects the relevant orientations as well, it follows that 
$$\big\{\ff_{(\ell+1)^{\pm}}^{\R}\big\}^{\!*}\R E_{J,K}
=\R E_{J\cup\{\ell+1\},K}\!+\!\R E_{J,K\cup\{\ell+1\}}
\in H^1\big(\R\ov\cM_{0,\ell+1};\Z\big)\,.$$
Thus, the $\ell\!\ge\!3$ cases of~\eref{RcM04rel_e2a} follow inductively
from the $\ell\!=\!2$ case of this relation.\\

\noindent
If $\{J,K\}\!\in\!\cP(\ell)$, interchanging the roles of~$J$ and~$K$ in~\eref{RboundidentE_e}
preserves the orientation on the left-hand side of~\eref{RboundidentE_e} induced by the complex orientation
of the right-hand side via the identification in~\eref{RboundidentE_e}
if and only~if
$$\dim_{\C}\!\ov\cM_{\{\nod\}\sqcup[\ell]}\!=\!\ell\!-\!2\in2\Z.$$
Taking into account the orientation correction by the parity of~$|K|$ and that \hbox{$|J|\!+\!|K|\!=\!\ell$}, 
we conclude that  the orientation on~$\R E_{J,K}$
constructed above does not depend on the choice of \hbox{$J\!\in\!\{J,K\}$};
this is reflected in the first claim of Lemma~\ref{DEorient_lmm} below.\\

\noindent
For $(I,\{J,K\})\!\in\!\wt\cP_{\bu}(\ell)$ with $\min(J\!\cup\!K)\!\in\!J$,
we take the second identification in~\eref{Rboundident_e} to be given by the projections
\BE{RboundidentD_e}
\R D_{I;J,K} \stackrel{\approx}{\lra} \R\ov\cM_{0,\{\nod\}\sqcup I}
\!\times\!\ov\cM_{\{\nod\}\sqcup(J\cup K)}\EE
to the minimal connected union~$\cC_0$ of the irreducible components of each curve~$\cC$
carrying the marked points $z_i^+(\cC),z_i^-(\cC)$ with $i\!\in\!I$
and to the minimal connected union~$\cC^+$ of the irreducible components of~$\cC$
carrying the marked points $z_j^+(\cC)$ with $j\!\in\!J$ and
the marked points $z_k^-(\cC)$ with $k\!\in\!K$.
The former inherits an involution from~$\cC$ which interchanges the two nodes
separating~$\cC_0$ from the remaining components of~$\cC$.
These two nodes thus form a conjugate pair $(z_{\nod}^+(\cC_0),z_{\nod}^-(\cC_0))$ of marked points
of~$\cC_0$, with $z_{\nod}^+(\cC_0)$ separating~$\cC_0$ from~$\cC^+$.
We take the marked point~$z_{\nod}(\cC^+)$ of~$\cC^+$ to be the node separating it from~$\cC_0$.
We define the orientation on~$\R D_{I;J,K}$ to be induced via 
the identification in~\eref{RboundidentD_e} from 
the orientation~$\fo_{0,\{\nod\}\sqcup I}$ on~$\R\ov\cM_{0,\{\nod\}\sqcup I}$
and the complex orientation of~$\ov\cM_{\{\nod\}\sqcup(J\cup K)}$ if $|K|\!\in\!2\Z$
and from the opposite orientation of~$\ov\cM_{\{\nod\}\sqcup(J\cup K)}$ if $|K|\!\not\in\!2\Z$.\\

\noindent
The bottom diagram of Figure~\ref{TopolRel_fig} indicates the resulting orientation
on the open subspace~$\R D_{I;J,K}^{\tau}$ of each of the six  boundary divisors
$\R D_{I;J,K}\!\subset\!\R\ov\cM_{0,3}$ 
which ($\R D_{I;J,K}^{\tau}$, that is) parametrizes smooth curves with nonempty real locus.
The above orientations on these six boundary divisors are opposite of those defined
in Section~3 of~\cite{RealEnum}.
Thus, the three $\ell\!=\!3$ cases of~\eref{RcM04rel_e2b}
are equivalent to the conclusion of the first statement of Remark~3.4 in~\cite{RealEnum}.
For any $(I',\{J',K'\})\!\in\!\wt\cP_{\bu}(\ell+1)$ with $I'\!\supset\!I$,
$J'\!\supset\!J$, and $K'\!\supset\!K$, the short exact sequence
\BE{ffDses_e}0\lra \ker\nd\ff_{(\ell+1)^{\pm}}^{\R}\lra T\big(\R D_{I';J',K'}\big)
\xlra{\nd\ff_{(\ell+1)^{\pm}}^{\R}} 
\big\{\ff_{(\ell+1)^{\pm}}^{\R}\big\}^{\!*}T\big(\R D_{I;J,K}\big)\lra0\EE
of vector bundles over the main strata of $\R D_{I';J',K'}$
respects the canonical orientation on the vertical tangent bundle
and the specified orientations on~$\R D_{I';J',K'}$ and~$\R D_{I;J,K}$.
Since
$$\big\{\ff_{(\ell+1)^{\pm}}^{\R}\big\}^{\!-1}\big(\R D_{I;J,K}\big)
=\R D_{I\cup\{\ell+1\};J,K}\!\cup\!\R D_{I;J\cup\{\ell+1\},K}\!\cup\!\R D_{I;J,K\cup\{\ell+1\}}$$
and~\eref{cMorientR_e} respects the relevant orientations as well, it follows that 
$$\big\{\ff_{(\ell+1)^{\pm}}^{\R}\big\}^{\!*}\R D_{I;J,K}
=\R D_{I\cup\{\ell+1\};J,K}\!+\!\R D_{I;J\cup\{\ell+1\},K}\!+\!\R D_{I;J,K\cup\{\ell+1\}}
\in H^2\big(\R\ov\cM_{0,\ell+1};\Z\big)\,.$$
Thus, the $\ell\!\ge\!4$ cases of~\eref{RcM04rel_e2b} with $\{a,b,c\}\!=\![3]$ 
follow inductively from the $\ell\!=\!3$ cases of this relation.\\

\noindent
If \hbox{$(I,\{J,K\})\!\in\!\wt\cP_{\bu}(\ell)$}, 
interchanging the roles of~$J$ and $K$ in~\eref{RboundidentD_e}
reverses the orientation of the first factor on the right-hand side
(this interchanges the marked points~$z_{\nod}^+(\cC_0)$ and~$z_{\nod}^-(\cC_0)$).
Thus,  interchanging the roles of~$J$ and $K$ in~\eref{RboundidentD_e}
preserves the orientation on the left-hand side of~\eref{RboundidentD_e} induced 
by the orientation of the right-hand side via the identification in~\eref{RboundidentD_e}
if and only~if
$$1\!+\!\dim_{\C}\!\ov\cM_{\{\nod\}\sqcup(J\cup K)}\!=\!|J|\!+\!|K|\!-\!1\in2\Z.$$
Taking into account the orientation correction by the parity of~$|K|$, 
we conclude that interchanging the roles of~$J$ and $K$ in~\eref{RboundidentD_e}
reverses the orientation on~$\R D_{I;J,K}$ constructed above;
this is reflected in the second claim of Lemma~\ref{DEorient_lmm}.
This also implies that the automorphism of~$\R\ov\cM_{0,\ell}$ interchanging two pairs of conjugate
marked points and sending a boundary divisor~$\R D_{I;J,K}$ to a boundary divisor~$\R D_{I';J',K'}$
respects the orientations of these divisors if and only if $\ep_{J,K}\!=\!\ep_{J',K'}$.
Thus, \eref{RcM04rel_e2b} follows from its case with $(a,b,c)\!=\!(1,2,3)$.\\

\noindent
For a finite set $S$ and $i\!\in\!S$, denote by
$$L_i\lra \ov\cM_S \qquad\hbox{and}\qquad  L_i^+\lra \R\ov\cM_{0,S}$$
the universal tangent line bundle for the $i$-th marked point and 
the universal tangent line bundle for the $+$~marked point in the $i$-th conjugate pair, respectively.
With the assumptions as in~\eref{RboundidentE_e} and~\eref{RboundidentD_e}, let 
\begin{equation*}\begin{split}
\cN_{\nod;[\ell]}=L_{\nod}\!\otimes_{\C}\!\ov{L_{\nod}}&\lra \ov\cM_{\{\nod\}\sqcup[\ell]} 
\quad\hbox{and}\\
\cN_{\nod;[\ell],I}=
\pi_1^*L_{\nod}^+\!\otimes_{\C}\!\pi_2^*L_{\nod}&\lra\R\ov\cM_{0,\{\nod\}\sqcup I}
\!\times\!\ov\cM_{\{\nod\}\sqcup(J\cup K)}\,,
\end{split}\end{equation*}
where 
$$\pi_1,\pi_2\!:\R D_{I;J,K} \lra \R\ov\cM_{0,\{\nod\}\sqcup I},\ov\cM_{\{\nod\}\sqcup(J\cup K)}$$
are the two component projections.
The identification~\eref{RboundidentE_e} canonically lifts to an identification
\BE{RboundidentE_e2} 
\cN_{\R\ov\cM_{0,\ell}}(\R E_{J,K}) \stackrel{\approx}{\lra} 
\cN_{\nod;[\ell]}^{\,\R}
\equiv \big\{v\!\otimes\!w\!\in\!\cN_{\nod;[\ell]}\!:v\!\otimes\!w\!=\!w\!\otimes\!v\big\}\,.\EE
The identification~\eref{RboundidentD_e} canonically lifts to an identification
\BE{RboundidentD_e2}\begin{split}
\cN_{\R\ov\cM_{0,\ell}}(\R D_{I;J,K}) \stackrel{\approx}{\lra} 
\big(\cN_{\nod;[\ell],I}\!\oplus\!\ov{\cN_{\nod;[\ell],I}}\big)^{\R}
&\equiv
\big\{\!(v,w)\!\in\!\cN_{\nod;[\ell],I}\!\oplus\!\ov{\cN_{\nod;[\ell],I}}\!:v\!=\!w\big\}\\
&\approx \cN_{\nod;[\ell],I}\,.
\end{split}\EE
The orientations on $\R\ov\cM_{0,\ell}$, $\R E_{J,K}$, and~$\R D_{I;J,K}$ specified above
determine orientations on the normal bundles on the left-hand sides 
in~\eref{RboundidentE_e2} and~\eref{RboundidentD_e2}.
The real line bundle on  the right-hand side in~\eref{RboundidentE_e2} is canonically oriented
by the ray $v\!\otimes_{\C}\!v$ with $v\!\in\!L_{\nod}$.
The complex line bundle on the right-hand side in~\eref{RboundidentD_e2} is canonically oriented
as well.

\begin{lmm}\label{DEorient_lmm}
Suppose $\ell\!\in\!\Z^+$ with $\ell\!\ge\!2$.
If $\{J,K\}\!\in\!\cP(\ell)$ with $J\!\in\!1$,
the isomorphism~\eref{RboundidentE_e2} is orientation-reversing.
If $(I,\{J,K\})\!\in\!\wt\cP_{\bu}(\ell)$ with $\min(J\!\cup\!K)\!\in\!J$,
the isomorphism~\eref{RboundidentD_e2} is orientation-preserving.
\end{lmm}

\begin{proof}
If $k\!\in\!K$, interchanging the two marked points in the $k$-th conjugate pair reverses
the orientations of $\R\ov\cM_{0,\ell}$, $\R E_{J,K}$, and~$\R D_{I;J,K}$ 
and thus preserves the orientations of the normal bundles.
We can thus assume that $K\!=\!\eset$.
Since the exact sequences~\eref{cMorientR_e}, \eref{ffEses_e}, and~\eref{ffDses_e} 
respect the relevant orientations,
we can also assume that $\ell\!=\!2$ in the first case and $\ell\!=\!3$ in the second case.\\

\noindent
In the first case, $\R E_{J,K}$ is then the one-point set consisting of the last curve
in Figure~\ref{M02_fig}.
The positive direction of the real line bundle on the right-hand side of~\eref{RboundidentE_e2} 
at this curve points into~$\cM_{0,2}^{\tau}$.
Thus, it is the opposite of the positive direction determined by the specified
orientation of~$\R\ov\cM_{0,2}$.\\

\noindent
If $l\!=\!3$ and $K\!=\!\eset$ in the second case, 
the interchange of two pairs of conjugate points preserves 
the orientations of $\R\ov\cM_{0,\ell}$ and~$\R D_{I;J,K}$ and thus of the normal bundles.
We can thus assume that $I\!=\!\{1\}$ and $J\!=\!\{2,3\}$.
The restriction of~$\ff_{3^{\pm}}^{\R}$ to~$\R D_{I;J,K}$ is then 
an orientation-preserving diffeomorphism
and $\ov\cM_{\{\nod\}\sqcup(J\cup K)}$ is a single-point.
The complex line bundle on the right-hand side of~\eref{RboundidentD_e2} 
in this case is isomorphic to
the restriction of the vertical tangent bundle in~\eref{cMorientR_e} to~$\R D_{I;J,K}$.
Thus, the former has the same orientation as the normal bundle on the left-hand side 
of~\eref{RboundidentD_e2}.  
\end{proof}

\noindent
For the purposes of orienting the boundary divisors 
\hbox{$\R D_{I';J',K'}\!\subset\!\R\ov\cM_{0,\{\nod\}\sqcup I}$} with \hbox{$I\!\subset\!\Z^+$},
we order the set $\{\nod\}\!\sqcup\!I$ by identifying~$\nod$ with $\min(\Z^+\!-\!I)$.
If \hbox{$(I,\{J,K\})\!\in\!\wt\cP_{\bu}(\ell)$} with \hbox{$\min(J\!\cup\!K)\!\in\!J$}  and 
\hbox{$(I',\{J',K'\})\!\in\!\wt\cP_{\bu}(\{\nod\}\!\sqcup\!I)$}
 with $\nod\!\in\!J'$, then
\BE{nodIord_e} \min(J'\!\cup\!K')\!\in\!J' ~~\Llra~~
\min\!\big(J\!\cup\!(J'\!-\!\{\nod\}\!)\!\cup\!(K\!\cup\!K')\!\big)
\in J\!\cup\!\big(J'\!-\!\{\nod\}\!\big)\,.\EE
The next two lemmas are used to establish the commutativity of the diagram~\eref{HRdiag_e}.

\begin{lmm}\label{RErestr_lmm}
Suppose $\ell\!\in\!\Z^+$ with $\ell\!\ge\!2$ and
\hbox{$\{J,K\}\!\in\!\cP(\ell)$} with $\min(J\!\cup\!K)\!\in\!J$.
Under the identification~\eref{RboundidentE_e}, 
$$D_{J',K'}=(-1)^{|K|}\ep_{J\cap K',K\cap K'}\R D_{J'-\{\nod\};J\cap K', K\cap K'}
\big|_{\R E_{J,K}}
\in H^2\big(\R E_{J,K};\Z\big)$$
for all \hbox{$\{J',K'\}\!\in\!\cP_{\bu}(\{\nod\}\!\sqcup\![\ell])$} with \hbox{$\nod\!\in\!J'$}.
\end{lmm}

\begin{proof} The submanifolds 
$\R E_{J,K},\R D_{J'-\{\nod\};J\cap K', K\cap K'}\!\subset\!\R\ov\cM_{0,\ell}$
are transverse and 
$$\R E_{J,K}\cap\R D_{J'-\{\nod\};J\cap K', K\cap K'}=D_{J',K'}$$
under the identification~\eref{RboundidentE_e}.
Thus, we only need to compare the orientations on the canonically identified bundles
\BE{RErestr_e5}\cN_{\ov\cM_{\{\nod\}\sqcup[\ell]}}D_{J',K'}\!=\!
\cN_{\R\ov\cM_{0,\ell}}\big(\R D_{J'-\{\nod\};J\cap K', K\cap K'}\big)
\big|_{D_{J',K'}}\lra D_{J',K'}\,.\EE
The orientation on the first bundle is given by its complex orientation 
if and only if $|K|\!\in\!2\Z$.
By the second statement of Lemma~\ref{DEorient_lmm}, 
the orientation on the second bundle is given by the orientation
of the restriction of the complex line bundle~$\cN_{\nod;[\ell],J'-\{\nod\}}$. 
This restriction is the complex line bundle $\cN_{\ov\cM_{\{\nod\}\sqcup[\ell]}}D_{J',K'}$
if 
$$\min(K')\!\equiv\!\min\big(\!(J\!\cap\!K')\!\cup\!(K\!\cap\!K')\!\big)\in J\!\cap\!K'$$
and the conjugate complex line bundle otherwise.
Thus, the orientations of the two bundles in~\eref{RErestr_e5} are the same if 
and only if \hbox{$(-1)^{|K|}\!=\!\ep_{J\cap K',K\cap K'}$}.
\end{proof}

\begin{lmm}\label{RDrestr_lmm}
Suppose $\ell\!\in\!\Z^+$ with $\ell\!\ge\!2$ and 
\hbox{$(I,\{J,K\})\!\in\!\wt\cP_{\bu}(\ell)$}  with $\min(J\!\cup\!K)\!\in\!J$.
Under the identification~\eref{RboundidentD_e}, 
$$1\!\otimes\!D_{J',K'}=(-1)^{|K|}\ep_{J\cap K',K\cap K'}
\R D_{I\cup(J'-\{\nod\});J\cap K', K\cap K'}
\big|_{\R D_{I;J,K}}\in H^2\big(\R D_{I;J,K};\Z\big)$$
for all \hbox{$\{J',K'\}\!\in\!\cP_{\bu}(\{\nod\}\!\sqcup\!(J\!\cup\!K))$} 
with \hbox{$\nod\!\in\!J'$}.
For all \hbox{$\{J',K'\}\!\in\!\cP_{\bu}(\{\nod\}\!\sqcup\!I)$} with \hbox{$\nod\!\in\!J'$},
$$\R E_{J',K'}\!\otimes\!1=\R E_{J\cup(J'-\{\nod\}),K\cup K'}
\big|_{\R D_{I;J,K}}\in H^1\big(\R D_{I;J,K};\Z\big)$$
under this identification.
If  \hbox{$(I',\{J',K'\})\!\in\!\wt\cP_{\bu}(\{\nod\}\!\sqcup\!I)$}, then
$$\R D_{I';J',K'}\!\otimes\!1=\begin{cases}
\R D_{(I'-\{\nod\})\cup J\cup K;J',K'}\big|_{\R D_{I;J,K}},&\hbox{if}~\nod\!\in\!I';\\
\R D_{I';J\cup (J'-\{\nod\}),K\cup K'}\big|_{\R D_{I;J,K}},&\hbox{if}~\nod\!\in\!J';
\end{cases}$$
under the identification~\eref{RboundidentD_e}.
\end{lmm}

\begin{proof} The submanifold $\R D_{I;J,K}\!\subset\!\R\ov\cM_{0,\ell}$
is transverse to each submanifold on the right-hand sides of the three equations.
In the first case,
$$\R D_{I;J,K}\!\cap\!\R D_{I\cup(J'-\{\nod\});J\cap K', K\cap K'}
=\R\ov\cM_{0,\{\nod\}\sqcup I}\!\times\!D_{J',K'}$$
under the identification~\eref{RboundidentD_e}.
The sign in the first equation is obtained as in the proof of Lemma~\ref{RErestr_lmm}.
In the second case,
$$\R D_{I;J,K}\!\cap\!\R E_{J\cup(J'-\{\nod\}),K\cup K'}
=\R E_{J',K'}\!\times\!\ov\cM_{\{\nod\}\sqcup(J\cup K)}$$
under the identification~\eref{RboundidentD_e}.
By the first statement of Lemma~\ref{DEorient_lmm}, 
the two relevant oriented normal bundles are now the same.
In the first case of the last equation,
$$\R D_{I;J,K}\!\cap\!\R D_{(I'-\{\nod\})\cup J\cup K;J',K'}
=\R D_{I';J',K'}\!\times\!\ov\cM_{\{\nod\}\sqcup(J\cup K)}\,.$$
The relevant oriented normal bundles are now the same 
by the second statement of Lemma~\ref{DEorient_lmm}.
In the second case of the last equation,
$$\R D_{I;J,K}\!\cap\!\R D_{I';J\cup (J'-\{\nod\}),K\cup K'}
=\R D_{I';J',K'}\!\times\!\ov\cM_{\{\nod\}\sqcup(J\cup K)}\,.$$
The relevant oriented normal bundles in this case are the same 
by the second statement of Lemma~\ref{DEorient_lmm} and~\eref{nodIord_e}.
\end{proof}

\subsection{The real case: topology}
\label{Rhomol_subs2}

\noindent 
For $i\!\in\![\ell]$ and $\bu\!=\!\pm$, let
\BE{wtDi0dfn_e}
\wt{D}_{[\ell^{\pm}]-\{i^{\bu}\}}^0\equiv  
\begin{cases}
\R D_{[\ell]-\{i\};\{i,\ell+1\},\eset},&\hbox{if}~\bu\!=\!+;\\
\R D_{[\ell]-\{i\};\{i\},\{\ell+1\}},&\hbox{if}~\bu\!=\!-;
\end{cases}\EE
be as in~\eref{wtDidfn_e}.
For $\ka\!\in\!H_p(\ov\cM_{\ell};R)$, define
\begin{equation*}\begin{split}
\R\ov\cM_{0,\ell+1}\big|_{\ka}
&=\PD_{\R\ov\cM_{0,\ell+1}}\!\big(\!
\big\{\ff_{(\ell+1)^{\pm}}^{\R}\big\}^{\!*}
\big(\PD_{\R\ov\cM_{0,\ell}}\!(\ka)\!\big)\!\big)
\in H_{p+2}\big(\R\ov\cM_{0,\ell+1};R\big),\\
\wt{D}_{[\ell^{\pm}]-\{1^+\}}^0\big|_{\ka}
&=\wt{D}_{[\ell^{\pm}]-\{1^+\}}^0\!\cap\!\big(\R\ov\cM_{0,\ell+1}\big|_{\ka}\big)
\in H_p\big(\R\ov\cM_{0,\ell+1};R\big).
\end{split}\end{equation*}
Suppose $\{J,K\}\!\in\!\cP(\ell)$. 
If $J\!\ni\!1$, define
\BE{wtREdfn_e}
\R\wt{E}_{J,K}^+=\R E_{J\cup\{\ell+1\},K}, \quad
\R\wt{E}_{J,K}^0=\R D_{\{\ell+1\};J,K}, \quad
\R\wt{E}_{J,K}^-=\R E_{J,K\cup\{\ell+1\}}\,.\EE
If $\ka_{J,K}\!\in\!H_{p-2+\de_{0\bu}}(\R E_{J,K};R)$ and $\bu\!=\!+,0,-$
with $\de_{0\bu}\!=\!1$ if $\bu\!=\!0$ and  $\de_{0\bu}\!=\!0$ otherwise, let
$$\wt\ka_{J,K}^{\bu}=\big\{\wt\io_{J,K}^{\bu}\big\}_{\!*}
\Big(\PD_{\R\wt{E}_{J,K}^{\bu}}\!\big(\big\{\ff_{(\ell+1)^{\pm}}^{\R}
\big|_{\R\wt{E}_{J,K}^{\bu}}\big\}^{\!*}(\PD_{\R E_{J,K}}\!(\ka_{J,K})\!\big)\!\big)
\!\Big) \in H_p\big(\R\ov\cM_{0,\ell+1};R\big),$$
where $\wt\io_{J,K}^{\bu}\!:\R\wt{E}_{J,K}^{\bu}\!\lra\!\R\ov\cM_{0,\ell+1}$ is the inclusion map.\\

\noindent
If $(I,\{J,K\})\!\in\!\wt\cP_{\bu}(\ell)$ with $\min(J\!\cup\!K)\!\in\!J$, define
\BE{wtRDdfn_e}\begin{split}
\big(\R\wt{D}_{I;J,K}^+,\R\wt{D}_{I;J,K}^0\big)&=
\begin{cases}
\big(\R D_{I\cup\{\ell+1\};J,K},\R D_{I;J\cup\{\ell+1\},K}\big),
&\hbox{if}~1\!\in\!I;\\
\big(\R D_{I;J\cup\{\ell+1\},K},\R D_{I\cup\{\ell+1\};J,K}\big),
&\hbox{if}~1\!\in\!J;
\end{cases}\\
\R\wt{D}_{I;J,K}^-&=\R D_{I;J,K\cup\{\ell+1\}}\,.
\end{split}\EE
If $\ka_{I;J,K}\!\in\!H_{p-2}(\R D_{I;J,K};R)$ and $\bu\!=\!+,0,-$, let
$$\wt\ka_{I;J,K}^{\bu}=\big\{\wt\io_{I;J,K}^{\bu}\big\}_{\!*}
\Big(\PD_{\R\wt{D}_{I;J,K}^{\bu}}\!\big(\big\{\ff_{(\ell+1)^{\pm}}^{\R}
\big|_{\R\wt{D}_{I;J,K}^{\bu}}\big\}^{\!*}(\PD_{\R D_{I;J,K}}\!(\ka_{I;J,K})\!\big)\!\big)
\!\Big) \in H_p\big(\R\ov\cM_{0,\ell+1};R\big),$$
where $\wt\io_{I;J,K}^{\bu}\!:\R\wt{D}_{I;J,K}^{\bu}\!\lra\!\R\ov\cM_{0,\ell+1}$ is the inclusion map.\\

\noindent
For the remainder of this section and in~\eref{HRdiag_e},
we take the homology and cohomology groups in a commutative ring~$R$ with unity,
which we drop from the notation for space reasons.

\begin{prp}\label{Rhomol_prp}
For every $\ell\!\in\!\Z^+$ with $\ell\!\ge\!2$ and $p\!\in\!\Z$, the homomorphism
\BE{Rhomol_e}\begin{split}
&\Psi\!:H_p\big(\R\ov\cM_{0,\ell}\big)\!\oplus\!H_{p-2}\big(\R\ov\cM_{0,\ell}\big)
\!\oplus\!\bigoplus_{(I,\{J,K\})\in\wt\cP_{\bu}(\ell)}\hspace{-.34in}
H_{p-2}\big(\R D_{I;J,K}\big)^{\!\oplus2}\\
&\hspace{1in}
\oplus\!\bigoplus_{\{J,K\}\in\cP(\ell)}\hspace{-.24in}
\big(H_{p-1}\big(\R E_{J,K}\big)\!\oplus\!H_{p-2}\big(\R E_{J,K}\big)\!\big)
\lra H_p\big(\R\ov\cM_{0,\ell+1}\big),\\
&
\Psi\big(\ka_0,\ka,(\ka_{I;J,K},\ka_{I;J,K}')_{(I,\{J,K\})\in\wt\cP_{\bu}(\ell)},
(\ka_{J,K},\ka_{J,K}')_{\{J,K\}\in\cP(\ell)}\big)\\
&\hspace{.2in}
=\wt{D}_{[\ell^{\pm}]-\{1^+\}}^0\big|_{\ka_0}\!+\!\R\ov\cM_{0,\ell+1}\big|_{\ka}
\!+\!\sum_{(I,\{J,K\})\in\wt\cP_{\bu}(\ell)}\hspace{-.34in}
\big(\wt\ka_{I;J,K}^0\!+\!\wt\ka_{I;J,K}'^-\big)
\!+\!\sum_{\{J,K\}\in\cP_{\bu}(\ell)}\hspace{-.24in}
\big(\wt\ka_{J,K}^0\!+\!\wt\ka_{J,K}'^-\big)\,,
\end{split}\EE
of $R$-modules is injective.
If 2 is a unit in the coefficient ring~$R$,
then this homomorphism is also surjective.
\end{prp}

\begin{proof} By Proposition~\ref{Ralg_prp}, 
the top homomorphism in the diagram~\eref{HRdiag_e} is surjective.
By the generation property~\ref{HRgen_it} on page~\pageref{HRgen_it} with~$\ell$ replaced 
by~$\ell\!+\!1$ and its proof, 
so is the right homomorphism if 2 is a unit in the coefficient ring~$R$.
By Theorem~\ref{HC_thm} with~$\ell$ replaced by each $\ell'\!\le\!\ell\!+\!1$
and Theorem~\ref{HC_thm} with~$\ell$ replaced by each $\ell'\!\le\!\ell$,
the left homomorphism is an isomorphism under this assumption on~$R$.
Along with the commutativity of the diagram~\eref{HRdiag_e},
this implies the surjectivity claim of the proposition.\\

\noindent
By the same reasoning as in the proof of Proposition~\ref{Chomol_prp},
the restrictions of~$\Psi$ to~$H_p(\R\ov\cM_{0,\ell})$ and~$H_{p-2}(\R\ov\cM_{0,\ell})$
are injective, $\Psi(H_p(\R\ov\cM_{0,\ell})\!)$ is disjoint from the image of~$\Psi$
on the rest of the direct sum in the domain of~$\Psi$,
and $\Psi(H_{p-2}(\R\ov\cM_{0,\ell})\!)$ is disjoint from the image of~$\Psi$
on the direct sum of the two big direct sums.
By Lemma~\ref{Rhomol_lmm1} and the Poincare Duality for~$\R E_{J,K}$, 
the restriction of~$\Psi$ to the direct sum of the $H_{p-1}\big(\R E_{J,K})$ terms 
in its domain is injective and
$$\Psi\bigg(\!\bigoplus_{\{J,K\}\in\cP(\ell)}\hspace{-.22in}H_{p-1}\big(\R E_{J,K}\big)\!\!\bigg)
\cap\Psi\bigg(\!\bigoplus_{(I,\{J,K\})\in\wt\cP_{\bu}(\ell)}\hspace{-.34in}
H_{p-2}\big(\R D_{I;J,K}\big)^{\!\oplus2}
\oplus\!\bigoplus_{\{J,K\}\in\cP(\ell)}\hspace{-.22in}H_{p-2}\big(\R E_{J,K}\big)\!\!\bigg)
=\{0\}\,.$$
By Lemmas~\ref{Rhomol_lmm1} and~\ref{Rhomol_lmm2} and the Poincare Duality for~$\R E_{J,K}$ again,
the restriction of~$\Psi$ to the direct sum of the $H_{p-2}\big(\R E_{J,K})$ terms 
in its domain is injective and
$$\Psi\bigg(\!\bigoplus_{\{J,K\}\in\cP(\ell)}\hspace{-.22in}H_{p-2}\big(\R E_{J,K}\big)\!\!\bigg)
\cap\Psi\bigg(\!\bigoplus_{(I,\{J,K\})\in\wt\cP_{\bu}(\ell)}\hspace{-.32in}
H_{p-2}\big(\R D_{I;J,K}\big)^{\!\oplus2}\!\bigg)=\{0\}\,.$$
By Lemma~\ref{Rhomol_lmm3} and the Poincare Duality for~$\R D_{I;J,K}$,
the restriction of~$\Psi$ to the direct sum over~$\wt\cP_{\bu}(\ell)$ is injective 
as well.
Thus, $\Psi$ is injective.
\end{proof}

\noindent
We denote by $\cdot_{\R E_{J,K}}$,  $\cdot_{\R D_{I;J,K}}$, and 
$\cdot_{\R\ov\cM_{0,\ell+1}}$ the homology intersection pairings
on~$\R E_{J,K}$, $\R D_{I;J,K}$, and $\R\ov\cM_{0,\ell+1}$, respectively.
The next three lemmas are analogues of Lemma~\ref{Chomol_lmm} in the present setting.

\begin{lmm}\label{Rhomol_lmm1}
Suppose $\ell\!\in\!\Z^+$ with $\ell\!\ge\!2$, $p\!\in\!\Z$, 
$\{J',K'\}\!\in\!\cP(\ell)$, 
\hbox{$\ka'\!\in\!H_{2\ell-3-p}(\R E_{J',K'})$},
and $\bu\!=\!0,-$.
If \hbox{$\{J,K\}\!\in\!\cP(\ell)$}
and \hbox{$\ka\!\in\!H_{p-2+\de_{0\bu}}(\R E_{J,K})$}, then
\BE{Rhomol1_e1}\wt\ka^{\bu}\!\cdot_{\R\ov\cM_{0,\ell+1}}\!\wt\ka'^+
=\pm\begin{cases}\ka\!\cdot_{\R E_{J,K}}\!\ka',
&\hbox{if}~\{J,K\}\!=\!\{J',K'\},\,\bu\!=\!0;\\
0,&\hbox{otherwise}.
\end{cases}\EE
If \hbox{$(I,\{J,K\})\!\in\!\wt\cP_{\bu}(\ell)$} 
and \hbox{$\ka\!\in\!H_{p-2}(\R D_{I;J,K})$}, then
\BE{Rhomol1_e2}\wt\ka^{\bu}\!\cdot_{\R\ov\cM_{0,\ell+1}}\!\wt\ka'^+=0.\EE
\end{lmm}

\begin{proof} The second case of~\eref{Rhomol1_e2} follows from
$$\R\wt{E}_{J,K}^{\bu}\!\cap\!\R\wt{E}_{J',K'}^+=\eset
\qquad\hbox{if}\quad \{J,K\}\!\neq\!\{J',K'\}
~\hbox{or}~\bu\!=\!-\,.$$
The proof of~\eref{Rhomol1_e2} is similar to the proof of
the $(J,K)\!\neq\!(J',K')$ case of Lemma~\ref{Chomol_lmm}.
The first case of~\eref{Rhomol1_e1} is obtained similarly to 
the first case of Lemma~\ref{Chomol_lmm}.
The signs of the intersection points~$(z_r,z_r')$ and 
$(\!(z_r,\wt\cC_r),(z_r',\wt\cC_r)\!)$ contributing to  
$\ka\!\cdot_{\R E_{J,K}}\!\ka'$ and
$\wt\ka^0\!\cdot_{\R\ov\cM_{0,\ell+1}}\!\wt\ka'^+$
may now be different,
but whether they are the same or opposite is the same for all such pairs
(this is determined by the exact definition of the intersection products and
the parity of~$p$).
\end{proof}

\begin{lmm}\label{Rhomol_lmm2}
Let $\ell$, $p$, $\{J',K'\}$, and $\ka'$ be as in Lemma~\ref{Rhomol_lmm1}.
If \hbox{$\{J,K\}\!\in\!\cP(\ell)$}
and \hbox{$\ka\!\in\!H_{p-2+\de_{0\bu}}(\R E_{J,K})$}, then
$$\wt\ka^-\!\cdot_{\R\ov\cM_{0,\ell+1}}\!\wt\ka'^0
=\pm\begin{cases}\ka\!\cdot_{\R E_{J,K}}\!\ka',
&\hbox{if}~\{J,K\}\!=\!\{J',K'\};\\
0,&\hbox{otherwise}.\end{cases}$$
If \hbox{$(I,\{J,K\})\!\in\!\wt\cP_{\bu}(\ell)$}, 
\hbox{$\ka\!\in\!H_{p-2}(\R D_{I;J,K})$},
and \hbox{$\bu\!=\!0,-$}, then
$$\wt\ka^{\bu}\!\cdot_{\R\ov\cM_{0,\ell+1}}\!\wt\ka'^0=0.$$
\end{lmm}

\begin{proof}
The proof is similar to Lemma~\ref{Rhomol_lmm2}.
\end{proof}

\begin{lmm}\label{Rhomol_lmm3}
Suppose $\ell\!\in\!\Z^+$ with $\ell\!\ge\!2$, $p\!\in\!\Z$, 
\hbox{$(I,\{J,K\}),(I',\{J',K'\})\!\in\!\wt\cP_{\bu}(\ell)$}, 
\hbox{$\ka\!\in\!H_{p-2}(\R D_{I;J,K})$}, and 
\hbox{$\ka'\!\in\!H_{2\ell-2}(\R D_{I';J',K'})$}.
If $1\!\not\in\!I$ and $\bu,\circ\!=\!0,-,+$ are distinct, then
\BE{Rhomol3_e1}\wt\ka^{\bu}\!\cdot_{\R\ov\cM_{0,\ell+1}}\!\wt\ka'^{\circ}
=\pm\begin{cases}\ka\!\cdot_{\R D_{I;J,K}}\!\ka',
&\hbox{if}~(I,\{J,K\})\!=\!(I',\{J',K'\}),\,
\{\bu,\circ\}\!=\!\{+,0\},\{0,-\};\\
0,&\hbox{otherwise}.\end{cases}\EE
If $1\!\in\!I$ and $\bu,\circ\!=\!0,-$, then
\BE{Rhomol3_e2}\wt\ka^{\bu}\!\cdot_{\R\ov\cM_{0,\ell+1}}\!\wt\ka'^{\circ}
=\pm\begin{cases}\ka\!\cdot_{\R D_{I;J,K}}\!\ka',
&\hbox{if}~(I,\{J,K\})\!=\!(I',\{J',K'\}),\,
\bu\!=\!\circ;\\
0,&\hbox{otherwise}.\end{cases}\EE
\end{lmm}

\begin{proof} The proof of the $(I,\{J,K\})\!\neq\!(I',\{J',K'\})$ cases 
of~\eref{Rhomol3_e1} and~\eref{Rhomol3_e2}  is as in the second
case of Lemma~\ref{Chomol_lmm}.
Furthermore,
$$\R\wt{D}_{I;J,K}^+\!\cap\!\R\wt{D}_{I;J,K}^-=\eset~~\hbox{if}~1\!\not\in\!I,
\qquad
\R\wt{D}_{I;J,K}^0\!\cap\!\R\wt{D}_{I;J,K}^-=\eset~~\hbox{if}~1\!\in\!I.$$
This establishes the claims in the second cases of~\eref{Rhomol3_e1} and~\eref{Rhomol3_e2}.
The first case of~\eref{Rhomol3_e1} is obtained similarly to
the first cases of Lemma~\ref{Chomol_lmm} and~\eref{Rhomol1_e1}.\\

\noindent
Suppose now that $(I,\{J,K\})\!=\!(I',\{J',K'\})$ and $1\!\in\!J$.
We can choose pseudocycle representatives
$$f\!:Z\lra \R D_{I;J,K} \qquad\hbox{and}\qquad f'\!:Z'\lra \R D_{I;J,K}$$
for~$\ka$ and~$\ka'$, respectively, 
so that they are transverse in~$\R D_{I;J,K}$ and 
thus intersect at a finite subset \hbox{$Z_{\cap}\!\subset\!Z\!\times\!Z'$} 
of signed points.
If the two pseudocycles are chosen generically, we can identity~$Z_{\cap}$
with a subset of~$\R D_{I;J,K}$ via~$f$. 
Via the second identification in~\eref{Rboundident_e},
we write the elements of~$Z_{\cap}$ as~pairs
$$(\cC_1,\cC_2)\in  \R\ov\cM_{0,\{\nod\}\sqcup I}\!\times\!\ov\cM_{\{\nod\}\sqcup(J\cup K)\cup\{\ell+1\}}\,.$$
The pseudocycle representatives
\begin{equation*}\begin{split}
\R\wt{D}_{I;J,K}^{\bu}\big|_f\!:f^*\R\wt{D}_{I;J,K}^{\bu}&\lra 
\R\wt{D}_{I;J,K}^{\bu}\subset \R\ov\cM_{0,\ell+1} 
\qquad\hbox{and}\\
\R\wt{D}_{I;J,K}^{\bu}\big|_{f'}\!:f'^*\R\wt{D}_{I;J,K}^{\bu}&\lra \R\wt{D}_{I;J,K}^{\bu}
\subset  \R\ov\cM_{0,\ell+1}
\end{split}\end{equation*}
for~$\wt\ka^{\bu}$ and~$\wt\ka'^{\bu}$, respectively, then intersect transversely 
in~$\R\wt{D}_{I;J,K}^{\bu}$ along~$\R\wt{D}_{I;J,K}^{\bu}\big|_{Z_{\cap}}$.
Thus,
\BE{Rhomol3_e3}\wt\ka^{\bu}\!\cdot_{\R\ov\cM_{0,\ell+1}}\!\wt\ka'^{\bu}
=\blr{e\big(\cN_{\R\ov\cM_{0,\ell+1}}(\R\wt{D}_{I;J,K}^{\bu})\!\big),
\R\wt{D}_{I;J,K}^{\bu}\big|_{Z_{\cap}}}\,.\EE

\vspace{.15in}

\noindent
Under the analogue of the second identification in~\eref{Rboundident_e},
$$\R\wt{D}_{I;J,K}^{\bu}\approx
\R\ov\cM_{0,\{\nod\}\sqcup I}\!\times\!\ov\cM_{\{\nod\}\sqcup(J\cup K\cup\{\ell+1\})}
\quad\hbox{and}\quad
\R\wt{D}_{I;J,K}^{\bu}\big|_{Z_{\cap}}\approx 
\bigsqcup_{(\cC_1,\cC_2)\in Z_{\cap}}\hspace{-.2in}
\{\cC_1\}\!\times\!\ff_{\ell+1}^{\,-1}(\cC_2),$$
where 
$$\ff_{\ell+1}\!: \ov\cM_{\{\nod\}\sqcup(J\cup K\cup\{\ell+1\})}\lra
\ov\cM_{\{\nod\}\sqcup(J\cup K)}$$
is the forgetful morphism dropping the last marked point.
With the notation as above Lemma~\ref{DEorient_lmm}, 
$$c_1(L_{\nod})=\ff_{\ell+1}^{\,*}c_1(L_{\nod})\!-\!D_{J\cup K,\{\nod,\ell+1\}}\in 
H^2\big(\ov\cM_{\{\nod\}\sqcup(J\cup K\cup\{\ell+1\})}\big).$$
Along with~\eref{Rhomol3_e3} and the second statement of Lemma~\ref{DEorient_lmm}, this gives
\begin{equation*}\begin{split}
\wt\ka^{\bu}\!\cdot_{\R\ov\cM_{0,\ell+1}}\!\wt\ka'^{\bu}
&=\sum_{(\cC_1,\cC_2)\in Z_{\cap}}\hspace{-.2in}
\blr{c_1\big( \cN_{\nod;[\ell],I}\!\big),\{\cC_1\}\!\times\!\ff_{\ell+1}^{\,-1}(\cC_2)\!}\\
&=-\sum_{(\cC_1,\cC_2)\in Z_{\cap}}\hspace{-.2in}
\blr{1\!\times\!D_{J\cup K,\{\nod,\ell+1\}},\{\cC_1\}\!\times\!\ff_{\ell+1}^{\,-1}(\cC_2)\!}
=\pm\big(\ka\!\cdot_{\R D_{I;J,K}}\!\ka'\big);
\end{split}\end{equation*}
whether the sign of the $(\cC_1,\cC_2)$ summand above agrees with the sign of
$(\cC_1,\cC_2)$ as an element of~$Z_{\cap}$ depends only on the parity of~$K$ and
on the choice of~$\bu$.
This establishes the claim in the first case of~\eref{Rhomol3_e2}.
\end{proof}

\subsection{The real case: algebra}
\label{Ralg_subs}

\noindent
For a finite set~$S$, let 
$$\wt\cP(S)=\big\{\!(I,\{J,K\})\!\!:S\!=\!I\!\sqcup\!J\!\sqcup\!K\big\}, \quad
\wt\cP_{\bu}(S)=\big\{\!(I,\{J,K\})\!\in\!\wt\cP(S)\!\!:
\,1\!\le\!|I|\!\le\!|S|\!-\!2\big\}.$$
For a commutative ring~$R$ with unity, denote~by
$$\cI_S\subset\!R_S\!\equiv\!R\big[(D_{J,K})_{\{J,K\}\in\cP_{\bu}(S)}\big]$$
the ideal generated by the left-hand sides of~\eref{cM04rel_e1a} and~\eref{cM04rel_e2}
with~$\cP_{\bu}(\ell)$ replaced by~$\cP_{\bu}(S)$ and $a,b,c,d\!\in\!S$ distinct.
If in addition $S$ is ordered, denote~by
$$\cI_{0,S}\subset R_{0,S}\!\equiv\!R\big[(\R E_{J,K})_{\{J,K\}\in\cP(S)},
(\R D_{I;J,K})_{(I,\{J,K\})\in\wt\cP_{\bu}(S)}\big]$$
the ideal generated by the left-hand sides of~\eref{RcM04rel_e1a}-\eref{RcM04rel_e2b}
with $\cP(\ell),\wt\cP_{\bu}(\ell)$ replaced by $\cP(S),\wt\cP_{\bu}(S)$
and $a,b,c\!\in\!S$ distinct.
Let
$$\cR_S=R_S/\cI_S \quad\hbox{and}\quad
\cR_{0,S}=R_{0,S}/\cI_{0,S}\,.$$
Thus, the left-hand sides of~\eref{HCthm_e} and~\eref{HRthm_e} 
are the rings $\cR_{\ell}\!\equiv\!\cR_{[\ell]}$ with $R\!=\!\Z$
and $\cR_{0,\ell}\!\equiv\!\cR_{0,[\ell]}$ with $R\!=\!\Q$.
Let $\cI_{0,\ell}\!=\!\cI_{0,[\ell]}$ and $R_{0,\ell}\!=\!R_{0,[\ell]}$.\\

\noindent
Let $\ell\!\in\!\Z^+$ with $\ell\!\ge\!2$.
Define an $R$-algebra homomorphism
\BE{FdfnR_e}\begin{aligned}
F\!:R_{0,\ell}\lra R_{0,\ell+1}, \quad
F\big(\R E_{J,K}\big)&=\R\wt{E}_{J,K}^+\!+\!\R\wt{E}_{J,K}^-
&&\forall\,\{J,K\}\!\in\!\cP(\ell),\\
F\big(\R D_{I;J,K}\big)&=\R\wt{D}_{I;J,K}^+\!+\!\R\wt{D}_{I;J,K}^0\!+\!\R\wt{D}_{I;J,K}^-
&&\forall\,(I,\{J,K\})\!\in\!\wt\cP_{\bu}(\ell),
\end{aligned}\EE
with  $\R\wt{E}_{J,K}^{\bu},\R\wt{D}_{I;J,K}^{\bu}\!\in\!R_{0,\ell+1}$, $\bu\!=\!+,0,-$,
defined as in~\eref{wtREdfn_e} and~\eref{wtRDdfn_e}.
For \hbox{$\{J,K\}\!\in\!\cP(\ell)$} with \hbox{$\min(J\!\cup\!K)\!\in\!J$}, 
define an $\R$-algebra homomorphism
\begin{gather*}
F_{J,K}\!:R_{\{\nod\}\sqcup[\ell]}\lra R_{0,\ell+1} \qquad\hbox{by}\\
F_{J,K}\big(D_{J',K'}\big)=(-1)^{|K|}\ep_{J\cap K',K\cap K'}\R D_{J'-\{\nod\};J\cap K', K\cap K'} 
\quad\forall~\{J',K'\}\!\in\!\cP_{\bu}\!\big(\{\nod\}\!\sqcup\![\ell]\big),\,J'\!\ni\!\nod.
\end{gather*}
For $(I,\{J,K\})\!\in\!\wt\cP_{\bu}(\ell)$  with $\min(J\!\cup\!K)\!\in\!J$, 
define $\R$-algebra homomorphisms
\begin{gather*}
F_{I;J,K}^2\!:R_{\{\nod\}\sqcup(J\cup K)}\lra R_{0,\ell+1}
\quad\hbox{and}\quad
F_{J,K}\!:R_{0,\{\nod\}\sqcup I}\lra R_{0,\ell+1} 
\qquad\hbox{by}\\
F_{I;J,K}^2\big(D_{J',K'}\big)=(-1)^{|K|}\ep_{J\cap K',K\cap K'}\R D_{I\cup(J'-\{\nod\});J\cap K', K\cap K'} 
\quad\forall~\{J',K'\}\!\in\!\cP_{\bu}\!\big(\{\nod\}\!\sqcup\!(J\!\cup\!K)\!\big),\,J'\!\ni\!\nod,\\
F_{J,K}\big(\R E_{J',K'}\big)=\R E_{J\cup(J'-\{\nod\}),K\cup K'} 
\quad\forall~\{J',K'\}\!\in\!\cP\big(\{\nod\}\!\sqcup\!I\big),\,J'\!\ni\!\nod,\\
F_{J,K}\big(\R D_{I';J',K'}\big)=\begin{cases}
\R D_{(I'-\{\nod\})\cup J\cup K;J',K'},
&\hbox{if}~(I',\{J',K'\})\!\in\!\wt\cP_{\bu}(\{\nod\}\!\sqcup\!I),\,I'\!\ni\!\nod;\\
\R D_{I';J\cup (J'-\{\nod\}),K\cup K'},
&\hbox{if}~(I',\{J',K'\})\!\in\!\wt\cP_{\bu}(\{\nod\}\!\sqcup\!I),\,J'\!\ni\!\nod.
\end{cases}
\end{gather*}
We denote by 
$$F_{I;J,K}\!:R_{0,\{\nod\}\sqcup I}\!\otimes_R\!R_{\{\nod\}\sqcup(J\cup K)}\lra R_{0,\ell+1}$$
the $R$-algebra homomorphism induced by~$F_{J,K}$ and~$F_{I;J,K}^2$.
Let $\wt{D}_{[\ell^{\pm}]-\{1^+\}}^0\!\in\!R_{0,\ell+1}$ be defined as in~\eref{wtDi0dfn_e}.

\begin{prp}\label{Ralg_prp}
For every $\ell\!\in\!\Z^+$ with $\ell\!\ge\!2$, the homomorphism
\BE{Ralg_e}\begin{split}
&\Phi\!:\cR_{0,\ell}^{\,\oplus2}\!\oplus\!
\bigoplus_{(I,\{J,K\})\in\wt\cP_{\bu}(\ell)}\hspace{-.34in}
\big(\cR_{0,\{\nod\}\sqcup I}\!\otimes_R\!\cR_{\{\nod\}\sqcup(J\cup K)}\big)^{\!\oplus2}
\!\oplus\!\bigoplus_{\{J,K\}\in\cP(\ell)}\hspace{-.22in}
\cR_{\{\nod\}\sqcup(J\cup K)}^{\,\oplus2} \lra \cR_{0,\ell+1},\\
&\Phi\big([\ka_0],[\ka],([\ka_{I;J,K}],[\ka_{I;J,K}'])_{(I,\{J,K\})\in\wt\cP_{\bu}(\ell)},
([\ka_{J,K}],[\ka_{J,K}'])_{\{J,K\}\in\cP(\ell)}\big)\\
&\quad=\wt{D}_{[\ell^{\pm}]-\{1^+\}}^0F(\ka_0)\!+\!F(\ka)
+\sum_{\{J,K\}\in\cP_{\bu}(\ell)}\hspace{-.22in}
\Big(\R\wt{E}_{J,K}^0F\big(F_{J,K}(\ka_{J,K})\!\big)\!+\!\R\wt{E}_{J,K}^-F\big(F_{J,K}(\ka_{J,K}')\!\big)
\!\!\Big)\\
&\hspace{.8in}+\sum_{(I,\{J,K\})\in\wt\cP_{\bu}(\ell)}\hspace{-.34in}
\Big(\R\wt{D}_{I;J,K}^0F\big(F_{I;J,K}(\ka_{I;J,K})\!\big)\!+
\!\R\wt{D}_{I;J,K}^-F\big(F_{I;J,K}(\ka_{I;J,K}')\!\big)\!\!\Big)
+\cI_{0,\ell+1}\,,
\end{split}\EE
of $R$-modules is well-defined and surjective.
\end{prp}

\noindent
With $\Psi$ denoting the composition of the homomorphism~\eref{Rhomol_e}
with the Poincare duals on both sides, the diagram
\BE{HRdiag_e}\begin{split}
\xymatrix{\tn{domain of~$\Phi$ in~\eref{Ralg_e}}
\ar[rr]^>>>>>>>>>>>>>>>{\Phi}\ar[d]_{\eref{HRthm_e}}^{\eref{HCthm_e}}
&&\cR_{0,\ell+1}\ar[d]^{\eref{HRthm_e}} \\
\tn{domain of~$\Psi$ in~\eref{Rhomol_e}}\ar[rr]^>>>>>>>>>{\Psi}
&& H^*\big(\R\ov\cM_{0,\ell+1})}
\end{split}\EE
commutes.
The commutativity on the $\{J,K\}$ and $(I,\{J,K\})$ summands in the domain of~$\Phi$ follows from 
Lemmas~\ref{RErestr_lmm} and~\ref{RDrestr_lmm}, respectively.
By Proposition~\ref{Ralg_prp} and the generation property~\ref{HRgen_it} on 
page~\pageref{HRgen_it} with~$\ell$ replaced by~$\ell\!+\!1$, 
the top and right homomorphisms in this diagram are surjective.
By Theorem~\ref{HC_thm} with~$\ell$ replaced by each $\ell'\!\le\!\ell\!+\!1$,
Theorem~\ref{HR_thm} with~$\ell$ replaced by each $\ell'\!\le\!\ell$,
and Proposition~\ref{Chomol_prp}, 
the left and bottom homomorphisms are isomorphisms.
Along with the commutativity of the diagram~\eref{HRdiag_e},
this implies that the top and right homomorphisms are isomorphisms as well
and establishes the injectivity property~\ref{HRrel_it} on 
page~\pageref{HRrel_it} with~$\ell$ replaced by~$\ell\!+\!1$.

\begin{proof}[{\bf{\emph{Proof of Proposition~\ref{Ralg_prp}: well-defined}}}]
The $R$-algebra homomorphism~$F$ sends the left-hand side of~\eref{RcM04rel_e1a},
(resp.~\eref{RcM04rel_e1b}, \eref{RcM04rel_e1c}) to the sum of 4
(resp.~6,9) analogous terms in its target.
It sends  the left-hand sides of~\eref{RcM04rel_e2a} and~\eref{RcM04rel_e2b}
to their analogues in the target.
Thus, \hbox{$F(\cI_{0,\ell})\!\subset\!\cI_{0,\ell+1}$} and
$\Phi$ is well-defined on the $\cR_{0,\ell}^{\,\oplus2}$ summand in~\eref{Ralg_e}.
Below we denote by~$\eref{cM04rel_e2}_{abcd}^S$ 
the left-hand side of~\eref{cM04rel_e2} with~$[\ell]$ replaced by a finite set~$S$,
by~$\eref{RcM04rel_e2b}_{abc}^S$ 
the left-hand side of~\eref{RcM04rel_e2b} with~$[\ell]$ replaced by a finite ordered set~$S$,
and by~$\eref{RcM04rel_e2b2}_{abc}^S$ 
the left-hand side of~\eref{RcM04rel_e2b2} with~$[\ell]$ replaced by a finite ordered set~$S$.
For $s_0\!\in\!S$ and \hbox{$a,b,c\!\in\!S\!-\!\{s_0\}$} distinct,
\BE{RcM04rel_e2bgen}  \eref{RcM04rel_e2b}_{abc}^S=
\eref{RcM04rel_e2b}_{s_0bc}^S+\eref{RcM04rel_e2b}_{bas_0}^S
- \eref{RcM04rel_e2b}_{cas_0}^S\,.\EE
In particular, the collection of 
the left-hand sides of~\eref{RcM04rel_e2b} with~$[\ell]$ replaced by a finite ordered set~$S$
is linearly generated by its subcollection with $s_0\!\in\!\{a,b,b\}$.\\

\noindent
Let $\{J,K\}\!\in\!\cP(\ell)$ with $\min(J\!\cup\!K)\!\in\!J$. 
Since~$F_{J,K}$ maps the left-hand side of~\eref{cM04rel_e1a}
in its domain to a product as in~\eref{RcM04rel_e1c},
$F\!\circ\!F_{J,K}$ sends this product to~$\cI_{0,\ell+1}$.
Let $b,c,d\!\in\![\ell]$ be distinct.
If $d\!\in\!J$,
$$(-1)^{|K|}F_{J,K}\Big(\eref{cM04rel_e2}_{(\nod)bcd}^{\{\nod\}\sqcup[\ell]}\Big)
-\eref{RcM04rel_e2b2}_{dcb}^{[\ell]}
=-\!\!\!\sum_{\begin{subarray}{c}(I',\{J',K'\})\in\wt\cP_{\bu}(\ell)\\
d\in J',\,c\not\in I',\,b\in I'\\ \{J',K'\}\not\preceq\{J,K\} \end{subarray}}
\hspace{-.39in}\ep_{J',K'}\R D_{I';J',K'}
+\!\!\! \sum_{\begin{subarray}{c}(I',\{J',K'\})\in\wt\cP_{\bu}(\ell) \\
d\in J',\,c\in I',\,b\not\in I'\\ 
\{J',K'\}\not\preceq\{J,K\}\end{subarray}}\hspace{-.39in}\ep_{J',K'}\R D_{I';J',K'}.$$
If $d\!\in\!K$, the same identity holds with $|K|$ replaced by $|K|\!+\!1$ above.
The products of~$\R\wt{E}_{J,K}^0$ and~$\R\wt{E}_{J,K}^-$ with the image of the right-hand side
of the last equation under~$F$ lie in the linear span of the left-hand sides 
of~\eref{RcM04rel_e1b} and~\eref{RcM04rel_e1c} with~$\ell$ replaced by~$\ell\!+\!1$.
Along with~\eref{cM04rel_e2gen}, this implies~that
$$\R\wt{E}_{J,K}^0F\big(F_{J,K}(\cI_{\{\nod\}\sqcup[\ell]})\!\big),
\R\wt{E}_{J,K}^-F\big(F_{J,K}(\cI_{\{\nod\}\sqcup[\ell]})\!\big)\subset \cI_{0,\ell+1}\,.$$
Thus, $\Phi$ is well-defined on each $\cR_{\{\nod\}\sqcup(J\cup K)}^{\,\oplus2}$-summand.\\

\noindent
Let $(I,\{J,K\})\!\in\!\wt\cP_{\bu}(\ell)$ with $\min(J\!\cup\!K)\!\in\!J$. 
Since~$F_{I;J,K}^2$ maps the left-hand side of~\eref{cM04rel_e1a}
in its domain to a product as in~\eref{RcM04rel_e1c},
$F\!\circ\!F_{I;J,K}$ sends this product tensored with~1 to~$\cI_{0,\ell+1}$.
Let $b,c,d\!\in\!J\!\cup\!K$ be distinct.
If $d\!\in\!J$,
$$(-1)^{|K|}F_{I;J,K}^{(2)}\Big(\eref{cM04rel_e2}_{(\nod)bcd}^{\{\nod\}\sqcup(J\cup K)}\Big)
-\eref{RcM04rel_e2b2}_{dcb}^{[\ell]}
=-\!\!\!\!\!\!\!\!\!\!\sum_{\begin{subarray}{c}(I',\{J',K'\})\in\wt\cP_{\bu}(\ell)\\
d\in J',\,c\not\in I',\,b\in I'\\ (I,\{J,K\})\not\cap(I',\{J',K'\}) 
\end{subarray}}\hspace{-.51in}\ep_{J',K'}\R D_{I';J',K'}
+\!\!\!\!\!\!\!\!\!\!\!\sum_{\begin{subarray}{c}(I',\{J',K'\})\in\wt\cP_{\bu}(\ell) \\
d\in J',\,c\in I',\,b\not\in I'\\ (I,\{J,K\})\not\cap(I',\{J',K'\})
\end{subarray}}\hspace{-.51in}\ep_{J',K'}\R D_{I';J',K'}.$$
If $d\!\in\!K$, the same identity holds with $|K|$ replaced by $|K|\!+\!1$ above.
The products of~$\R\wt{D}_{I;J,K}^0$ and~$\R\wt{D}_{I;J,K}^-$ with the image of the right-hand side
of the last equation under~$F$ lie in the linear span of the left-hand sides 
of~\eref{RcM04rel_e1c} with~$\ell$ replaced by~$\ell\!+\!1$.
Along with~\eref{cM04rel_e2gen}, this implies~that
$$\R\wt{D}_{I;J,K}^{\bu}F\big(F_{I;J,K}(R_{0,\{\nod\}\sqcup I}\!\otimes_R\!\cI_{\{\nod\}\sqcup(J\cup K)})\!\big)
\subset \cI_{0,\ell+1}\,,$$
where $\bu\!=\!0,-$.\\

\noindent
With $(I,\{J,K\})$ as above, $F_{J,K}$ maps the left-hand sides of~\eref{RcM04rel_e1a},
\eref{RcM04rel_e1b}, and~\eref{RcM04rel_e1c} in its domain to analogous products in the target. 
Thus, $F\!\circ\!F_{I;J,K}$ sends these products tensored with~1
to~$\cI_{0,\ell+1}$.
Furthermore,
$$F_{J,K}\bigg(\sum_{(J',K')\in\cP(\{\nod\}\sqcup I)}\hspace{-.41in}\R E_{J',K'}\bigg)
-\sum_{(J',K')\in\cP(\ell)}\hspace{-.24in}\R E_{J',K'}
=-\sum_{\begin{subarray}{c}(J',K')\in\cP(\ell)\\
\{J,K\}\not\preceq\{J',K'\} \end{subarray}}\hspace{-.28in}\R E_{J',K'}\,.$$
If $b,c\!\in\!I$ distinct and $a\!\in\!J$,
\begin{equation*}\begin{split}
F_{J,K}\Big(\eref{RcM04rel_e2b}_{(\nod)bc}^{\{\nod\}\sqcup I}\Big)
-\eref{RcM04rel_e2b}_{abc}^{[\ell]}
=-\!\!\!\sum_{\begin{subarray}{c}(I',\{J',K'\})\in\wt\cP_{\bu}(\ell)\\
a,b\in J',\,c\in I'\\ 
(I,\{J,K\})\not\cap(I',\{J',K'\}) \end{subarray}}
\hspace{-.5in}\ep_{J',K'}\R D_{I';J',K'}
&+\!\!\! \sum_{\begin{subarray}{c} (I,\{J',K'\})\in\wt\cP_{\bu}(\ell) \\
a,c\in J',\,b\in I'\\ (I,\{J,K\})\not\cap(I',\{J',K'\})\end{subarray}}
\hspace{-.5in}\ep_{J',K'}\R D_{I';J',K'}\\
&+\!\!\! \sum_{\begin{subarray}{c} (I,\{J',K'\})\in\wt\cP_{\bu}(\ell) \\
a\in I',\,b\in J',\,c\in K'\\ (I,\{J,K\})\not\cap(I',\{J',K'\})\end{subarray}}
\hspace{-.5in}\ep_{J',K'}\R D_{I';J',K'}.
\end{split}\end{equation*}
Interchanging~$(\nod)$ and~$b$ on the left-hand side above produces the right-hand side
with~$a$ and~$b$ interchanged.
The products of~$\R\wt{D}_{I;J,K}^0$ and~$\R\wt{D}_{I;J,K}^-$ with the images of the right-hand sides
of the last two equations under~$F$ lie in the linear span of the left-hand sides 
of~\eref{RcM04rel_e1c} with~$\ell$ replaced by~$\ell\!+\!1$.
Along with~\eref{RcM04rel_e2bgen}, this implies~that
$$\R\wt{D}_{I;J,K}^{\bu}F\big(F_{I;J,K}(\cI_{0,\{\nod\}\sqcup I}\!\otimes_R\!R_{\{\nod\}\sqcup(J\cup K)})\!\big)
\subset \cI_{0,\ell+1}\,,$$
where $\bu\!=\!0,-$.
Thus, $\Phi$ is well-defined on each 
$(R_{0,\{\nod\}\sqcup I}\!\otimes_R\!R_{\{\nod\}\sqcup(J\cup K)})^{\!\oplus2}$-summand
as well.
\end{proof}

\begin{proof}[{\bf{\emph{Proof of Proposition~\ref{Ralg_prp}: surjective}}}]
The $R$-algebra $\cR_{0,\ell+1}$ is algebraically generated by the equivalence classes of 
\begin{gather}\label{Ralg_e15}
F(\R D_{I;J,K}),F(\R E_{J',K'}),\R\wt{D}_{I;J,K}^0,\R\wt{D}_{I;J,K}^-,
\R\wt{E}_{J',K'}^0,\R\wt{E}_{J',K'}^-,\wt{D}_{[\ell]-\{i^{\circ}\}}^0\in R_{0,\ell+1}\\
\notag
\hbox{with}\qquad
\big(I,\{J,K\}\big)\in\wt\cP_{\bu}(\ell),~
\{J',K'\}\in\cP(\ell),~i\!\in\![\ell],~\circ=+,-.
\end{gather}
Since the left-hand side of~\eref{RcM04rel_e2b} with $(\ell,b,c)$
replaced by $(\ell\!+\!1,i,\ell\!+\!1)$ and any $a\!\in\![\ell]$ distinct from~$i$
lies in~$\cI_{0,\ell+1}$, it is enough to take $\circ\!=\!+$ in the above list of algebraic generators.
Since the left-hand side of~\eref{RcM04rel_e2b} with $(\ell,a,b,c)$
replaced by $(\ell\!+\!1,\ell\!+\!1,1,i)$ with $i\!\neq\!1$ also
lies in~$\cI_{0,\ell+1}$, it is enough to take $i\!=\!1$ in the resulting list of algebraic generators.\\

\noindent
By~\eref{RcM04rel_e1b} and~\eref{RcM04rel_e1c} with~$\ell$ replaced by~$\ell\!+\!1$, 
$$\wt{D}_{[\ell]-\{1^+\}}^0\R\wt{D}_{I;J,K}^0,\wt{D}_{[\ell]-\{1^+\}}^0\R\wt{D}_{I;J,K}^-,
\wt{D}_{[\ell]-\{1^+\}}^0\R\wt{E}_{J',K'}^0,\wt{D}_{[\ell]-\{1^+\}}^0\R\wt{E}_{J',K'}^-\in \cI_{0,\ell+1}\,.$$
Along with~\eref{RcM04rel_e1b} and~\eref{RcM04rel_e2b}, this implies that 
$\cR_{0,\ell+1}$ is {\it linearly} generated~by 
\begin{enumerate}[label=(\arabic*),leftmargin=*]

\item\label{RgenD_it} arbitrary products of the elements of the first four types in~\eref{Ralg_e15}, 

\item\label{RgenE_it} 
these products with a {\it single} factor of~$\R\wt{E}_{J',K'}^0$ or~$\R\wt{E}_{J',K'}^-$,
but not both, and

\item arbitrary products of the elements of the first two types
with a {\it single} factor of~$\wt{D}_{[\ell]-\{1^+\}}^0$.
\end{enumerate}
Since $F$ is a ring homomorphism, the algebraic span of the elements of the first two types
in~\eref{Ralg_e15} is~$\Im\,F$.
We denote~by \hbox{$R_{0,\ell+1}^{\circ}\!\subset\!R_{0,\ell+1}$} the $(\Im\,F)$-submodule
linearly generated by~1, $\wt{D}_{[\ell]-\{1^+\}}^0$,
$\R\wt{E}_{J',K'}^0$ and~$\R\wt{E}_{J',K'}^-$
with \hbox{$\{J',K'\}\!\in\!\cP(\ell)$}, and 
$\R\wt{D}_{I;J,K}^0$ and~$\R\wt{D}_{I;J,K}^-$ with \hbox{$(I,\{J,K\})\!\in\!\wt\cP_{\bu}(\ell)$}.\\

\noindent
Let $\{J,K\}\!\in\!\cP(\ell)$ with $1\!\in\!J$ and 
$(I',\{J',K'\})\!\in\!\wt\cP_{\bu}(\ell)$ with $\min(J'\!\cup\!K')\!\in\!J'$.
By~\eref{RcM04rel_e1c},
$$\R\wt{E}_{J,K}^0\R\wt{D}_{I';J',K'}^-\in\cI_{0,\ell+1}\,.$$
If $\{J',K'\}\!\not\preceq\!\{J,K\}$, then 
$$\R\wt{E}_{J,K}^{\bu}\R\wt{D}_{I';J',K'}^{\circ}\in\cI_{0,\ell+1}
\quad\forall\,\bu,\circ=+,0,-$$
by~\eref{RcM04rel_e1b} and~\eref{RcM04rel_e1c}.
Suppose $\{J',K'\}\!\preceq\!\{J,K\}$.
If $J'\!\subset\!J$ and $K'\!\subset\!K$, then 
\begin{gather*}
\R\wt{E}_{J,K}^-\R\wt{D}_{I';J',K'}^-\!-\!
F\big(\R E_{J,K}\big)\R\wt{D}_{I';J',K'}^-
=-\R\wt{E}_{J,K}^+\R\wt{D}_{I';J',K'}^-\in\cI_{0,\ell+1},\\
\R\wt{E}_{J,K}^{\bu}\R\wt{D}_{I';J',K'}^{\circ}\in\cI_{0,\ell+1}
\quad\hbox{if}~~\bu=0,-,~\circ=\begin{cases}+,&\hbox{if}~1\!\in\!J';\\
0,&\hbox{if}~1\!\in\!I'.\end{cases}
\end{gather*}
If $1\!\in\!J'$, then $J'\!\subset\!J$, $K'\!\subset\!K$,
\begin{gather*}
\R\wt{E}_{J,K}^0\R\wt{D}_{I';J',K'}^0\!-\!
\R\wt{E}_{J,K}^0F\big(\R D_{I',J', K'}\big)
=-\R\wt{E}_{J,K}^0\big(\R\wt{D}_{I';J',K'}^+\!+\!\R\wt{D}_{I';J',K'}^-\big)
\in\cI_{0,\ell+1},\\ 
\R\wt{E}_{J,K}^-\big(\R\wt{D}_{I';J',K'}^0\!+\!\R\wt{D}_{I';J',K'}^-\big)\!-\!
\R\wt{E}_{J,K}^-F\big(\R D_{I',J', K'}\big)
=-\R\wt{E}_{J,K}^-\R\wt{D}_{I';J',K'}^+\in\cI_{0,\ell+1}\,.
\end{gather*}
If $J'\!\subset\!K$ and $K'\!\subset\!J$, then $1\!\in\!I'$, 
\begin{gather*}
\R\wt{E}_{J,K}^+\R\wt{D}_{I';J',K'}^0,
\R\wt{E}_{J,K}^0\R\wt{D}_{I';J',K'}^0,\R\wt{E}_{J,K}^-\R\wt{D}_{I';J',K'}^-\in\cI_{0,\ell+1}\,,\\
\R\wt{E}_{J,K}^-\R\wt{D}_{I';J',K'}^0\!-\!F\big(\R E_{J,K}\big)\R\wt{D}_{I';J',K'}^0
=-\R\wt{E}_{J,K}^+\R\wt{D}_{I';J',K'}^0\in\cI_{0,\ell+1}\,.
\end{gather*}
Thus, the products with one factor of $\R\wt{E}_{J,K}^0$ or~$\R\wt{E}_{J,K}^-$
and one factor of~$\R\wt{D}_{I';J',K'}^0$ or~$\R\wt{D}_{I';J',K'}^-$
lie in~$R_{0,\ell+1}^{\circ}\!+\!\cI_{0,\ell+1}$.\\

\noindent
Let $(I,\{J,K\}),(I',\{J',K'\})\!\in\!\wt\cP_{\bu}(\ell)$ be distinct
with \hbox{$\min(J\!\cup\!K)\!\in\!J$} and \hbox{$\min(J'\!\cup\!K')\!\in\!J'$}.
If \hbox{$(I,\{J,K\})\!\!\not\!\cap(I',\{J',K'\})$}, then 
\BE{Rgen_e25}\R\wt{D}_{I;J,K}^{\bu}\R\wt{D}_{I';J',K'}^{\circ}\in\cI_{0,\ell+1}
\quad\forall\,\bu,\circ=+,0,-\EE
by~\eref{RcM04rel_e1c}.
Suppose $\{J',K'\}\!\preceq\!\{J,K\}$ and and thus $I'\!\supsetneq\!I$.
If $1\!\in\!J$, the reasoning in the previous paragraph with~$\R\wt{E}_{J,K}^{\bu}$ 
replaced by~$\R\wt{D}_{I;J,K}^{\bu}$ shows that the images of 
$\R\wt{D}_{I;J,K}^{\bu}\R\wt{D}_{I';J',K'}^{\circ}$ with \hbox{$\bu,\circ\!=\!0,-$}
in~$\cR_{0,\ell+1}$
lie in the linear span of $\R\wt{D}_{I;J,K}^{\bu}$ and~$\R\wt{D}_{I';J',K'}^{\circ}$
over~$\Im\,F$. 
If $1\!\in\!I$, then $1\!\in\!I'$ and
$$\R\wt{D}_{I;J,K}^+\R\wt{D}_{I';J',K'}^0,\R\wt{D}_{I;J,K}^+\R\wt{D}_{I';J',K'}^-\in\cI_{0,\ell+1}\,.$$
If in addition $J'\!\subset\!J$ and $K'\!\subset\!K$, then
\begin{gather*}
\R\wt{D}_{I;J,K}^-\R\wt{D}_{I';J',K'}^0,
\R\wt{D}_{I;J,K}^0\R\wt{D}_{I';J',K'}^0\!-\!F\big(\R D_{I,J,K}\big)\R\wt{D}_{I';J',K'}^0
\in\cI_{0,\ell+1},\\
\R\wt{D}_{I;J,K}^0\R\wt{D}_{I';J',K'}^-,
\R\wt{D}_{I;J,K}^-\R\wt{D}_{I';J',K'}^-\!-\!F\big(\R D_{I,J,K}\big)\R\wt{D}_{I';J',K'}^-
\in\cI_{0,\ell+1}\,.
\end{gather*}
If instead $J'\!\subset\!K$ and $K'\!\subset\!J$, then
the above two lines hold with~$\R\wt{D}_{I';J',K'}^0$ and~$\R\wt{D}_{I';J',K'}^-$
interchanged.\\

\noindent
Suppose now that $J'\!\cup\!K'\!\subset\!I$ and thus $J\!\cup\!K\!\subset\!I'$.
If $1\!\in\!J$, then
$$\R\wt{D}_{I;J,K}^{\pm}\R\wt{D}_{I';J',K'}^{\bu},
\R\wt{D}_{I;J,K}^0\R\wt{D}_{I';J',K'}^{\bu}
\!-\!F\big(\R D_{I,J,K}\big)\R\wt{D}_{I';J',K'}^{\bu}
\in\cI_{0,\ell+1}$$ 
for $\bu\!=\!0,-$.
If $1\!\in\!I\!\cap\!I'$, then \eref{Rgen_e25} with \hbox{$\bu,\circ\!=\!0,-$} still holds.
Thus, the products with one factor of $\R\wt{D}_{I;J,K}^0$ or~$\R\wt{D}_{I;J,K}^-$
and one factor of~$\R\wt{D}_{I';J',K'}^0$ or~$\R\wt{D}_{I';J',K'}^-$
with distinct elements~$(I,\{J,K\})$ and~$(I',\{J',K'\})$ of~$\wt\cP_{\bu}(\ell)$
lie in~$R_{0,\ell+1}^{\circ}\!+\!\cI_{0,\ell+1}$.\\

\noindent
Let $(I,\{J,K\})\!\in\!\wt\cP_{\bu}(\ell)$ with \hbox{$\min(J\!\cup\!K)\!\in\!J$}.
If $1\!\in\!I$ and $\bu\!=\!0,-$, then
$$\R\wt{D}_{I;J,K}^{\bu}\R\wt{D}_{I;J,K}^{\bu}\!+\!\R\wt{D}_{I;J,K}^+\R\wt{D}_{I;J,K}^{\bu}
\!-\!F\big(\R D_{I;J,K}\big)\R\wt{D}_{I;J,K}^{\bu}
=-\R\wt{D}_{I;J,K}^0\R\wt{D}_{I;J,K}^-\in\cI_{0,\ell+1}\,.$$
If in addition $b\!\in\!J$, then
\begin{equation*}\begin{split}
&\ep_{J,K}\R\wt{D}_{I;J,K}^+\R\wt{D}_{I;J,K}^0-F\big(\R D_{I;J,K}\big)\wt{D}_{[\ell]-\{1^+\}}^0
-\R\wt{D}_{I;J,K}^+\eref{RcM04rel_e2b}_{(\ell+1)b1}^{[\ell+1]}\\
&\qquad+\R\wt{D}_{I;J,K}^+\bigg(\!\!\!\!\!
\sum_{\begin{subarray}{c}~~(I',\{J',K'\})\in\wt\cP_{\bu}(\ell)\\
b\in J'\neq J,\,1\in I' \end{subarray}}\hspace{-.44in}\ep_{J',K'}\R D_{I';J'\cup\{\ell+1\},K'}
-\!\!\! 
\sum_{\begin{subarray}{c}(I',\{J',K'\})\in\wt\cP_{\bu}(\ell)\\
1\in J',\,b\in I \end{subarray}}\hspace{-.38in}\R\wt{D}_{I';J',K'}^+
+\!\!\! 
\sum_{\begin{subarray}{c}(I',\{J',K'\})\in\wt\cP_{\bu}(\ell+1)\\
\ell+1\in I',\,b\in J',\,1\in K' \end{subarray}}\hspace{-.47in}\R D_{I';J',K'}\!\!\bigg)\\
&=-\big(\R\wt{D}_{I;J,K}^0\!+\!\R\wt{D}_{I;J,K}^-\big)\wt{D}_{[\ell]-\{1^+\}}^0
\!-\!\R\wt{D}_{I;J,K}^+\wt{D}_{[\ell]-\{b^+\}}^0
\in\cI_{0,\ell+1}\,.
\end{split}\end{equation*}
Each term in the first sum above is either $\R\wt{D}_{I';J',K'}^0$ or $\R D_{I';J',K'}^-$
with~$(I',\{J',K'\})$ {\it distinct} from~$(I,\{J,K\})$.
Each term in the last sum is either $\R\wt{E}_{J',K'}^0$ for some \hbox{$\{J',K'\}\!\in\!\cP(\ell)$}
or $\R\wt{D}_{I';J',K'}^0$ for some \hbox{$(I',\{J',K'\})\!\in\!\wt\cP_{\bu}(\ell)$}
{\it distinct} from~$(I,\{J,K\})$.
By~\eref{FdfnR_e} and the preceding discussion,
the last product on the left-hand side above thus lies 
in~$R_{0,\ell+1}^{\circ}\!+\!\cI_{0,\ell+1}$.
By the same reasoning with~$\eref{RcM04rel_e2b}_{(\ell+1)b1}^{[\ell+1]}$
and~$\wt{D}_{[\ell]-\{b^+\}}^0$
replaced by~$\eref{RcM04rel_e2b}_{1(\ell+1)b}^{[\ell+1]}$ and~$\wt{D}_{[\ell]-\{b^-\}}^0$,
$\R\wt{D}_{I;J,K}^-\R\wt{D}_{I;J,K}^-$ lies in~$R_{0,\ell+1}^{\circ}\!+\!\cI_{0,\ell+1}$ 
as~well.\\

\noindent
If $1\!\in\!J$, then 
$$\R\wt{D}_{I;J,K}^-\R\wt{D}_{I;J,K}^-\!+\!\R\wt{D}_{I;J,K}^0\R\wt{D}_{I;J,K}^-
\!-\!F\big(\R D_{I;J,K}\big)\R\wt{D}_{I;J,K}^-
=-\R\wt{D}_{I;J,K}^+\R\wt{D}_{I;J,K}^-\in\cI_{0,\ell+1}.$$
If in addition $b\!\in\!J\!-\!\{1\}$, then
\begin{gather*}\begin{split}
\R\wt{D}_{I;J,K}^0\R\wt{D}_{I;J,K}^0
&-\R\wt{D}_{I;J,K}^0\eref{RcM04rel_e2b}_{1b(\ell+1)}^{[\ell+1]}\\
&+\R\wt{D}_{I;J,K}^0
\bigg(\!\!\!
\sum_{\begin{subarray}{c}~(I',\{J',K'\})\in\wt\cP_{\bu}(\ell)\\
1,b\in J'\neq J \end{subarray}}\hspace{-.41in}\R\wt{D}_{I';J',K'}^0
+\sum_{\begin{subarray}{c}\{J',K'\}\in\cP(\ell)\\
1,b\in J'\end{subarray}}\hspace{-.24in}\R\wt{E}_{J',K'}^0\!\!\bigg)
\in\cI_{0,\ell+1}\,,
\end{split}\\
\begin{split}
\R\wt{D}_{I;J,K}^+\R\wt{D}_{I;J,K}^0
&-F\big(\R D_{I;J,K}\big)\wt{D}_{[\ell]-\{1^+\}}^0
-\R\wt{D}_{I;J,K}^+\eref{RcM04rel_e2b}_{1b(\ell+1)}^{[\ell+1]}\\
&+\R\wt{D}_{I;J,K}^+\bigg(\!\!
\sum_{\begin{subarray}{c}\,(I',\{J',K'\})\in\wt\cP_{\bu}(\ell)\\
1,b\in J'\subsetneq J \end{subarray}}\hspace{-.41in}\R\wt{D}_{I';J',K'}^0
\!-\!\!
\sum_{\begin{subarray}{c}(I',\{J',K'\})\in\wt\cP_{\bu}(\ell)\\
1\in J'\subsetneq J,\,b\in I' \end{subarray}}\hspace{-.38in}\R\wt{D}_{I';J',K'}^+\!\!\bigg)
\in\cI_{0,\ell+1}\,.
\end{split}\end{gather*}
If instead $b\!\in\!K$, then
\begin{gather*}\begin{split}
\R\wt{D}_{I;J,K}^0\R\wt{D}_{I;J,K}^0
&+\R\wt{D}_{I;J,K}^0\eref{RcM04rel_e2b}_{(\ell+1)1b}^{[\ell+1]}\\
&+\R\wt{D}_{I;J,K}^0\bigg(\!\!\!
\sum_{\begin{subarray}{c}\,\,(I',\{J',K'\})\in\wt\cP_{\bu}(\ell)\\
1\in J'\neq J,\,b\in K'\end{subarray}}\hspace{-.41in}\R\wt{D}_{I';J',K'}^0
+\sum_{\begin{subarray}{c}\{J',K'\}\in\cP(\ell)\\
1\in J',\,b\in K'\end{subarray}}\hspace{-.24in}\R\wt{E}_{J',K'}^0\!\!\bigg)
\in\cI_{0,\ell+1}\,,
\end{split}\\
\begin{split}
\R\wt{D}_{I;J,K}^+\R\wt{D}_{I;J,K}^0
&+F\big(\R D_{I;J,K}\big)\wt{D}_{[\ell]-\{1^+\}}^0
+\R\wt{D}_{I;J,K}^+\eref{RcM04rel_e2b}_{(\ell+1)1b}^{[\ell+1]}\\
&+\R\wt{D}_{I;J,K}^+\bigg(\!\!\!
\sum_{\begin{subarray}{c}\,\,(I',\{J',K'\})\in\wt\cP_{\bu}(\ell)\\
1\in J'\subsetneq J,\,b\in K \end{subarray}}\hspace{-.41in}\R\wt{D}_{I';J',K'}^0
\!-\!\!
\sum_{\begin{subarray}{c}(I',\{J',K'\})\in\wt\cP_{\bu}(\ell)\\
1\in J'\subsetneq J,\,b\in I' \end{subarray}}\hspace{-.38in}\R\wt{D}_{I';J',K'}^+\!\!\bigg)
\in\cI_{0,\ell+1}\,.
\end{split}\end{gather*}
By the above and~\eref{FdfnR_e}, 
$\R\wt{D}_{I;J,K}^0\R\wt{D}_{I;J,K}^0$, $\R\wt{D}_{I;J,K}^0\R\wt{D}_{I;J,K}^-$,
and $\R\wt{D}_{I;J,K}^-\R\wt{D}_{I;J,K}^-$
lie in \hbox{$R_{0,\ell+1}^{\circ}\!+\!\cI_{0,\ell+1}$} in either case.\\ 

\noindent
Thus, $\cR_{0,\ell+1}\!=\!R_{0,\ell+1}^{\circ}\!+\!\cI_{0,\ell+1}$.
Along with Lemma~\ref{Ralg_lmm} below, 
this establishes the surjectivity of the homomorphism~\eref{Ralg_e}.
\end{proof}

\begin{lmm}\label{Ralg_lmm}
Let $\ell\!\in\!\Z^+$ with $\ell\!\ge\!2$ and $\circ\!=\!0,-$.
For every \hbox{$\{J,K\}\!\in\!\cP(\ell)$},
$$\R\wt{E}_{J,K}^{\circ}(\Im\,F)\subset
\R\wt{E}_{J,K}^{\circ}\big(\Im\{F\!\circ\!F_{J,K}\}\!\big)+\cI_{0,\ell+1}\,.$$
For every \hbox{$(I,\{J,K\})\!\in\!\wt\cP_{\bu}(\ell)$},
$$\R\wt{D}_{I;J,K}^{\circ}(\Im\,F)\subset
\R\wt{D}_{I;J,K}^{\circ}\big(\Im\{F\!\circ\!F_{I;J,K}\}\!\big)+\cI_{0,\ell+1}\,.$$
\end{lmm}

\begin{proof}
Let $\{J,K\}\!\in\!\cP(\ell)$.
By~\eref{RcM04rel_e2a}, the domain of~$F$ is algebraically generated by~$\R E_{J',K'}$
with \hbox{$\{J',K'\}\!\in\!\cP(\ell)$} distinct from~$\{J,K\}$,
$\R D_{I';J',K'}$ with \hbox{$(I',\{J',K'\})\!\in\!\wt\cP_{\bu}(\ell)$},
and~$\cI_{0,\ell}$.
By~\eref{RcM04rel_e1a}, \eref{RcM04rel_e1b}, and~\eref{RcM04rel_e1c},
\begin{alignat*}{2}
\R\wt{E}_{J,K}^{\circ}F(\R E_{J',K'})&\in\cI_{0,\ell+1}
&\quad&\forall\,\{J',K'\}\!\in\!\cP(\ell),\,\{J',K'\}\!\neq\!\{J,K\}, \\
\R\wt{E}_{J,K}^{\circ}F(\R D_{I';J',K'})&\in\cI_{0,\ell+1}
&\quad&\forall\,(I',\{J',K'\})\!\in\!\wt\cP_{\bu}(\ell),\,
\{J',K'\}\!\not\preceq\!\{J,K\}\,.
\end{alignat*}
If $(I',\{J',K'\})\!\in\!\wt\cP_{\bu}(\ell)$, $J'\!\subset\!J$, and $K'\!\subset\!K$, then 
$$\R D_{I';J',K'}=
(-1)^{|K|}\ep_{J',K'}F_{J,K}(D_{\{\nod\}\sqcup I',J'\cup K'}\big)$$
by the definition of~$F_{J,K}$.
Since $F(\cI_{0,\ell})\!\subset\!\cI_{0,\ell+1}$, the first claim follows.\\

\noindent
Let \hbox{$(I,\{J,K\})\!\in\!\wt\cP_{\bu}(\ell)$}.
By~\eref{RcM04rel_e2b}, the domain of~$F$ is algebraically generated by~$\R E_{J',K'}$
with \hbox{$\{J',K'\}\!\in\!\cP(\ell)$},
$\R D_{I';J',K'}$ with \hbox{$(I',\{J',K'\})\!\in\!\wt\cP_{\bu}(\ell)$}
distinct from~$(I,\{J,K\})$, and~$\cI_{0,\ell}$.
By~\eref{RcM04rel_e1b} and~\eref{RcM04rel_e1c},
\begin{alignat*}{2}
\R\wt{D}_{I;J,K}^{\circ}F(\R E_{J',K'})&\in\cI_{0,\ell+1}
&\quad&\forall\,\{J',K'\}\!\in\!\cP(\ell),\,\{J,K\}\!\not\preceq\!\{J',K'\}, \\
\R\wt{D}_{I;J,K}^{\circ}F(\R D_{I';J',K'})&\in\cI_{0,\ell+1}
&\quad&\forall\,(I',\{J',K'\})\!\in\!\wt\cP_{\bu}(\ell),\,
\big(I',\{J',K'\}\!\big)\!\not\!\cap\big(I,\{J,K\}\!\big)\,.
\end{alignat*}
If $\{J',K'\}\!\in\!\cP(\ell)$, $J'\!\supset\!J$, and $K'\!\supset\!K$, then
$$\R E_{J',K'}=
F_{I;J,K}\big(\R E_{\{\nod\}\sqcup(I\cap J'),I\cap K'}\!\otimes\!1\big)$$
by the definition of~$F_{I;J,K}$.
For $(I',\{J',K'\})\!\in\!\wt\cP_{\bu}(\ell)$,
$$\R D_{I';J',K'}=
\begin{cases}(-1)^{|K|}\ep_{J',K'}
F_{I;J,K}(1\!\otimes\!D_{\{\nod\}\sqcup(I'-I),J'\cup K'}\!),
&\hbox{if}~J'\!\subset\!J,\,K'\!\subset\!K;\\
F_{I;J,K}(\R D_{I',\{\nod\}\sqcup(J'-J),K'-K}\!\otimes\!1),
&\hbox{if}~J'\!\supset\!J,\,K'\!\supset\!K;\\
F_{I;J,K}(\R D_{\{\nod\}\sqcup(I\cap I'),J',K'}\!\otimes\!1),
&\hbox{if}~J',K'\!\supset\!I;
\end{cases}$$
by the definition of~$F_{I;J,K}$ again.
Since $F(\cI_{0,\ell})\!\subset\!\cI_{0,\ell+1}$, the second claim follows.
\end{proof}

\vspace{.2in}

\noindent
{\it Department of Mathematics, Harvard University, Cambridge, MA 02138\\
xujiachen@g.harvard.edu}\\

\noindent
{\it  Institut de Math\'ematiques de Jussieu - Paris Rive Gauche, Paris, 75005, France\\
penka.georgieva@imj-prg.fr}\\

\noindent
{\it Department of Mathematics, Stony Brook University, Stony Brook, NY 11794\\
azinger@math.stonybrook.edu}\\

\end{document}